\documentclass[final,times,twocolumn]{elsarticle}

\usepackage{amssymb}
\usepackage{amsmath}

\usepackage{color}
\usepackage{mathtools}
\usepackage{comment}

\journal{Probabilistic Engineering Mechanics}

\begin{document}

\begin{frontmatter}



\title{Algebraic conditions for second-moment stability boundaries of linear,  time-invariant stochastic delay-differential equations} 


\author{Zsolt Iklodi \& Harry Dankowicz} 

\affiliation{organization={Department of Mechanical Engineering, University of Maryland},
            city={College Park},
            postcode={20742}, 
            state={MD},
            country={USA}}

\begin{abstract}
For linear, time-invariant stochastic delay-differential equations with a single constant delay and both multiplicative and additive noise, this paper derives optimal semi-analytic algebraic equality conditions that can be used to identify second-moment stability boundaries without the use of problem discretization. Successful validation against Monte Carlo simulations and published results for several low-dimensional models clarifies limitations of stability conditions proposed in the literature and demonstrates considerable savings in computational effort relative to discretization-based approaches. In particular, using the theory derived in this paper, second-moment stability boundaries are shown to be computable using parameter continuation techniques applied to discretization-free equality conditions that scale only with the square of the problem dimension. For the case of one-dimensional stochastic delay-differential equations, in particular, the analysis is entirely closed form with a stability condition expressed entirely in terms of elementary functions. These results are enabled by the derivation of an advection-type boundary-value problem with non-local boundary conditions for a three-variable correlation function followed by a reduction to a delay-differential boundary-value problem for a two-variable correlation function. For the former problem, observations regarding the spectral abscissa of the discretization of the corresponding infinitesimal generator, particularly that second-moment stability is lost when a real eigenvalue passes through the origin, motivate identification of second-moment stability boundaries with a loss of uniqueness of stationary solutions to the latter problem.

\end{abstract}


\begin{keyword}
Time delay \sep Stochastic dynamics \sep Correlation functions \sep Moment stability \sep Parameter continuation


\end{keyword}

\end{frontmatter}

\section{Introduction}

Delay differential equations (DDEs) are employed in a wide variety of mathematical models that aim to capture the influence of time lags or memory effects on the behavior of real-world dynamical systems. This class of equations has been shown to facilitate elegant, efficient, and insightful analytical and numerical study of complex problems, including phantom traffic jams \cite{orosz2010traffic}, machine tool chatter \cite{munoa2016chatter}, biological networks \cite{glass2021nonlinear}, and advanced control applications \cite{michiels2010control}. Recently, there has been growing interest in extending models with delay to also account for inherent noise, uncertainties, and/or unpredictable disturbances to the system dynamics. The resultant stochastic delay differential equations (SDDEs) have been employed to study traffic dynamics with stochastic delays \cite{sykora2020moment}, metal cutting with uncertain force characteristics \cite{fodor2020stochastic}, disease models with incubation delays \cite{babasola2023stochastic}, and energy harvesting applications with delayed feedback control \cite{yang2022stochastic}. 

Although SDDEs present a powerful tool for modeling complex dynamics, their widespread adoption has been slowed by excessive mathematical intricacy. Analytical solutions are rarely available \cite{frank2001stationary,sykora2019stochastic}, and most numerical methods developed for their study are computationally expensive and scale poorly with discretization resolution \cite{fodor2023collocation,rene2017mean}. Attempts to approximate the solutions of SDDEs, which have met with some success, include the derivation of moment bounds~\cite{wang2014moment}, perturbation analysis using the method of multiple scales~\cite{kuske2010competition}, and analysis of an approximate Fokker-Planck equation~\cite{ohira2000delayed}. These techniques, however, only apply to a limited subset of SDDEs and are feasible only for low-dimensional systems. For more complex delay-stochastic problems, often the most rewarding approach is \textit{brute force} forward-time Monte-Carlo simulation~\cite{rackauckas2017differentialequations,higham2001algorithmic}.

This paper is focused on a subset of SDDEs, namely linear SDDEs with constant coefficients, albeit without any limitations on dimension. While the restriction to linearity increases the chance of deriving meaningful insights, the combination of delay and stochasticity still poses significant challenges to an analytical treatment. For the target class of SDDEs, the most important objects of study are stationary distributions of stochastic trajectories near the equilibrium points of the underlying deterministic DDEs, given a well-defined source of noise or uncertainty. The existence and sensitivity of such distributions to local perturbations is of great practical interest as it corresponds to the stochastic stability of the equilibrium points (the literature on stochastic ordinary differential equations includes a variety of relevant definitions~\cite{oksendal2003stochastic,sun2006stochastic}).

The main contribution of this paper is the formulation of an algebraic condition for identifying optimal second-moment stability boundaries (we consider first-moment stability to be a solved problem), without the need for problem discretization and unimpeded by the major scaling issues characteristic of most numerical techniques. An enabling step in the derivation of this condition is the discovery of a semi-analytic, closed-form expression for a family of correlation functions of noisy trajectories near the deterministic equilibrium point, expressed in terms of eigenvalues and eigenvectors of known matrices.

The proposed formulation advances the application of SDDE models to several of the examples mentioned above by allowing reliable, optimal stability maps to be constructed (on a personal computer) in seconds compared to, at best, hours with discretization-based techniques~\cite{sykora2019stochastic,fodor2023collocation,torkamani2014numerical}. It overcomes the limitations of suboptimal guarantees of second-moment stability obtained using asymptotic~\cite{mao1992robustness,tuncc2019asymptotic} and Lyapunov-type~\cite{shaikhet2021one,samiei2013lyapunov,havskovec2022asymptotic} methods, and is complementary to theoretical stability guarantees~\cite{appleby2009geometric} that only translate to conditions on the coefficient matrices for special subclasses.

The rest of the article is structured as follows. Section~\ref{sec:mathematical fundamentals} presents the derivation of two deterministic boundary-value problems that govern the time evolution of suitably defined correlation functions which, in turn, characterize the second-moment dynamics of the most general linear, autonomous SDDE. This derivation relies on the inclusion of additional \textit{space-like} dimensions to resolve a lack of equation closure that arises in their absence. In Sec.~\ref{sec:discretization}, a possible discretization along the space-like dimensions brings to the fore several key insights from the theory of cone-preserving semi-groups on finite-dimensional spaces, most notably the observation that second-moment stability is lost only when a real eigenvalue of the corresponding infinitesimal generator crosses zero and that second-moment stability is equivalent to the existence of a positive-definite stationary covariance kernel. 

Section~\ref{sec:analytical results} relies upon these insights (without proving their generalization to the infinite-dimensional case) to characterize second-moment stability boundaries for the original SDDE without discretization. For the case of scalar SDDEs, and assuming first-moment stability, the analysis predicts an optimal closed-form inequality for second-moment stability and compares this result to edge cases documented in the literature. Section~\ref{sec:generalcase} shows how the analysis may be generalized to arbitrary number of dimensions so as to associate, without approximation, both first- and second-moment stability boundaries with the vanishing of a determinant of a matrix that can be constructed in terms of the spectral properties of suitable combinations of the coefficient matrices. 

Validations of the theoretical predictions against discretization-based approaches and stochastic integration, as well as comparisons against several incorrect results in the literature, are presented in Sec.~\ref{sec:numerical examples}. These include applications to robot control and turning processes, respectively, for which second-moment stability diagrams have practical industrial significance. Without apparent limitations, the numerical results substantiate the theoretical claims of optimality and their advantage in terms of computational efficiency relative to the use of discretization. We offer reflections and notes on possible generalizations of the presented results in Sec.~\ref{sec:conclusions}.

\section{Mathematical fundamentals}\label{sec:mathematical fundamentals}

The target of analysis in this paper is the Itô SDDE
\begin{align}
    &\mathrm{d}x(t) = \left(ax(t)+bx(t-\tau)\right) \mathrm{d}t \nonumber\\
    &\qquad+ \left(\alpha x(t)+\beta x(t-\tau)+\gamma\right) \mathrm{d}W_t\label{eq:Ito_sdde}
\end{align}
where $x:\mathbb{R}\rightarrow\mathbb{R}^n$, and $a,b,\alpha,\beta\in\mathbb{R}^{n\times n}$ and $\gamma\in\mathbb{R}^n$ are constant. This equation represents, for example, the most general form of the linearization of a non-linear autonomous SDDE with a single (constant) point delay about a deterministic equilibrium solution. When $\gamma=0$, this is a special case of the more general formulation in Appleby, Mao, and Riedle~\cite{appleby2009geometric} of a \textit{geometric Brownian motion with delay}. Without loss of generality, we omit an additive constant in the drift term (the coefficient of $\mathrm{d}t$), since this may be removed through an appropriate coordinate transformation. In Eq.~\eqref{eq:Ito_sdde}, $\mathrm{d}W_t$ represents increments of a standard Wiener process, such that $\mathbb{E}[\mathrm{d}W_t]=0$, $\mathbb{V}[\mathrm{d}W_t]=\mathrm{d}t$, $\mathbb{E}[\mathrm{d}t\,\mathrm{d}W_t]=0$, and $\mathbb{E}[\mathrm{d}W_t\mathrm{d}W_{t+\vartheta}] = \delta(\vartheta)\mathrm{d}t$. To ensure causality, we restrict attention to the case when $\tau>0$. 

Given the presence of noisy perturbations to the underlying deterministic dynamics, we are concerned with a characterization of key statistics of ensembles of trajectories, including average deviations from the trivial deterministic equilibrium at the origin and correlations among such deviations. As we discuss in this section, such characterization may be made possible by a reformulation of Eq.~\eqref{eq:Ito_sdde} in terms of an advective transport equation.

\subsection{An advection formalism}
\label{sec:advection formalism}
For $b=\beta=0$, Eq.~\eqref{eq:Ito_sdde} reduces to a linear SODE. In this case, the first and second moments, $m(t) = \mathbb{E}[x(t)]$ and $C(t) = \mathbb{E}[x(t)x^\mathsf{T}(t)]$, satisfy a cascading system of linear, autonomous ODEs~\cite{mao2007stochastic}. Algebraic conditions describing first- and second-moment stability boundaries may then be easily derived.

The case with nonzero $b$ and/or $\beta$ is considerably more difficult to analyze analytically. Although the equation for $m(t)$ remains closed (as it is identical to the underlying deterministic system), it is now a DDE. Assessment of first-moment stability, therefore, usually requires some form of numerical approximation technique~\cite{insperger2011semi,breda2014stability}. 

For nonzero $b$, analysis of second-moment stability leads to even deeper problems. Here, due to the presence of delayed terms on the right-hand side of \eqref{eq:Ito_sdde}, it is no longer possible to derive a deterministic differential equation for $C(t)$ that is expressible solely in terms of present or delayed versions of $m$ and $C$. Instead, terms containing either $\mathbb{E}[x(t)x^\mathsf{T}(x-\tau)]$ or $\mathbb{E}[x(t-\tau)x^\mathsf{T}(t)]$ lead to a lack of closure. In this paper, we resolve this lack of closure by considering, in lieu of \eqref{eq:Ito_sdde}, an advection PDE that makes explicit the state space of time histories and the corresponding evolution operator. We also retain this formulation when considering first- and second-moment dynamics. 

At the core of this approach is the introduction of one or two \textit{space-like} dimensions that parameterize the time histories of solution trajectories. As is common practice in the literature on deterministic DDEs \cite{insperger2011semi,breda2014stability,novivcenko2012phase}, we make use of the notation $x_t:\vartheta\mapsto x(t+\vartheta)$ for $\vartheta\le 0$ such that the restriction to $\vartheta \in [-\tau, \, 0]$ describes the state of the delayed dynamical system at time $t$. Eq.~\eqref{eq:Ito_sdde} may then be restated as a stochastic, first-order transport equation with boundary noise:
\begin{align}
    &\left(1+\delta(\vartheta)\right)\mathrm{d}x_t(\vartheta)\nonumber\\
    &\qquad=\left(\partial_\vartheta x_t(\vartheta)+\delta(\vartheta)(ax_t(0)+bx_t(-\tau))\right)\mathrm{d}t\nonumber\\
    &\qquad\qquad+\delta(\vartheta)(\alpha x_t(0)+\beta x_t(-\tau)+\gamma)\mathrm{d}W_t,
    \label{eq:Itotransport}
\end{align}
where $\partial_\vartheta$ is understood as a weak derivative and $\delta$ denotes the Dirac delta function.

\subsection{The mean boundary-value problem}
\label{sec:mean bvp}

We consider the assessment of first-moment stability of a linear SDDE to be a solved problem. Nevertheless, as it is an indispensable prerequisite for second-moment stability, it is worth reflecting on how first-moment dynamics may be studied using the advection formalism. To this end, let $m_t:\vartheta\mapsto\mathbb{E}[x_t(\vartheta)]$ describe the \textit{mean history function}. As long as $\vartheta+\sigma\le 0$, it follows directly from the definition that $m_{t+\sigma}(\vartheta)=m_t(\vartheta+\sigma)$. Contrary to $x_t(\vartheta)$, $m_t(\vartheta)$ is continuously differentiable in both $t$ and $\vartheta$. Indeed, from \eqref{eq:Itotransport},
\begin{align}
&\left(1+\delta(\vartheta)\right)\mathrm{d}m_t(\vartheta)\nonumber\\
&\qquad=\mathbb{E}\left[\left(1+\delta(\vartheta)\right)\mathrm{d}x_t(\vartheta)\right]\nonumber\\&\qquad=\mathbb{E}\left[\partial_\vartheta x_t(\vartheta)\right]\mathrm{d}t+\delta(\vartheta)\left(am_t(0)+bm_t(-\tau)\right)\mathrm{d}t.
\end{align}
By considering the definition of $\partial_\vartheta x_t(\vartheta)$ in terms of its action on a suitable test function, we obtain the homogeneous distributional equality
\begin{align}
    &\left(1+\delta(\vartheta)\right)\partial_t m_t(\vartheta)\nonumber\\
    &\qquad=\partial_\vartheta m_t(\vartheta)+\delta(\vartheta)\left(am_t(0)+bm_t(-\tau)\right)
\end{align}
or, equivalently, the homogeneous advection PDE
\begin{equation}
    \frac{\partial m_t(\vartheta)}{\partial t} = \frac{\partial m_t(\vartheta)}{\partial \vartheta},\,\vartheta<0\label{eq:meanbvpode}
\end{equation}
with non-local, homogeneous boundary condition
\begin{equation}
    \frac{\partial m_t(0)}{\partial t} = am_t(0)+bm_t(-\tau)\label{eq:meanbpvbc}.
\end{equation}
This is identical to the boundary-value problem obtained by considering only the noise-free limit of \eqref{eq:Itotransport}. By homogeneity, $m_t\equiv 0$ is a trivial stationary solution.

The \textit{mean boundary-value problem} in Eqs.~\eqref{eq:meanbvpode}-\eqref{eq:meanbpvbc} defines a strongly continuous semigroup on the space of differentiable functions $\psi$ on $[-\tau,0]$ with $\psi'(0)=a\psi(0)+b\psi(-\tau)$ in terms of the infinitesimal generator $\mathcal{A}:\psi\mapsto\psi'$~\cite{breda2014stability,hale2013introduction}. 

It is easy to see that a plane wave of the form $Ve^{\lambda t+\lambda_\vartheta\vartheta}$ belongs to the domain of $\mathcal{A}$ and corresponds to an eigenfunction of $\mathcal{A}$ for nontrivial $V$ provided that $\lambda_\vartheta=\lambda$ is a root of the quasi-polynomial
\begin{equation}
    \det\,(a + be^{-\lambda\tau}-\lambda I_n)\label{eq:mean_char_eq}
\end{equation}\cite{insperger2011semi,hayes1950roots}. Since these eigenfunctions form a Riesz basis for solutions to the mean boundary-value problem, first-moment stability is guaranteed by requiring that all such roots have negative real part. 

\subsection{The covariance boundary-value problem}
\label{sec:covariance boundary-value problem}
To investigate the second-moment dynamics, we define the \textit{three-variable correlation function} $C_t(\theta,\vartheta) = \mathbb{E}[x_t(\theta)x_t^\mathsf{T}(\vartheta)]$ in terms of the two extra space-like dimensions  $\theta$ and $\vartheta$, each restricted to non-positive values. By construction,
\begin{equation}
    C_t^\mathsf{T}(\theta,\vartheta) = C_t(\vartheta,\theta)\label{eq:cov_transp}
\end{equation}
and, as long as $\theta+\sigma$ and $\vartheta+\sigma$ are both non-positive,
\begin{equation}
    C_{t+\sigma}(\theta,\vartheta) = C_{t}(\theta+\sigma,\vartheta+\sigma).\label{eq:cov_shift}
\end{equation}
We recover the covariance matrix $C(t)$ by letting $\theta=\vartheta=0$. Since, by Jensen's inequality, $C_t(0,0)\succeq m_t(0)m_t^\mathsf{T}(0)$, first-moment stability is a necessary condition for second-moment stability. 

Just like $m_t$ and contrary to $x_t$, the function $C_t$ is continuously differentiable in all its arguments. Indeed, by the Itô product rule \cite{oksendal2003stochastic,doering2018modeling},
\begin{align}
    &\left(1+\delta(\vartheta)\right)\left(1+\delta(\theta)\right)\mathrm{d}C_t(\theta,\vartheta) \nonumber\\&\qquad=\left(1+\delta(\vartheta)\right)\mathbb{E}\left[\left(1+\delta(\theta)\right)\mathrm{d}x_t(\theta)\,x_t^\mathsf{T}(\vartheta)\right]\nonumber\\&\qquad\quad+\left(1+\delta(\theta)\right)\mathbb{E}\left[x_t(\theta)\,\left(1+\delta(\vartheta)\right)\mathrm{d}x_t^\mathsf{T}(\vartheta)\right]\nonumber\\&\qquad\quad+\mathbb{E}\left[\left(1+\delta(\theta)\right)\mathrm{d}x_t(\theta)\,\left(1+\delta(\vartheta)\right)\mathrm{d}x_t^\mathsf{T}(\vartheta)\right].
\label{eq:prod_rule}
\end{align}
Substitution of \eqref{eq:Itotransport} and neglecting quadratic terms in $\mathrm{d}t$ then yields the homogeneous first-order advection PDE
\begin{equation}
    \frac{\partial C_t(\theta,\vartheta)}{\partial t} = \frac{\partial C_t(\theta,\vartheta)}{\partial \theta} + \frac{\partial C_t(\theta,\vartheta)}{\partial \vartheta},\,\theta< 0,\vartheta< 0\label{eq:cov_pde}
\end{equation}
with non-local, non-homogeneous boundary conditions
\begin{equation}
    \frac{\partial C_t(0,\vartheta)}{\partial t}-\frac{\partial C_t(0,\vartheta)}{\partial \vartheta}=aC_t(0,\vartheta)+bC_t(-\tau,\vartheta)\label{eq:bcvartheta}
\end{equation}
for $\vartheta<0$,
\begin{equation}
    \frac{\partial C_t(\theta,0)}{\partial t}-\frac{\partial C_t(\theta,0)}{\partial \theta}=C_t(\theta,0)a^\mathsf{T}+C_t(\theta,-\tau)b^\mathsf{T}\label{eq:bctheta}
\end{equation}
for $\theta<0$, and
\begin{align}
    \frac{\partial C_t(0,0)}{\partial t}&=2\left[aC_t(0,0)+C_t(0,-\tau)b^\mathsf{T}+\alpha C_t(0,-\tau)\beta^\mathsf{T}\right]_s\nonumber\\&\quad+\alpha C_t(0,0)\alpha^\mathsf{T}+\beta C_t(-\tau,-\tau)\beta^\mathsf{T}\nonumber\\&\quad+2\left[\alpha m_t(0)\gamma^\mathsf{T}+\beta m_t(-\tau)\gamma^\mathsf{T}\right]_s+\gamma\gamma^\mathsf{T},\label{eq:bccorner}
\end{align}
where $\left[\cdot\right]_s$ denotes the symmetric part of the matrix argument and we have used \eqref{eq:cov_transp} to deduce that $C_t(0,0)$ is symmetric and $C_t(-\tau,0)=C_t^\mathsf{T}(0,-\tau)$. 

Analogously to the mean boundary-value problem, we anticipate (but do not prove) that the homogeneous form of the \textit{covariance boundary-value problem} in Eqs~\eqref{eq:cov_pde}-\eqref{eq:bccorner} is associated with a strongly continuous semigroup on the space of differentiable functions $\psi$ on $[-\tau,0]\times[-\tau,0]$ with $\psi^\mathsf{T}(\theta,\vartheta)=\psi(\vartheta,\theta)$,
\begin{align}
    0&=-\frac{\partial\psi}{\partial\theta}(0,\vartheta)+a\psi(0,\vartheta)+b\psi(-\tau,\vartheta), \, \vartheta<0\label{eq:edgedomain}\\
    0&=-\frac{\partial\psi}{\partial\vartheta}(\theta,0)+\psi(\theta,0)a^\mathsf{T}+\psi(\theta,-\tau)b^\mathsf{T}, \, \theta<0\\
    0&=-\frac{\partial\psi}{\partial\theta}(0,0) -\frac{\partial\psi}{\partial\vartheta}(0,0)+ 2\big[a\psi(0,0)+\psi(0,-\tau)b^\mathsf{T}\nonumber\\
    &\quad+\alpha \psi(0,-\tau)\beta^\mathsf{T}\big]_s +\alpha \psi(0,0)\alpha^\mathsf{T}+\beta \psi(-\tau,-\tau)\beta^\mathsf{T}\label{eq:cornerdomain}
\end{align}
and infinitesimal generator $\mathcal{F}:\psi(\theta,\vartheta)\mapsto\partial\psi/\partial\theta+\partial\psi/\partial\vartheta$. 

Let the symbols $\oplus$ and $\otimes$ denote the Kronecker sum and product, respectively, and define $\otimes_s$ so that $A\otimes_s B$ equals the symmetrized Kronecker product $\frac{1}{2}(A \otimes B + B \otimes A)$. Then, by substitution into \eqref{eq:edgedomain}-\eqref{eq:cornerdomain}, plane waves of the form $Ve^{\lambda t+\lambda_\theta\theta+\lambda_\vartheta\vartheta}$ can be shown to belong to the domain of $\mathcal{F}$ and correspond to an eigenfunction of $\mathcal{F}$ for nontrivial symmetric $V$ if and only if $2\lambda_\theta=2\lambda_\vartheta=\lambda$, the matrices 
\begin{equation}
    a+be^{-\lambda\tau}-\lambda I_n\label{eq:planewavecondition1}
\end{equation}
and
\begin{equation}
    \alpha\otimes\alpha+2e^{-\lambda\tau}\alpha\otimes_s \beta+e^{-2\lambda\tau}\beta\otimes\beta\label{eq:planewavecondition}
\end{equation}
are singular, and the columns of $V$ lie in the nullspace of \eqref{eq:planewavecondition1} while its vectorization lies in the nullspace of \eqref{eq:planewavecondition}.

For the case of only additive noise ($\alpha=\beta=0$), this implies that first-moment stability is a necessary and sufficient condition of second-moment stability. Consider, instead, the case with $\alpha\ne 0$ and $\beta=0$. Eq.~\eqref{eq:planewavecondition} then implies that plane wave solutions of the stipulated form belong to the domain of $\mathcal{F}$ only if $\alpha$ has at least one zero eigenvalue. This is clearly impossible in the scalar case, in contradiction to the claim in \cite{mackey1995solution}, where the simultaneous vanishing of \eqref{eq:planewavecondition1} and \eqref{eq:planewavecondition} is studied only through the implied (necessary) condition
\begin{align}
    2a+2be^{-\lambda\tau}-2\lambda+\alpha^2=0,\label{eq:Mackey}
\end{align}
which follows by consideration only of \eqref{eq:cornerdomain}. This condition, by itself, appears to predict stability provided that first-moment stability is achieved with $a$ replaced by $a+\alpha^2/2$. As we show in Sec.~\ref{sec:numerical examples}, it is straightforward to construct counterexamples to this assertion.

For the general matrix case, let $P_n$ be the unique $n^2\times n^2$ permutation matrix that maps the vectorization of an arbitrary $n\times n$ matrix onto the vectorization of its transpose. Then, by a similar study of the implied necessary condition that follows by consideration only of \eqref{eq:cornerdomain} and partially ignoring the symmetry of $V$, \cite{buckwar2013note} predicts that second-moment stability is assured provided that the roots of the quasi-polynomial
\begin{align}
    &\det\bigg(-2\lambda I_{n^2}+a\oplus a+\alpha\otimes\alpha\nonumber\\
    &\quad+e^{-\lambda\tau}\left(I_{n^2}+P_n\right)(b\otimes I_n+\beta\otimes\alpha)+e^{-2\lambda\tau}\beta\otimes\beta\bigg)\label{eq:Buckwar}
\end{align}
all have negative real part. Accounting for the symmetry of $V$ and the invariance of its vectorization under $P_n$, this simplifies to the identical condition on the roots of
\begin{align}
    &\det\bigg(-2\lambda I_{n^2}+a\oplus a+\alpha\otimes\alpha\nonumber\\
    &\qquad+e^{-\lambda\tau}(b\oplus b+2\alpha\otimes_s\beta)+e^{-2\lambda\tau}\beta\otimes\beta\bigg).\label{eq:Buckwar2}
\end{align}As we show in Sec.~\ref{sec:numerical examples}, neither condition, by itself, is sufficient to guarantee stable second-moment dynamics.

\subsection{The correlation boundary-value problem}
\label{sec:correlation bvp}

The concept of a three-variable correlation function like $C_t(\theta,\vartheta)$ is rarely used in the literature on stochastic differential equations to characterize the distribution of an ensemble of noisy trajectories. A frequently used alternative is in terms of expression of the form $\mathbb{E}[x(t)x^\mathsf{T}(t+\vartheta)]$ \cite{sun2006stochastic,milton2021mathematics}. Following this practice, we define the \textit{two-variable correlation function} $\phi_t(\vartheta) = C_t(0,\vartheta)= \mathbb{E}[x_t(0)x_t^\mathsf{T}(\vartheta)]$ for $\vartheta\le 0$. Then, by construction, 
\begin{equation}
    \phi_{t+\theta}(\vartheta-\theta)=C_{t+\theta}(0,\vartheta-\theta)=C_t(\theta,\vartheta)\label{eq:cov_reconstlower}
\end{equation}
for $\theta\ge\vartheta$ and
\begin{equation}
    \phi_{t+\vartheta}^\mathsf{T}(\theta-\vartheta)=C_{t+\vartheta}^\mathsf{T}(0,\theta-\vartheta)=C_t^\mathsf{T}(\vartheta,\theta)=C_t(\theta,\vartheta)\label{eq:cov_reconstupper}
\end{equation}
for $\vartheta\ge\theta$, whenever these are defined. It follows that the three-variable correlation function $C_t(\theta,\vartheta)$ can be fully reconstructed from time histories of $\phi_t(\vartheta)$.

\begin{figure}[ht]
\centering
\includegraphics[width=0.85\columnwidth]{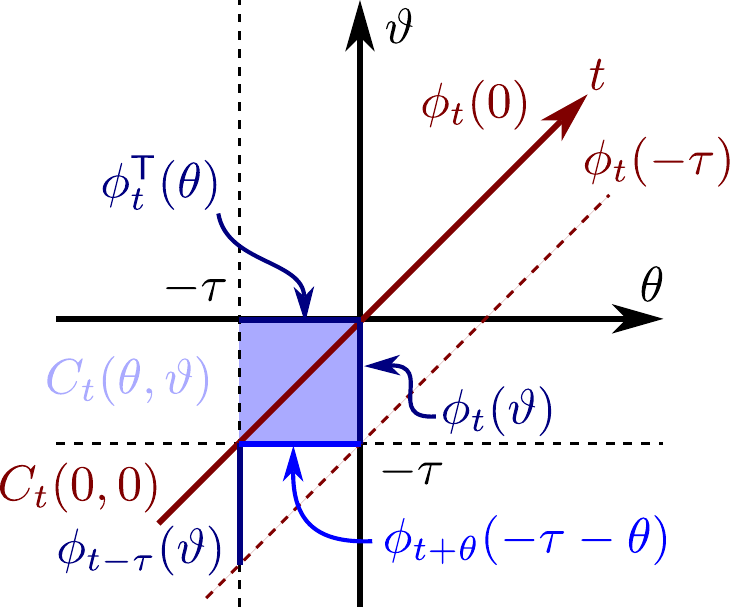}
\caption{Domains of definition for the solutions to the covariance  and correlation boundary-value problems. Here, the fundamental domain $[-\tau,0]\times[-\tau,0]$ is shaded. The functional expressions $\phi_t(\vartheta)$, $\phi_{t-\tau}(\vartheta)$, $\phi_{t+\theta}(-\tau-\theta)$, and $\phi_t^\mathsf{T}(\theta)$ are associated with values of $C_t(\theta,\vartheta)$ along the corresponding line segments. The 45$^\circ$ arrow represents the advective property $C_{t+\sigma}(\theta,\vartheta)=C_t(\theta+\sigma,\vartheta+\sigma)$.}\label{fig:pde_domain}
\end{figure}

Figure \ref{fig:pde_domain} provides further insights into the relationship between $C_t(\theta,\vartheta)$ and $\phi_t(\vartheta)$ by visualizing the domains of these two mathematical objects on the $\theta$ - $\vartheta$ plane. Here, the shaded area defined by $(\theta,\vartheta)\in[-\tau, \, 0]\times[-\tau,0]$ is the \textit{fundamental domain} for the covariance boundary-value problem. 
By Eq.~\eqref{eq:cov_transp}, three-variable correlation functions are symmetric under transposition and the variable interchange $\theta\leftrightarrow\vartheta$. Moreover, by Eq.~\eqref{eq:cov_shift} it follows that evolution in $t$ is equivalent to movement along lines that make an angle of $45^\circ$ with each of the positive $\theta$ and $\vartheta$ axes. Consequently, by Eq.~\eqref{eq:cov_reconstlower}, $\phi_{t+\theta}(\vartheta-\theta)$ encodes $C_t(\theta,\vartheta)$ on the lower triangle half of the fundamental domain, where  $\theta > \vartheta$, while by Eq.~\eqref{eq:cov_reconstupper}, $\phi_{t+\vartheta}^\mathsf{T}(\theta-\vartheta)$ encodes $C_t(\theta,\vartheta)$ on  the upper triangle half of the fundamental domain, where $\theta < \vartheta$.

From Eqs.~\eqref{eq:cov_reconstlower} and \eqref{eq:cov_reconstupper}, it follows that 
\begin{equation}
    C_t(-\tau,\vartheta)=\left\{\begin{array}{ll}\phi_{t-\tau}(\vartheta+\tau) & \vartheta\le-\tau\\\phi_{t+\vartheta}^\mathsf{T}(-\tau-\vartheta) & \vartheta\ge-\tau\end{array}\right.
\end{equation}
Consequently, by substitution of the definition of $\phi_t(\vartheta)$ into \eqref{eq:bcvartheta}, we obtain
\begin{align}
    \frac{\partial \phi_t(\vartheta)}{\partial t} -\frac{\partial \phi_t(\vartheta)}{\partial \vartheta} = a\phi_t(\vartheta) + b\phi^\mathsf{T}_{t+\vartheta}(-\tau-\vartheta)\label{eq:corr_pde}
\end{align}
for $\vartheta\in[-\tau,0)$ and
\begin{align}
    \frac{\partial\phi_t(\vartheta)}{\partial t}-\frac{\partial\phi_t(\vartheta)}{\partial\vartheta}=a\phi_t(\vartheta)+b\phi_{t-\tau}(\vartheta+\tau)\label{eq:corr_pde2}
\end{align}
for $\vartheta\in(-\infty,-\tau)$, while Eq.~\eqref{eq:bccorner} yields the non-local boundary condition
\begin{align}
    \frac{\partial \phi_t(0)}{\partial t}&=2\left[a\phi_t(0)+\phi_t(-\tau)b^\mathsf{T}+\alpha \phi_t(-\tau)\beta^\mathsf{T}\right]_s\nonumber\\&\qquad+\beta \phi_{t-\tau}(0)\beta^\mathsf{T}+\alpha \phi_t(0)\alpha^\mathsf{T}\nonumber\\&\qquad+2\left[\alpha m_t(0)\gamma^\mathsf{T}+\beta m_t(-\tau)\gamma^\mathsf{T}\right]_s+\gamma\gamma^\mathsf{T},\label{eq:corr_pde_bc}
\end{align}
where, by construction, $\phi_t(0)$ is symmetric. In contrast to the covariance boundary problem, this \textit{correlation boundary-value problem} includes delay terms in the PDE as well as the boundary conditions. At first glance it might not seem like the benefits of reducing the problem dimension by one outweigh the burdens of these additional complexities. Nonetheless, in Sec.~\ref{sec:analytical results} below, we show that this form is more amenable to analytical treatment, especially when searching for stationary solutions.

\subsection{Discretization}
\label{sec:discretization}
The primary contribution of this paper lies in the derivation of algebraic conditions for second-moment stability boundaries directly from the correlation boundary-value problem. Nevertheless, we consider in this section a natural discretization of the covariance boundary-value problem in the space-like variables, as properties of the spatially discretized dynamics inform our subsequent analysis. (A similar discretization of the correlation boundary-value problem \eqref{eq:corr_pde}-\eqref{eq:corr_pde_bc} is performed in \ref{app:num_corr} and reflected upon in Sec.~\ref{sec:numerical examples}.)

Recall from Sec.~\ref{sec:covariance boundary-value problem} the covariance boundary-value problem as the functional differential equation
\begin{align}
    &\left(1+\delta(\theta)\right)\left(1+\delta(\vartheta)\right)\left(\frac{\partial C_t}{\partial t}-\mathcal{F}C_t\right)=\delta(\theta)\delta(\vartheta)G_t\label{eq:cov_op_eq}
\end{align}
on the domain of the infinitesimal generator $\mathcal{F}$ with nonhomogeneity
\begin{equation}
     G_t =2\left[\alpha m_t(0)\gamma^\mathsf{T}+\beta m_t(-\tau)\gamma^\mathsf{T}\right]_s +\gamma\gamma^\mathsf{T}.
\end{equation}
Inspired by \cite{breda2014stability}, consider a mesh $-\tau=s_1<\cdots<s_M=0$ for some integer $M$ and let bold-faced symbols represent matrices whose size depends on $M$. For example, let $\mathbf{e}_i$ denote the $M$-dimensional column matrix with $\left(\mathbf{e}_i\right)_j$ equal to the Kronecker delta tensor $\delta_{ij}$. Define $\mathbf{c}(t)$ as the vectorization of the block matrix
\begin{equation}
    \mathbf{C}_t=\begin{bmatrix}
        C_t(s_1,s_1) & \cdots & C_t(s_1,s_M)\\\vdots & \ddots & \vdots\\C_t(s_M,s_1) & \cdots & C_t(s_M,s_M)
    \end{bmatrix}\label{eq:vec_cov}.
\end{equation}
Then, we refer to the equation
\begin{equation}
    \dot{\mathbf{c}}(t) = \mathbf{F}\mathbf{c}(t) + \mathbf{g}(t)\label{eq:cov_disc}
\end{equation}
for suitably constructed matrices $\mathbf{F}$ and $\mathbf{g}(t)$ as the \textit{surrogate covariance ODE} associated with the covariance boundary-value problem. Indeed, by linearity, it is clear that $\mathbf{F}$ may be decomposed into a sum of two terms depending on the drift and diffusion coefficients, respectively.

For example, if $\mathbf{m}(t)$ denotes the vectorization of the block matrix
\begin{equation}
    \begin{bmatrix}
        m_t(s_1) & \cdots & m_t(s_M)
    \end{bmatrix},
\end{equation}
then $\mathbf{g}(t)$ is the vectorization of
\begin{equation}
\begin{split}
    \left(\mathbf{e}_{M}\mathbf{e}_{M}^\mathsf{T}\right)\otimes\left(\boldsymbol{\eta}^\mathsf{T}\mathbf{\mathbf{m}}(t)\gamma^\mathsf{T} + \gamma \mathbf{m}^\mathsf{T}(t)\boldsymbol{\eta} + \gamma\gamma^\mathsf{T}\right)
    \end{split},\label{eq:cov_disc_const}
\end{equation}
where
\begin{equation}
    \boldsymbol{\eta}^\mathsf{T} = \begin{bmatrix} \beta & \mathbf{0}_{n\times n(M-2)} & \alpha\end{bmatrix}
\end{equation}
and $\mathbf{0}_{k\times l}$ denotes a $k$ by $l$ null matrix when at least one dimension scales with $M$. Similarly, since
\begin{align}
    \boldsymbol{\eta}^\mathsf{T}\mathbf{C}_t\boldsymbol{\eta}&=\alpha C_t(0,0)\alpha^\mathsf{T}+\beta C_t(-\tau.-\tau)\beta^\mathsf{T}\nonumber\\
    &\qquad+\alpha C_t(0,-\tau)\beta^\mathsf{T}+\beta C_t(-\tau,0)\alpha^\mathsf{T},
\end{align}
the contribution to $\mathbf{F}$ from the presence of noise is given by the Kronecker product $\mathbf{B}\otimes\mathbf{B}$, where
\begin{equation}
    \mathbf{B} =\begin{bmatrix} \hat{\mathbf{0}}_{M\times M} \otimes I_n \\ \boldsymbol{\eta}^\mathsf{T} \end{bmatrix}.\label{eq:discB}
\end{equation}
Finally, let $\mathbf{D}_M \in\mathbb{R}^{M\times M}$ denote an appropriately defined approximate differentiation matrix (e.g., the one derived in Chapter 6 of \cite{trefethen2000spectral} for a Chebyshev maxima grid, as seen in \cite{torkamani2014numerical}) and let an accent $\hat{\;\;}$ indicate that the last row has been removed from a matrix. Then, the matrices $ \left(\left(\mathbf{D}_{M}\otimes I_{n}\right) \otimes \mathbf{I}_{Mn}\right)\mathbf{c}(t)$ and $\left(\mathbf{I}_{Mn} \otimes \left(\mathbf{D}_{M}\otimes I_{n}\right)\right)\mathbf{c}(t)$
are discretized approximations of the  vectorizations of the block matrices 
\begin{equation}
    \begin{bmatrix}
        \partial_\theta C_t(s_1,s_1) & \cdots & \partial_\theta C_t(s_1,s_M)\\\vdots & \ddots & \vdots\\\partial_\theta C_t(s_M,s_1) & \cdots & \partial_\theta C_t(s_M,s_M)
    \end{bmatrix}
\end{equation}
and
\begin{equation}
    \begin{bmatrix}
        \partial_\vartheta C_t(s_1,s_1) & \cdots & \partial_\vartheta C_t(s_1,s_M)\\\vdots & \ddots & \vdots\\\partial_\vartheta C_t(s_M,s_1) & \cdots & \partial_\vartheta C_t(s_M,s_M)
    \end{bmatrix},
\end{equation}
respectively. Accounting also for the boundary conditions \eqref{eq:bcvartheta} and \eqref{eq:bctheta}, the noise-independent contribution to $\mathbf{F}$ is, consequently, given by the Kronecker sum $\mathbf{A}\oplus \mathbf{A}$, where
\begin{equation}
    \mathbf{A} = \begin{bmatrix}\hat{\mathbf{D}}_M \otimes I_n \\ \begin{bmatrix} b & \mathbf{0}_{n\times n(M-2)} & a\end{bmatrix}\end{bmatrix}.\label{eq:discA}
\end{equation}

As an alternative to the ``discretization last'' approach of deriving \eqref{eq:cov_disc} from \eqref{eq:cov_pde}-\eqref{eq:bccorner}, one may pursue a ``discretization first'' approach, in which the original Itô SDDE \eqref{eq:Ito_sdde} is replaced by the surrogate high-dimensional system of SODEs
\begin{equation}
    \mathrm{d}\mathbf{x}(t) = \mathbf{A}\mathbf{x}(t) \,\mathrm{d}t + \left(\mathbf{B}\mathbf{x}(t) + \mathbf{\Gamma}\right) \mathrm{d}W_t.\label{eq:surr_SODE}
\end{equation}
where $\mathbf{x}(t)$ is the vectorization of 
\begin{equation}
    \begin{bmatrix}x_t(s_1) & \cdots & x_t(s_M)\end{bmatrix}
\end{equation}
and $\mathbf{\Gamma}^\mathsf{T} = \begin{bmatrix}\mathbf{0}_{1\times n(M-1)} & \gamma^\mathsf{T}\end{bmatrix}$. It is straightforward to show that the vectorized form of the associated covariance ODE is again the surrogate covariance ODE \eqref{eq:cov_disc}. Examples in the literature of using such a ``discretization first'' approach to assess second-moment stability of linear SDDEs include a finite difference scheme in \cite{rene2017mean} and a pseudo-spectral approximation technique in \cite{torkamani2014numerical}.

By the properties of linear systems of SODEs, it follows that first-moment stability of the discretized dynamics in \eqref{eq:surr_SODE} is guaranteed by $\mathbf{A}$ being Hurwitz. To investigate second-moment stability, it therefore suffices to consider the case when $\mathbf{m}=\mathbf{0}_{1,nM}$, such that $\mathbf{g}(t)\equiv\bar{\mathbf{g}}$ is the vectorization of $\mathbf{e}_M\mathbf{e}_M^\mathsf{T}\otimes\gamma\gamma^\mathsf{T}$. In this case, $\mathbf{B}=\mathbf{0}_{M\times M}$ implies that second-moment stability is guaranteed by $\mathbf{A}\oplus\mathbf{A}$ being Hurwitz. Since the eigenvalues of $\mathbf{A}\oplus\mathbf{A}$ equal all possible sums of pairs of eigenvalues of $\mathbf{A}$, we again conclude that first-moment stability implies second-moment stability in the presence of only additive noise.

For $\mathbf{B}\ne\mathbf{0}_{M\times M}$, second-moment stability is determined by the eigenvalues of the coefficient matrix $\mathbf{F} = \mathbf{A} \oplus \mathbf{A} + \mathbf{\textbf{B}} \otimes \mathbf{\textbf{B}}$, which is Hurwitz only if $\mathbf{A}$ is Hurwitz. By the invariance of positive semi-definite symmetric matrices under the semigroup $e^{t\mathbf{F}}$ and reference to the finite-dimensional Krein-Rutman theorem~\cite{Schneider65}, it follows that the spectral abscissa of $\mathbf{F}$ is an eigenvalue of $\mathbf{F}$. We conclude that second-moment stability is lost along a curve in parameter space only when the corresponding eigenvalue crosses through $0$ from the negative to the positive half plane, through a point where $\det\,(\mathbf{F})=0$. In other words, second-moment stability cannot be lost through a simple Hopf bifurcation. From \eqref{eq:cov_disc}, such a point coincides with a loss of uniqueness of the equilibrium solution $\bar{\mathbf{c}}$ for which $\mathbf{F}\bar{\mathbf{c}}=-\bar{\mathbf{g}}$. As such a point is approached for $\bar{\mathbf{g}}\ne\mathbf{0}$, the equilibrium solution $\bar{\mathbf{c}}$ must grow in norm beyond all bounds. Finally, we note that second-moment stability implies convergence of the second-moment dynamics to $\bar{\mathbf{c}}$, which must therefore be the vectorization of a positive semi-definite matrix. In fact, as shown by \cite{nandanoori2018mean}, second-moment stability of \eqref{eq:surr_SODE} with $\mathbf{\Gamma}=\mathbf{0}$ is equivalent to $\bar{\mathbf{c}}$ being positive definite for some $\mathbf{\Gamma}\ne\mathbf{0}$.

We proceed in the next section by anticipating, without proof, that these properties relating to second-moment stability translate also to the infinite-dimensional covariance/correlation boundary-value problems. Specifically, we seek to identify the locus of points within the region of first-moment stability where i) uniqueness of a stationary solution to either problem cannot be guaranteed and ii) a change occurs from positive to negative definiteness of the corresponding kernel.

\section{Analytical results}\label{sec:analytical results}

As noted in the previous section, the correlation boundary-value problem \eqref{eq:corr_pde}-\eqref{eq:corr_pde_bc} includes delay terms in the PDE and the boundary conditions. In spite of this added complexity, we show in this section that closed-form analysis of this delay-PDE enables the derivation of candidate algebraic equalities and inequalities defining second-moment stability boundaries and the region of second-moment stability.

To make our notation more concise, and to allow the obtained formulas to be reducible to a simple form, we restrict attention at first to scalar problems where $n = 1$, for which transposition is the identity. After a detailed comparison of our findings with analytic second-moment stability analysis results available in the literature, we show how our analytical results generalize to non-scalar problems.

\subsection{The scalar problem}
\label{sec:scalar problem}

We recall from \eqref{eq:bccorner} and \eqref{eq:corr_pde_bc} the dependence of the second-moment dynamics on the evolution of the first moment $m_t$. Therefore, before investigating the second-moment stability of the scalar form of \eqref{eq:Ito_sdde}, we cite the following well-known results for first-moment stability of the Hayes equation~\cite{hayes1950roots} $\dot{x}(t)=ax(t)+bx(t-\tau)$, in terms of the roots of its characteristic quasi-polynomial $p(\lambda)=a+be^{-\lambda\tau}-\lambda$.

\begin{figure}[t]
    \centering
    \includegraphics[width=0.95\columnwidth]{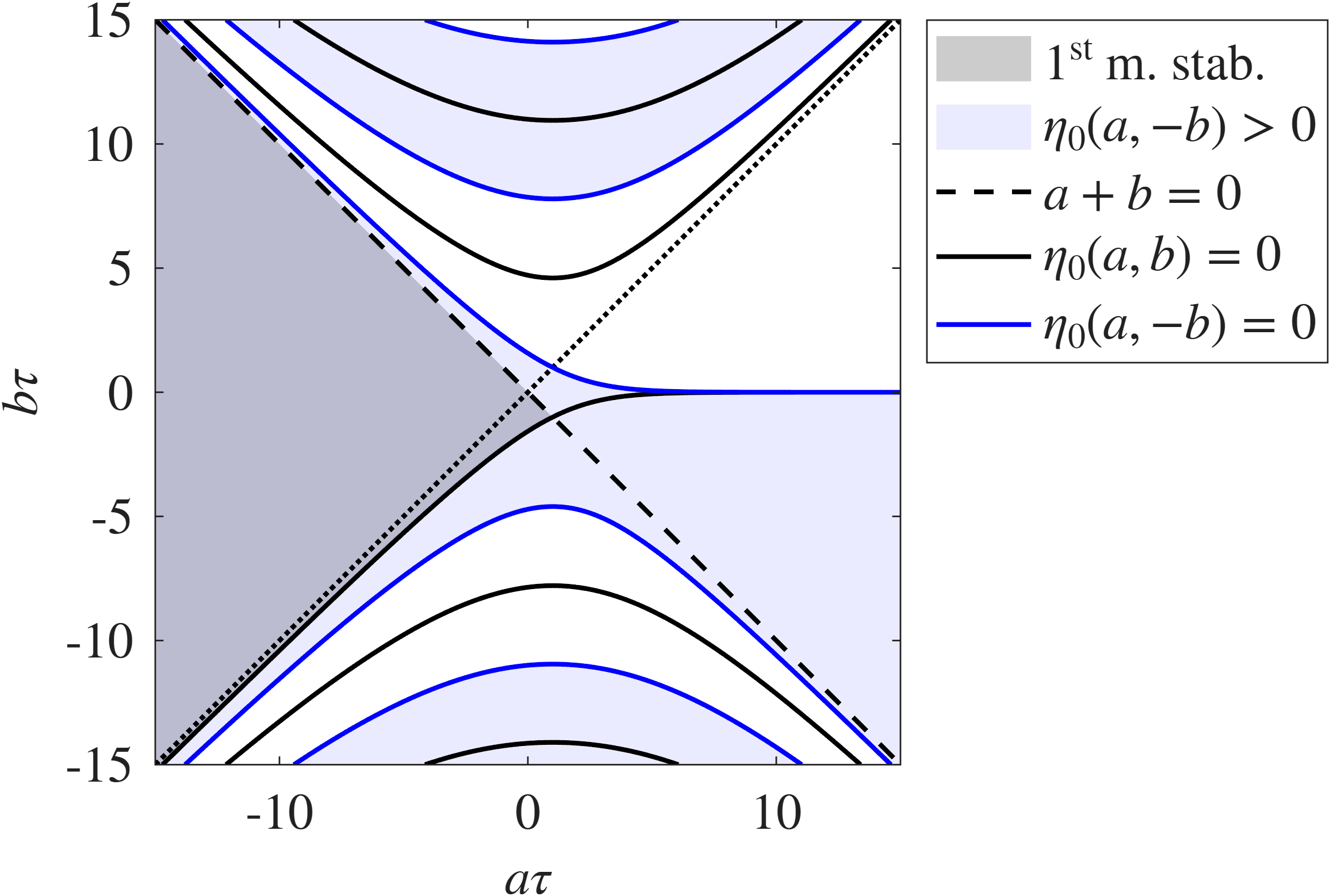}
    \caption{Analytic first-moment stability bounds of a scalar SDDE of the form \eqref{eq:Ito_sdde}. Here, first-moment stability is guaranteed by the simultaneous conditions $a+b<0$, $\eta_{\tau/2}(a,b)>0$, corresponding to the dark shaded region. The zero-level curves of $\eta_{\tau/2}(a,b)$ consist of Hopf bifurcations and neutral saddles. The zero-level curves of $\eta_{\tau/2}(a,-b)$ and regions where this quantity is positive play a role in the derivation of a condition for second-moment stability in Sec.~\ref{sec:stationary solutions}. In particular, we see that first-moment stability implies that $\eta_{\tau/2}(a,-b)>0$.}
    \label{fig:scalar_astab}
\end{figure}

First, we note that saddle-node bifurcations coincide with parameter pairs $(a,b)$ such that $\lambda=0$ is a root of $p(\lambda)$. Clearly, these occur along the line $a+b=0$. Second, Hopf bifurcations coincide with parameter pairs $(a,b)$ for which $\lambda=\pm\mathrm{i}\omega$ are a pair of conjugate imaginary roots of $p(\lambda)$ with $\omega> 0$. These occur along the family of non-intersecting and alternately nested (consider $\omega\tau=m\pi+\pi/2$ for integer values $m$) curves $\omega\mapsto(a,b)=(\omega\cot\omega\tau,-\omega\csc\omega\tau)$ for $\omega\tau\ne n\pi$ for nonzero $n\in\mathbb{Z}$, and coincide with points along the zero-level set of
\begin{equation}
    \eta_{\tau/2}(a,b)=\cos \frac{\tau\sqrt{b^2-a^2}}{2}+\frac{b-a}{\sqrt{b^2-a^2}} \sin \frac{\tau\sqrt{b^2-a^2}}{2}\label{eq:1stmomenthopf}
\end{equation}
where $b^2-a^2=\omega^2$. As also confirmed by the illustration in Fig.~\ref{fig:scalar_astab}, a unique intersection is obtained between the saddle-node bifurcation curve and the family of Hopf bifurcation curves at $a=-b=1/\tau$. Moreover, for $b=0$, stability is obtained for $a<0$. We conclude that first-moment stability is, in general, obtained for $a=\omega\cot\omega\tau<1/\tau$ and $b\in(-\omega\csc\omega\tau,-\omega\cot\omega\tau)$ for $\omega\tau\in(0,\pi)$, shaded in gray in Fig.~\ref{fig:scalar_astab}. In this region, $a+b<0$ and $\eta_{\tau/2}(a,b)>0$.

Recall again, by Jensen's inequality, the positive semi-definiteness of $\phi_t(0)-m_t(0)m_t^\mathsf{T}(0)$. Restricting attention, therefore, to the case of first-moment stability, we proceed to assume that enough time has passed since the last external perturbation such that $m_t(\vartheta)\approx 0$. Substitution into the scalar form of the delay-PDE in \eqref{eq:corr_pde}-\eqref{eq:corr_pde2} and its non-local boundary condition \eqref{eq:corr_pde_bc} then yields the following simplified version of the correlation boundary-value problem:
\begin{align}
    &\frac{\partial \phi_t(\vartheta)}{\partial t} -\frac{\partial \phi_t(\vartheta)}{\partial \vartheta} = a\phi_t(\vartheta) + b\phi_{t+\vartheta}(-\tau-\vartheta),\label{eq:corr_pde_scalar1}
\end{align}
for $\vartheta\in[-\tau,0)$,
\begin{align}
    &\frac{\partial\phi_t(\vartheta)}{\partial t}-\frac{\partial\phi_t(\vartheta)}{\partial\vartheta}=a\phi_t(\vartheta)+b\phi_{t-\tau}(\vartheta+\tau),\label{eq:corr_pde_scalar2}
\end{align}
for $\vartheta\in(-\infty,-\tau)$, and
\begin{align}
    \frac{\partial \phi_t(0)}{\partial t}&=(2a+\alpha^2)\phi_t(0)+(2b+2\alpha\beta)\phi_t(-\tau)\nonumber\\
    &\qquad+\beta^2 \phi_{t-\tau}(0)+\gamma^2.\label{eq:corr_pde_scalar3}
\end{align}
As a central result of this paper we show that there exists a subset of the region $\{(a,b)\,|\,a+b<0,\,\eta_{\tau/2}(a,b)>0\}$ of first-moment stability where this intricate boundary-value problem has a unique stationary solution that is available in closed form, such that parameter sets where uniqueness is not guaranteed coincide with second-moment stability boundaries.

\subsection{Stationary solutions}\label{sec:stationary solutions}

It follows from Eqs.~\eqref{eq:cov_reconstlower} and \eqref{eq:cov_reconstupper} that $\bar{C}(\theta,\vartheta)=\bar{\phi}\left(-|\theta-\vartheta|\right)$ is a stationary solution of the covariance boundary-value problem if $\bar{\phi}(\vartheta)$ is a stationary solution of the correlation boundary-value problem. In particular, $\bar{\phi}\left(-|\theta-\vartheta|\right)$ must be a valid covariance kernel on $[-\tau,0]\times[-\tau,0]$, i.e., every matrix constructed from this kernel must be positive (semi-)definite. As a special case, we must have $\bar{\phi}(0)\ge 0$. We refer to a solution $\bar{\phi}$ that violates this condition as \textit{spurious}. 

Assuming the existence of a stationary solution $\bar{\phi}$ of \eqref{eq:corr_pde_scalar1}-\eqref{eq:corr_pde_scalar3}, let $\hat{\phi}(\sigma)=\bar{\phi}(-\sigma-\tau)$ for $\sigma\in(0,\infty)$. Then, Eq.~\eqref{eq:corr_pde_scalar2} implies that
\begin{equation}
    \hat{\phi}'(\sigma)=a\hat{\phi}(\sigma)+b\hat{\phi}(\sigma-\tau),\label{eq:corr_stat_dde}
\end{equation}
where $'$ denotes differentiation with respect to $\sigma$. Since this is identical in form to the deterministic part of the original SDDE \eqref{eq:Ito_sdde} and since we already assume first-moment stability, it follows that
\begin{equation}
    \lim_{\vartheta \rightarrow  -\infty}\bar{\phi}(\vartheta) = \lim_{\vartheta \rightarrow  -\infty}\hat{\phi}(-\vartheta-\tau)=0.
\end{equation}

It remains to consider the boundary-value problem for $\bar{\phi}(\vartheta)$ in $\vartheta\in[-\tau,0]$. Substitution in \eqref{eq:corr_pde_scalar1} and \eqref{eq:corr_pde_scalar3} yields the scalar ODE
\begin{equation}
    \bar{\phi}'(\vartheta) = -a\bar{\phi}(\vartheta)-b\bar{\phi}(-\vartheta-\tau),
    \label{eq:corr_scalar_ode1}
\end{equation}
where $'$ now denotes differentiation with respect to $\vartheta\in[-\tau,0)$, and the boundary condition
\begin{equation}
    (2a + \alpha^2 + \beta^2)\bar{\phi}(0)+(2b + 2\alpha\beta)\bar{\phi}(-\tau)+\gamma^2 =0. \label{eq:corr_scalar_const_bc}
\end{equation}
Differentiation of \eqref{eq:corr_scalar_ode1} with respect to $\vartheta$ and simplification using \eqref{eq:corr_scalar_ode1} then yields the equivalent boundary-value problem
\begin{equation}
    \bar{\phi}''(\vartheta)=(a^2-b^2)\bar{\phi}(\vartheta),\,\bar{\phi}'(0)=-a\bar{\phi}(0)-b\bar{\phi}(-\tau),\label{eq:corr_2nd_ode}
\end{equation}
and \eqref{eq:corr_scalar_const_bc}. This has the unique solution
\begin{equation}
\begin{split}
        &\bar{\phi}(\vartheta) =
-\frac{\gamma^2}{\chi} \eta_{\vartheta+\tau/2}(a,-b),
    \label{eq:corr_scalar_solution}
\end{split}
\end{equation}
where (consistent with \eqref{eq:1stmomenthopf})
\begin{align}
\eta_{\xi}(a,-b)&=\cosh \mu \xi-\frac{a+b}{\mu} \sinh \mu \xi
\end{align}
and $\mu=\sqrt{a^2-b^2}$,
provided that
\begin{align}
    \chi &= \left((\alpha +\beta )^2+2 a+2 b\right) \cosh \frac{\mu\tau}{2}\nonumber\\
    &\qquad-\frac{a+b}{\mu} \left((\alpha -\beta )^2+2 a-2 b\right) \sinh \frac{\mu\tau}{2}
\label{eq:corr_chi}
\end{align}
is nonzero. When $\chi=0$, a solution exists only if $\gamma=0$, in which case $\bar{\phi}(\vartheta)$ is an arbitrary multiple of $\eta_{\vartheta+\tau/2}(a,-b)$. The solution in \eqref{eq:corr_scalar_solution} is spurious unless $\eta_{\tau/2}(a,-b)/\chi\le 0$. As demonstrated by the shaded blue regions in Fig.~\ref{fig:scalar_astab} $a+b<0$ and $\eta_{\tau/2}(a,b)>0$ together imply that $\eta_{\tau/2}(a,-b)>0$. First-moment stability thus implies the existence of a unique, non-spurious stationary solution of the correlation boundary-value problem if and only if $\chi<0$. By the analogy with the discrete case, we anticipate that this is the region of second-moment stability.

\begin{figure*}[ht]
    \centering
    \includegraphics[width=0.33\linewidth]{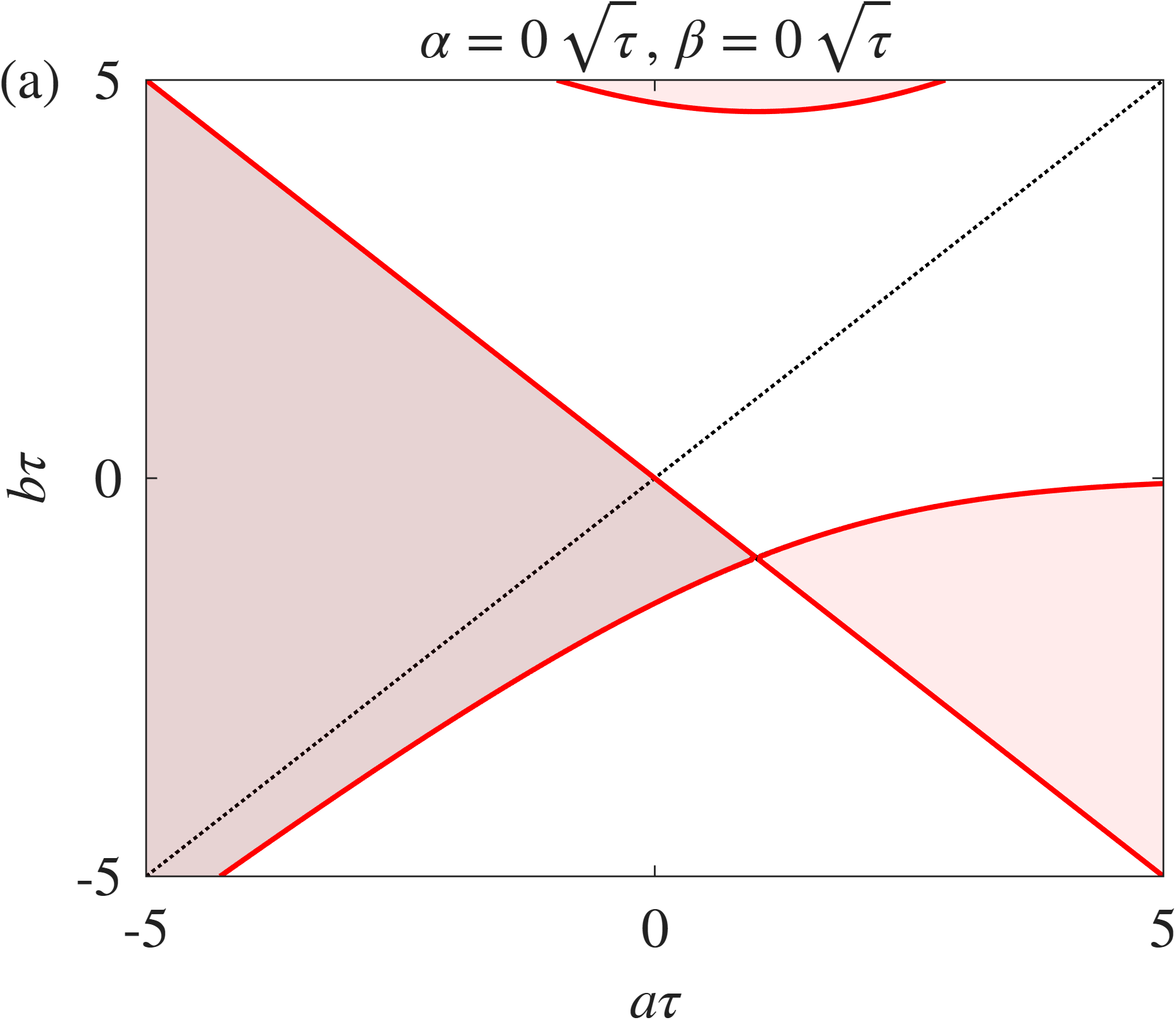}
    \includegraphics[width=0.33\linewidth]{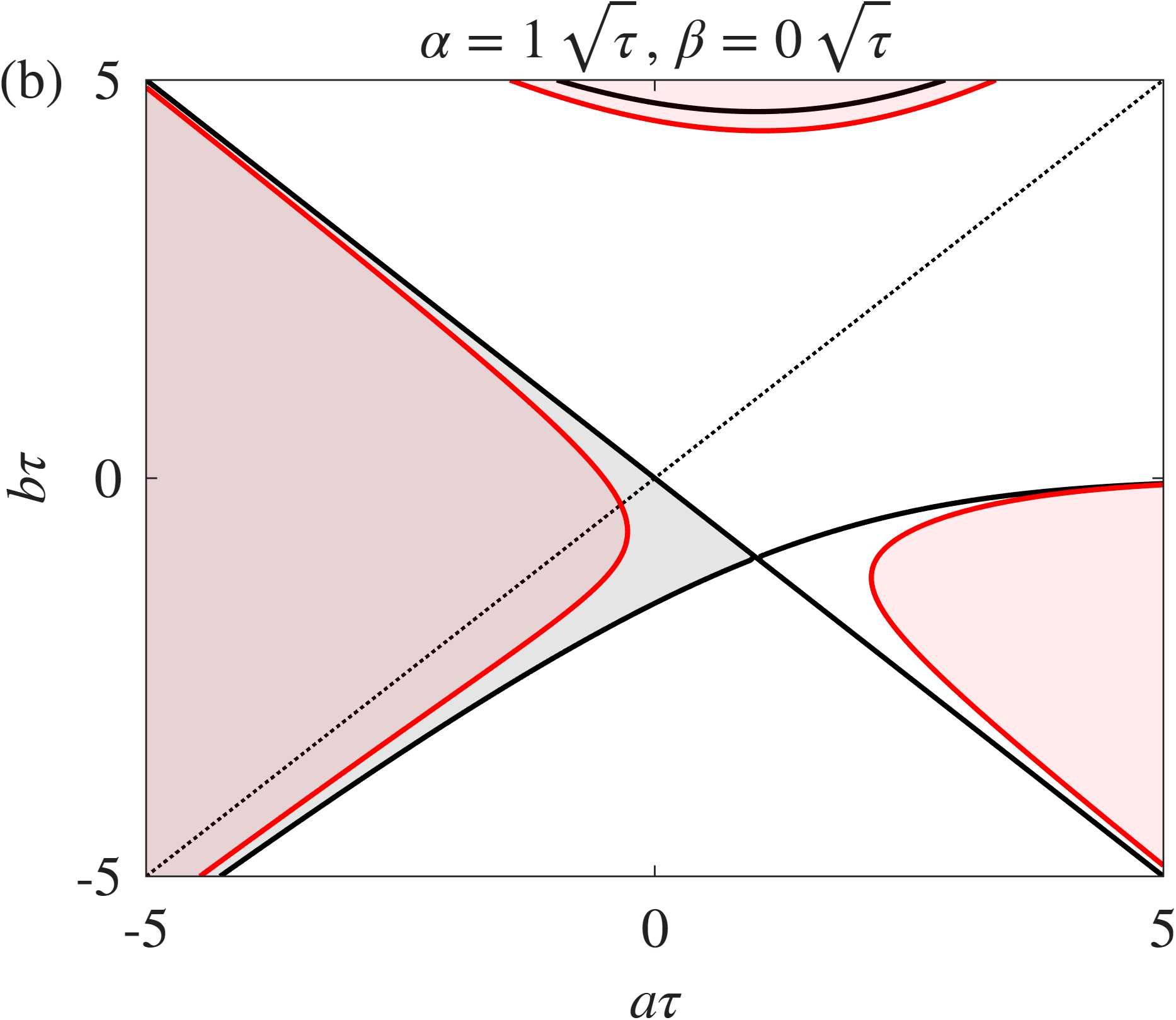}
    \includegraphics[width=0.33\linewidth]{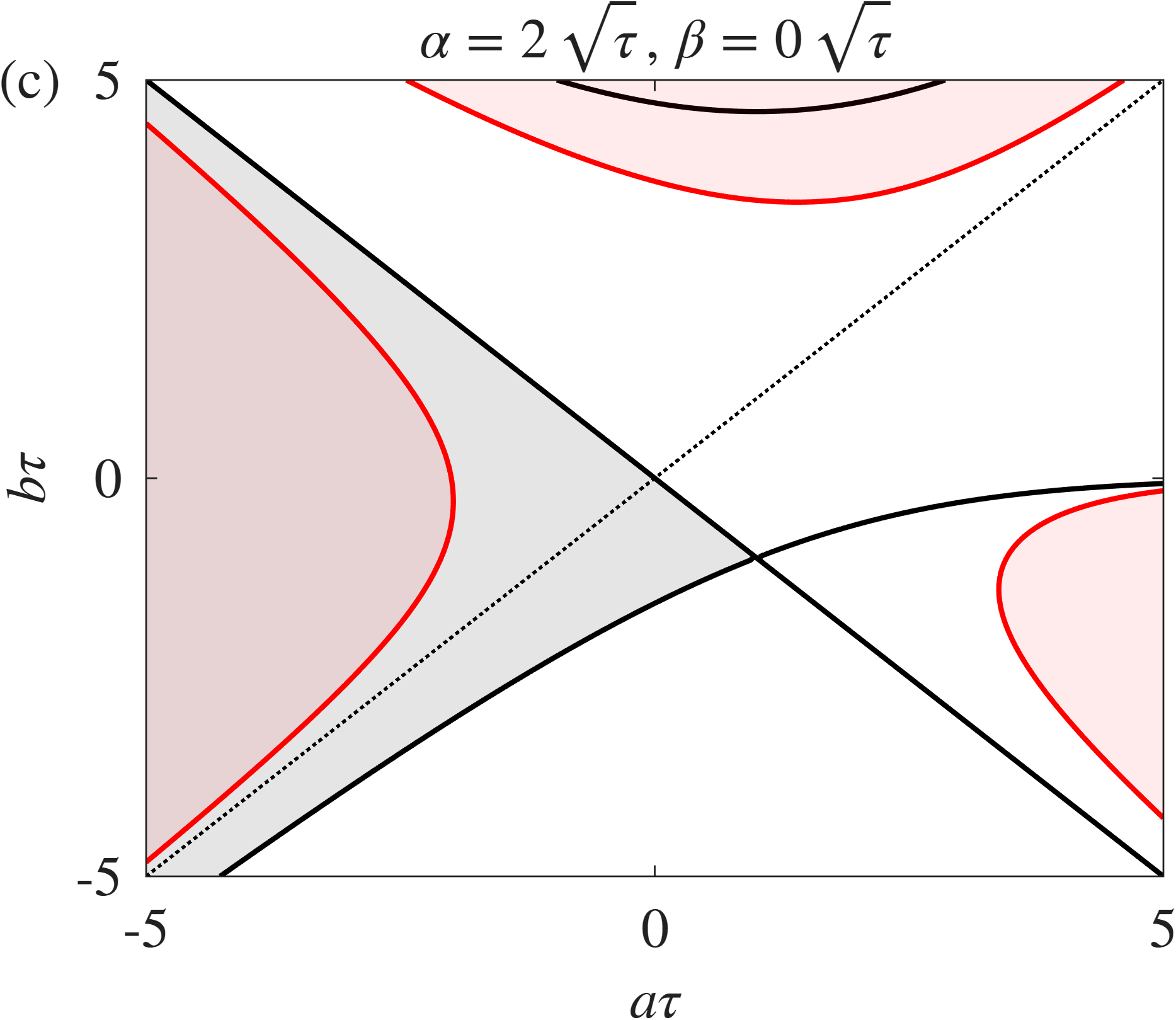}
    \includegraphics[width=0.33\linewidth]{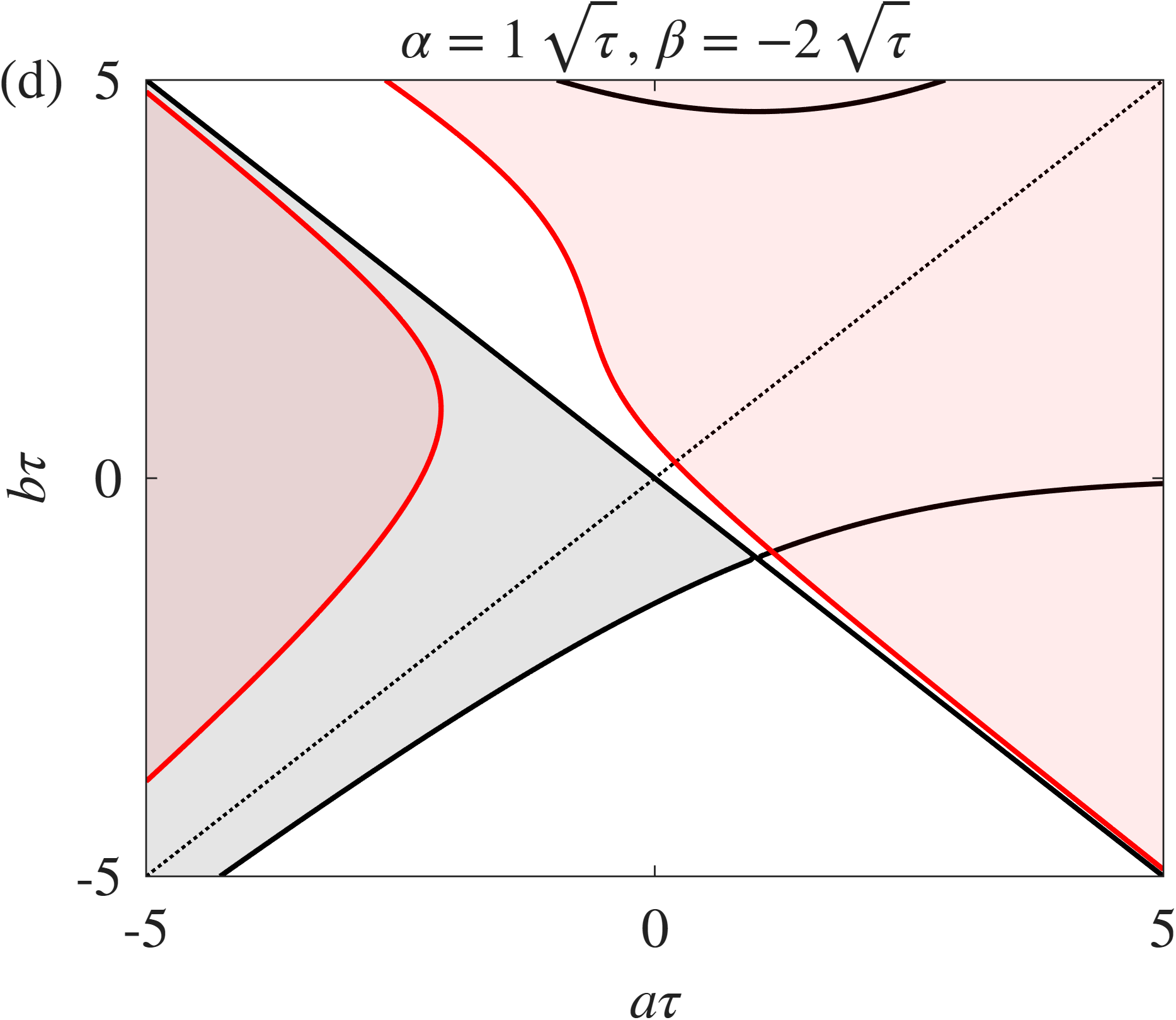}
    \includegraphics[width=0.33\linewidth]{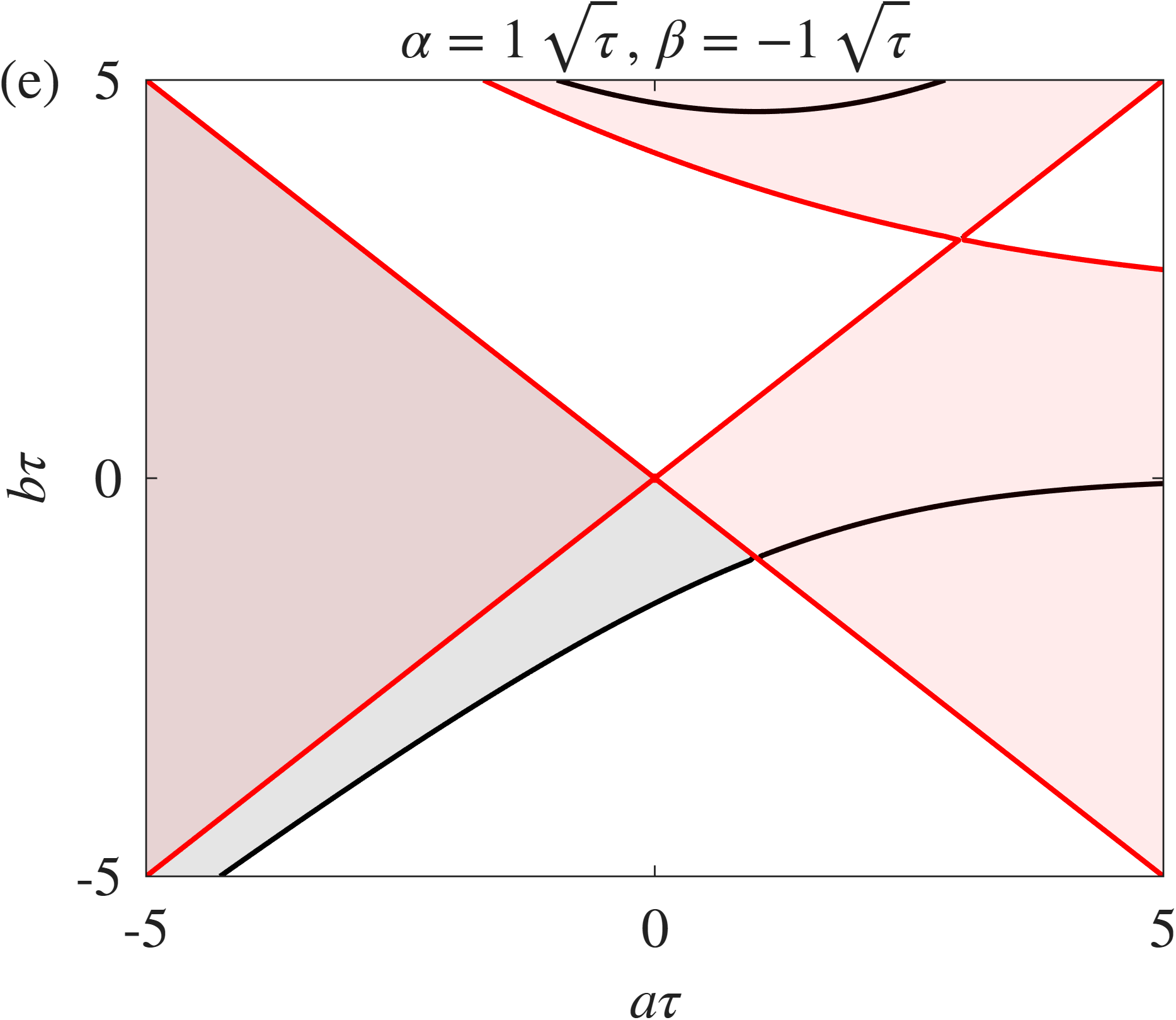}
    \includegraphics[width=0.33\linewidth]{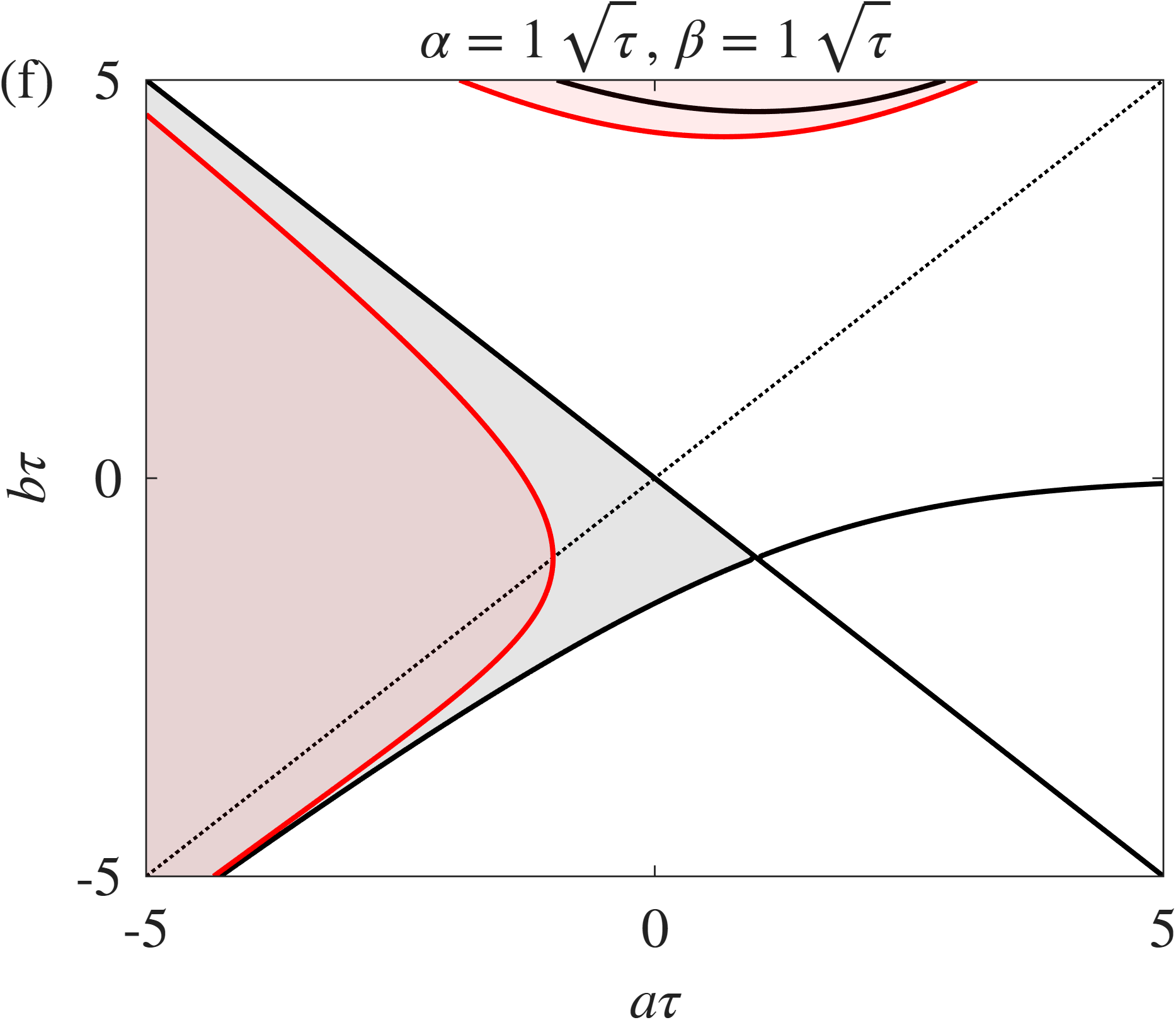}
    \caption{Parameter regions (shaded pink) in the $(a,b)$ plane for different combinations of $\alpha$ and $\beta$ where $\chi$, defined in Eq.~\eqref{eq:corr_chi}, is negative for the scalar problem. By symmetry, panels (a-c) show curves that are identical to those obtained by fixing $\alpha=0$ and setting $\beta=\{0,\,\sqrt{\tau},\,2\sqrt{\tau}\}$. Only the region with $\chi<0$ within the region of first-moment stability (shaded gray) corresponds to second-moment stability.}
    \label{fig:scalar_stab_ab}
\end{figure*}

It is straightforward to show that $\alpha=\beta=0$ implies $\chi=2(a+b)\eta_{\tau/2}(a,b)$. In this case, there exists a unique non-spurious stationary solution $\bar{\phi}$ throughout the entire region of first-moment stability.  This is consistent with the earlier observation that first-moment stability implies second-moment stability in the case of only additive noise. In general, along the saddle-node bifurcation curve $a+b=0$, $\chi=(\alpha+\beta)^2$, which is positive unless $\alpha=-\beta$, in which case it is identically equal to $0$. Similarly, along the Hopf bifurcation curve $\omega\mapsto(\omega\cot\omega\tau,-\omega\csc\omega\tau)$ for $\omega\tau\in[0,\pi)$, it holds that $\chi=(\alpha^2+\beta^2+2\alpha\beta\cos\omega\tau)\sec(\omega\tau/2)\ge0$ with equality only for $\alpha=-\beta$ at $\omega\tau=0$. It follows that the first-moment stability boundary does not intersect the region $\chi<0$. Fig.~\ref{fig:scalar_stab_ab} illustrates this observation for different choices of $\alpha$ and $\beta$.

We proceed to verify that the stationary solution
\begin{equation}
\bar{C}(\theta,\vartheta)=\bar{\phi}\left(-|\theta-\vartheta|\right).
\end{equation}
is a positive definite kernel provided that $\chi<0$ and first-moment stability holds. By a classical analysis, since $\bar{\phi}$ is continuous on the compact interval $[-\tau,0]$, this property of $\bar{C}$ can be established by showing that the eigenvalues of the integral operator
\begin{equation}
    \mathcal{I}[\varphi](s)=\int_{-\tau}^0\bar{\phi}(-|s-t|)\varphi(t)\,\mathrm{d}t
\end{equation}
are all positive \cite{ghanem2003stochastic}. Let $\lambda$ be such an eigenvalue with corresponding eigenfunction $\varphi_\lambda$. It follows from \eqref{eq:corr_scalar_solution} that $\varphi_\lambda$ must satisfy the second-order differential equation
\begin{equation}
    \varphi''(s)=\left(\mu^2-\frac{2\Delta\gamma^2}{\lambda\chi}\right)\varphi(s)
\end{equation}
and non-local boundary conditions
\begin{align}
    0&=\begin{pmatrix}a+b & a+b\\-\mu^2 & \mu^2\end{pmatrix}\begin{pmatrix}\varphi(0)\\\varphi(-\tau)\end{pmatrix}\nonumber\\
    &\qquad+\begin{pmatrix}-1 & 1\\a+b & a+b\end{pmatrix}\begin{pmatrix}\varphi'(0)\\\varphi'(-\tau)\end{pmatrix},
\end{align}
where $\Delta=(a+b)\eta_{\tau/2}(a,b)<0$. 

Consider, first, the case that $|a|>|b|$. Then, a negative eigenvalue must correspond to an eigenfunction of the form $\varphi_\lambda(s)=A\cosh \Omega s+B\sinh\Omega s$, where $\Omega^2=\mu^2-2
\Delta\gamma^2/\lambda\chi>\mu^2$. Substitution into the boundary conditions shows that nonzero values for $A$ and $B$ are possible only for $\Omega>0$ given by roots of the product
\begin{align}
    &4\cosh^2\frac{\Omega\tau}{2}\left((a+b)-\Omega\tanh\frac{\Omega\tau}{2}\right)\nonumber\\
    &\qquad\qquad\times\left((a+b)\Omega-\mu^2\tanh\frac{\Omega\tau}{2}\right).
\end{align}
Since both parenthetical factors are negative when $a+b<0$, we arrive at a contradiction.

Consider, next, the case when $|b|\ge|a|$, such that $\mu=\mathrm{i}\tilde{\mu}=\mathrm{i}\sqrt{b^2-a^2}$ and
\begin{equation}
    \Delta=(a+b)\cos\frac{\tilde{\mu}\tau}{2}+\tilde{\mu}\sin\frac{\tilde{\mu}\tau}{2}<0.\label{eq:Deltatrig}
\end{equation}
A negative eigenvalue must then correspond to an eigenfunction of the form $\varphi_\lambda(s)=A\cos\Omega s+B\sin\Omega s$, where $\Omega^2=\tilde{\mu}^2+2\Delta\gamma^2/\lambda\chi<\tilde{\mu}^2$. In this case, nonzero values for $A$ and $B$ are possible only for $\Omega\in(0,\tilde{\mu})$ given by roots of the product
\begin{align}
    &4\Omega\cos^2\frac{\Omega\tau}{2}\left((a+b)+\Omega\tan\frac{\Omega\tau}{2}\right)\nonumber\\
    &\qquad\qquad\times\left((a+b)+\frac{\tilde{\mu}^2}{\Omega}\tan\frac{\Omega\tau}{2}\right).
\end{align}
Since both $\Omega$-dependent terms in the parenthetical factors are increasing functions of $\Omega$, substitution of $\Omega=\tilde{\mu}$ and reference to \eqref{eq:Deltatrig} again leads to a contradiction.

\begin{figure*}[ht]
    \centering
    \includegraphics[width=0.9\linewidth]{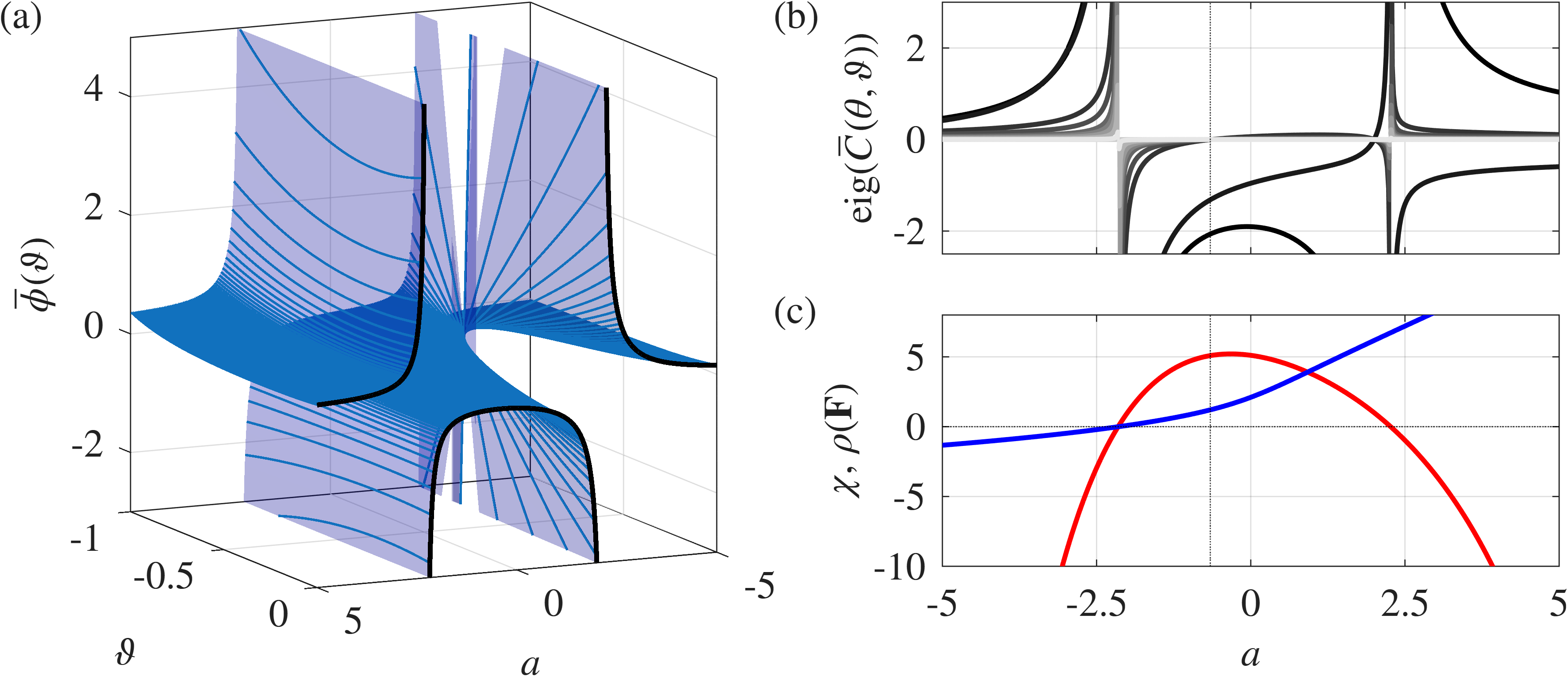}
    \caption{(a) Stationary solutions $\bar{\phi}(\vartheta)$ to the scalar correlation boundary-value problem; (b) eigenvalues of the corresponding discretized ($M=20$) covariance matrix in \eqref{eq:vec_cov} given $C_t(\theta,\vartheta)=\bar{C}(\theta,\vartheta)=\bar{\phi}\left(-|\theta-\vartheta|\right)$; and (c) variations of $\chi$ in Eq.~\eqref{eq:corr_chi} (red) and the spectral abscissa $\rho(\mathbf{F})$ (blue) for $\mathbf{F}=\mathbf{A}\oplus\mathbf{A}+\mathbf{B}\otimes\mathbf{B}$ and $\mathbf{A}$ and $\mathbf{B}$ in \eqref{eq:discA} and \eqref{eq:discB} under variations in $a$ and given $\tau = 1$, $b = -2$, $\alpha = -1.5$, $\beta = 0.5$, $\gamma = 1$. Second moment stability is guaranteed only for $a<-2.1543$, corresponding to the left-most zero-crossing of $\chi$ and $\rho(\mathbf{F})$. To the right of this point, $\bar{C}(\theta,\vartheta)$ is either negative definite or indefinite with a transition precisely where first-moment stability is lost at $a=-0.6380$.}
    \label{fig:scalar_sol_ina}
\end{figure*}

We conclude that $\bar{C}(\theta,\vartheta)$ is a positive definite kernel provided that $\chi<0$ and first-moment stability holds. This result is illustrated in Fig.~\ref{fig:scalar_sol_ina} which depicts $\bar{\phi}(\vartheta)$ and a discretized approximation of the spectrum of $\bar{C}(\theta,\vartheta)$ for $\tau=1$, $b=-2$, $\alpha=-1.5$, $\beta=0.5$, and $\gamma=1$ over a range of values of $a$. Here, $\chi<0$, $a+b<0$, and $\eta_{\tau/2}(a,b)>0$ provided that $a<-2.1543$. All other values of $a$ result in spurious solutions as confirmed by the negative definiteness or indefiniteness of the corresponding spectra.

\subsection{Special cases}
\label{sec:special cases}
Before generalizing our derivations to non-scalar problems, we consider several edge cases of the above analysis and tie these to results reported in the literature.

Consider, first, the scalar, linear, stochastic ODE
\begin{equation}
    \mathrm{d}x(t) = ax(t) \mathrm{d}t + \left(\alpha x(t)+\gamma\right) \mathrm{d}W_t. \label{eq:eq_scalar_sode}
\end{equation}
obtained by letting $b=\beta=0$ and, without loss of generality, $\tau=0$. In this case, a non-spurious stationary solution is given by
\begin{equation}
    \bar{\phi}(\vartheta)=c\left(\cosh|a|\vartheta-\mathrm{sign}(a)\sinh|a|\vartheta\right)=ce^{-a\vartheta},
\end{equation}
where $c$ is either uniquely defined by the ratio $-\gamma^2/(\alpha^2+2a)$ provided that $\chi=\alpha^2+2a<0$ or arbitrary non-negative if $\chi=\gamma=0$; no solution exists in the case that $\chi=0$ while $\gamma\neq0$. This result is in perfect agreement with those already available in the literature \cite{oksendal2003stochastic,sun2006stochastic}. The corresponding stationary solution $\bar{C}(\theta,\vartheta) = ce^{a|\theta-\vartheta|}$ is a multiple of the Ornstein-Uhlenbeck kernel \cite{williams2006gaussian}. 

Since first-moment stability here requires that $a<0$, we obtain a unique solution if and only if $\alpha^2<-2a$. This is equivalent to the well-known second-moment stability condition of the geometric Brownian motion (the special case of \eqref{eq:eq_scalar_sode} with $\gamma = 0$) \cite{baxendale2007stochastic,higham2000mean}. Since additive noise has no effect on the stability of stochastic differential equations, we conclude that the region in which the unique stationary solution exists and satisfies the positivity condition $\bar{\phi}(0)\ge 0$ coincides with the region of second-moment stability.

\begin{figure*}[ht]
    \centering
    \includegraphics[width=0.9\linewidth]{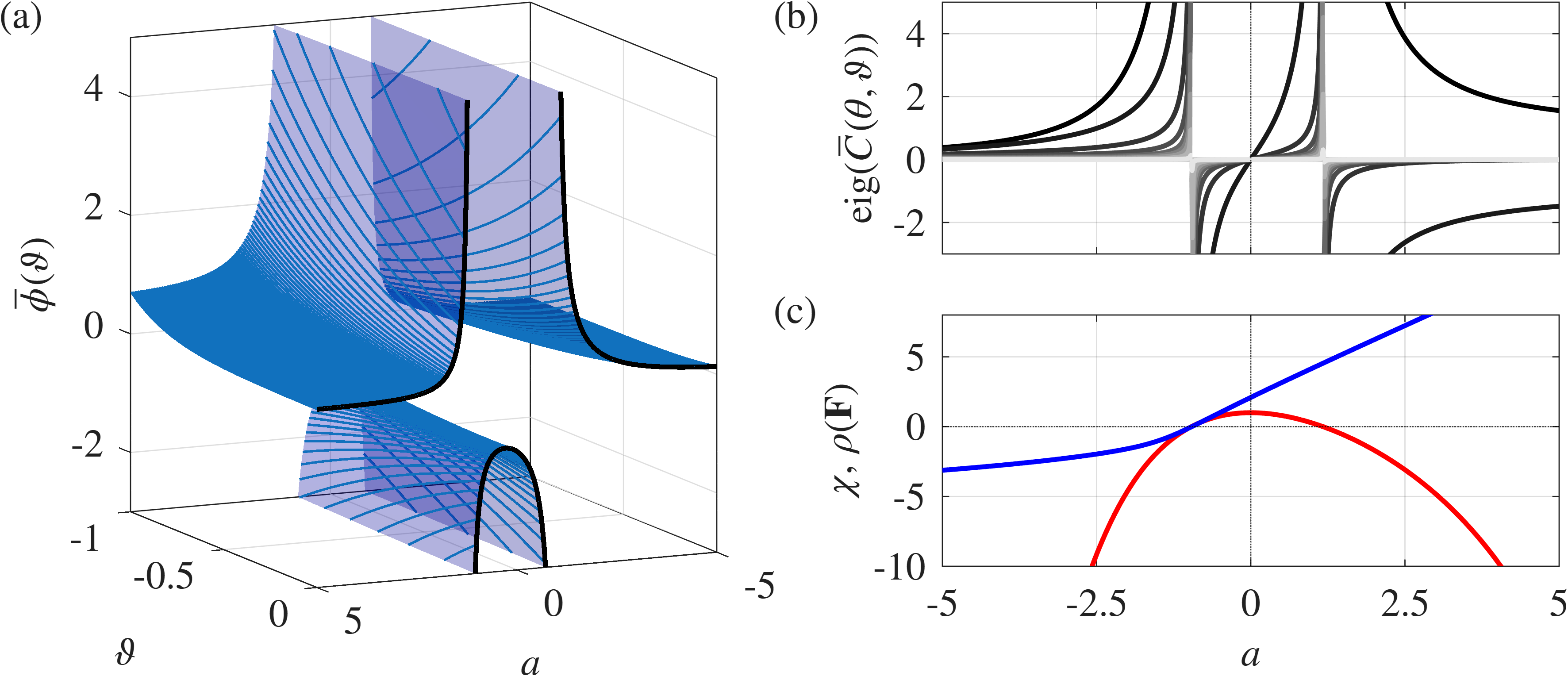}
    \caption{(a) Stationary solutions $\bar{\phi}(\vartheta)$ to the scalar correlation boundary-value problem; (b) eigenvalues of the corresponding discretized covariance  ($M=20$) matrix in \eqref{eq:vec_cov} given $C_t(\theta,\vartheta)=\bar{C}(\theta,\vartheta)=\bar{\phi}\left(-|\theta-\vartheta|\right)$; and (c) variations of $\chi$ in Eq.~\eqref{eq:corr_chi} (red) and the spectral abscissa $\rho(\mathbf{F})$ (blue) for $\mathbf{F}=\mathbf{A}\oplus\mathbf{A}+\mathbf{B}\otimes\mathbf{B}$ and $\mathbf{A}$ and $\mathbf{B}$ in \eqref{eq:discA} and \eqref{eq:discB}  under variations in $a$ and given $\tau = 1$, $b = 0$, $\alpha = -1.5$, $\beta = 0.5$, $\gamma = 1$ ($M=20$). Second moment stability is guaranteed only for $a<-0.9640$, corresponding to the left-most zero-crossing of $\chi$ and $\rho(\mathbf{F})$. To the right of this point, $\bar{C}(\theta,\vartheta)$ is either negative definite or indefinite with a transition precisely where first-moment stability is lost at $a=0$.}
    \label{fig:scalar_sol_ina_ex}
\end{figure*}

As a second edge case, consider the SDDE obtained when $b=\alpha=0$. In this case, the theory yields a non-spurious stationary solution of the form
\begin{equation}
    \bar{\phi}(\vartheta)=ce^{-a\left(\vartheta+\tau/2\right)},
\end{equation}
where $c$ is either uniquely defined by the ratio $-\gamma^2e^{a\tau/2}/(\beta^2+2a)$ provided that $e^{-a\tau/2}\chi=\beta^2+2a<0$
or arbitrary non-negative if $\beta^2+2a=\gamma=0$; no solution exists in the case that $\beta^2+2a=0$ while $\gamma\ne 0$. The corresponding stationary solution $\bar{C}(\theta,\vartheta)$ is another scaled version of the Ornstein-Uhlenbeck kernel. For this case, the authors of \cite{sykora2019stochastic} show that $\beta^2+2a<0$ guarantees second-moment stability. Here, again, the region of second-moment stability coincides with the region of existence of a unique stationary solution $\bar{\phi}$ that satisfies the positivity condition $\bar{\phi}(0)\ge 0$.

As a generalization of the previous two edge cases, consider the SDDE obtained when $b=0$ but $\alpha$ and $\beta$ are arbitrary. In this case, the theory predicts a unique stationary solution that satisfies the positivity condition $\bar{\phi}(0)\ge 0$ provided that $a<0$ and
\begin{equation}
    \chi=e^{-a\tau/2}\left(2a+\alpha^2+\beta^2+2\alpha\beta e^{a\tau}\right)<0.
\end{equation}
This is consistent with the result in Example 3.5 in Appleby, Mao, and Riedle~\cite{appleby2009geometric}, where second-moment stability follows from the necessary and sufficient condition
\begin{equation}
    \int_0^\tau \alpha^2e^{2at}\,\mathrm{d}t+\int_\tau^\infty\left(\alpha e^{at}+\beta e^{a(t-\tau)}\right)^2\,\mathrm{d}t<1.
\end{equation}
For $\tau=1$,$\alpha=-1.5$, and $\beta=0.5$, Fig.~\ref{fig:scalar_sol_ina_ex}(c) shows that this condition is satisfied for $a<0$ as long as $a<-0.9640$.

\begin{figure}[ht]
    \centering
    \includegraphics[width=0.9\linewidth]{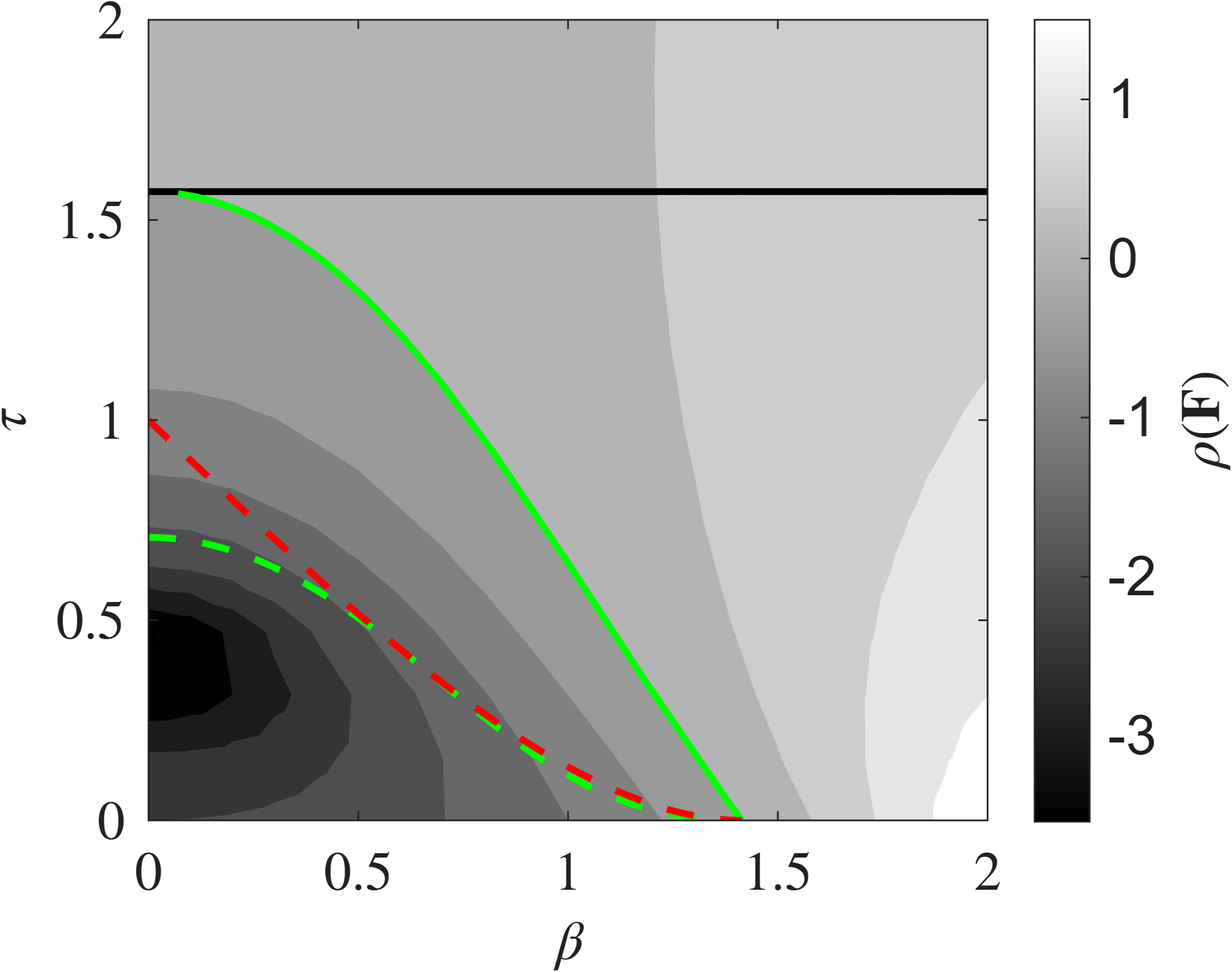}
    \caption{Candidate second-moment stability boundaries for the scalar problem with $a=\alpha=\gamma=0$ and $b=-1$ and the spectral abscissa $\rho(\mathbf{F})$ for $\mathbf{F}=\mathbf{A}\oplus\mathbf{A}+\mathbf{B}\otimes\mathbf{B}$ and $\mathbf{A}$ and $\mathbf{B}$ in \eqref{eq:discA} and \eqref{eq:discB}  under variations in $\tau$ and $\beta$  ($M=20$). Solid black line: first-moment stability boundary. Dashed red line: limit from Eq.~\eqref{eq:Haskovec_cond}. Dashed green line: limit from Eq.~\eqref{eq:HaskovecMao_cond}. Solid green line: limit from Eq.~\eqref{eq:Haskovec}. The figure confirms the suboptimality of the first two conditions.}
    \label{fig:scalar_map_atau}
\end{figure}

Finally, we consider the special case with $a=\alpha=0$ and $b=-1$, which sits securely within the region of first-moment stability. For this parameter combination, Haskovec~\cite{havskovec2022asymptotic} cites the sufficient condition
\begin{equation}
    \beta^2<2,\,\tau<1-\sqrt{\beta^2-\frac{\beta^4}{4}}\label{eq:Haskovec_cond}
\end{equation}
and Erban, Haskovec, and Sun~\cite{erban2016cucker} cite the sufficient condition
\begin{equation}
    \beta^2<2,\,\tau<\frac{1}{4}\left(-2\beta^2+\sqrt{4\beta^4+2(2-\beta^2)^2}\right)\label{eq:HaskovecMao_cond}
\end{equation}
for second-moment stability. By considering the noise-free limit, they note that neither condition is optimal, since first-moment stability in this case implies second-moment stability for $\tau<\pi/2$. In this case, the theory in the previous section reduces to the condition
\begin{equation}
    \chi=\left(\beta^2-2\right)\cos\frac{\tau}{2}+\left(\beta^2+2\right)\sin\frac{\tau}{2}<0.\label{eq:Haskovec}
\end{equation}
In particular, for $\beta^2<2$ and $\tau$ increasing from $0$, we predict that second-moment stability is lost at
\begin{equation}
    \tau=\tau_c=2\arctan\frac{2-\beta^2}{2+\beta^2}.
\end{equation}
In the noise-free limit, this reduces to $\tau<\pi/2$, suggesting that the condition in \eqref{eq:Haskovec} is optimal. Fig.~\ref{fig:scalar_map_atau} compares all three of these predicted second moment stability bounds against the numerically computed spectral abscissa of the coefficient matrix $\mathbf{F}=\mathbf{A}\oplus\mathbf{A}+\mathbf{B}\otimes\mathbf{B}$ of the surrogate covariance ODE in \eqref{eq:cov_disc}. 

\subsection{Generalization to non-scalar problems}
\label{sec:generalcase}

We proceed to consider the unrestricted form of Eq.~\eqref{eq:Ito_sdde} where $n>1$. In this case, stationary solutions of the correlation boundary-value problem \eqref{eq:corr_pde}-\eqref{eq:corr_pde_bc} with $m_t\equiv 0$ satisfy the differential equations
\begin{align}
    \phi'(\vartheta) = -a\phi(\vartheta) - b\phi^\mathsf{T}(-\tau-\vartheta)\label{eq:corr_pde_stat}
\end{align}
for $\vartheta\in[-\tau,0)$ and
\begin{align}
    \phi'(\vartheta)=-a\phi(\vartheta)-b\phi(\vartheta+\tau)\label{eq:corr_pde2_stat}
\end{align}
for $\vartheta\in(-\infty,-\tau)$, constrained by the non-local, symmetric boundary condition
\begin{align}
    0&=2\big[a\phi(0)+\phi(-\tau)b^\mathsf{T}+\alpha \phi(-\tau)\beta^\mathsf{T}\big]_s\nonumber\\
    &\quad+\beta \phi(0)\beta^\mathsf{T}+\alpha \phi(0)\alpha^\mathsf{T}+\gamma\gamma^\mathsf{T}\label{eq:corr_pde_bc_stat}
\end{align}
and the symmetry condition $\phi^\mathsf{T}(0)=\phi(0)$. 

Assuming the existence of a stationary solution $\bar{\phi}$ of \eqref{eq:corr_pde_stat}-\eqref{eq:corr_pde_bc_stat}, let $\hat{\phi}(\sigma)=\bar{\phi}(-\sigma-\tau)$ for $\sigma\in(0,\infty)$. Then, Eq.~\eqref{eq:corr_pde2_stat} implies that
\begin{equation}
    \hat{\phi}'(\sigma)=a\hat{\phi}(\sigma)+b\hat{\phi}(\sigma-\tau),\label{eq:corr_stat_dde}
\end{equation}
where $'$ denotes differentiation with respect to $\sigma$. Since this is identical in form to the deterministic part of the original SDDE \eqref{eq:Ito_sdde} and since we already assume first-moment stability, it follows that
\begin{equation}
    \lim_{\vartheta \rightarrow  -\infty}\bar{\phi}(\vartheta) = \lim_{\vartheta \rightarrow  -\infty}\hat{\phi}(-\vartheta-\tau)=0.
\end{equation}
With this in mind, we focus our attention on Eqs.~\eqref{eq:corr_pde_stat} and \eqref{eq:corr_pde_bc_stat}. Unfortunately, since $a$ and $b$ are matrices, it is no longer the case that $\bar{\phi}$ satisfies a second-order differential equation of the form in \eqref{eq:corr_2nd_ode}.

Notably, as $\vartheta$ increases from $-\tau$, the argument $-\vartheta-\tau$ in the second term on the right-hand side of \eqref{eq:corr_pde_stat} decreases from $0$ with equality achieved when $\vartheta=-\tau/2$. Consistent with this observation, we define the functions $\bar{f}(\vartheta)=\bar{\phi}(\vartheta-\tau/2)$ and $\bar{g}(\vartheta)=f(-\vartheta)$ such that $\bar{f}(0)=\bar{g}(0)$, $\bar{f}(\tau/2)$ is symmetric, and $(\bar{f}, \bar{g})$ satisfy the homogeneous LTI system of first-order ordinary differential equations
\begin{equation}
    \bar{f}' = -a\bar{f}-b\bar{g}^\mathsf{T},\,\bar{g}' = a\bar{g}+b\bar{f}^\mathsf{T}\label{eq:corr_sub_stat_ode}
\end{equation} 
with the symmetric terminal condition
\begin{align}
    0&=2\left[a\bar{f}\left(\frac{\tau}{2}\right)+\bar{g}\left(\frac{\tau}{2}\right)b^\mathsf{T}+\alpha\bar{g}\left(\frac{\tau}{2}\right)\beta^\mathsf{T}\right]_s\nonumber\\
    &\qquad+\alpha \bar{f}\left(\frac{\tau}{2}\right)\alpha^\mathsf{T}+\beta \bar{f}\left(\frac{\tau}{2}\right)\beta^\mathsf{T}+\gamma\gamma^\mathsf{T}.
\end{align}
It follows that $\bar{f}_\mathfrak{v}=\mathfrak{vec}(\bar{f})$ and $\bar{g}_\mathfrak{v}=\mathfrak{vec}(\bar{g})$ satisfy the LTI system
\begin{equation}    \begin{bmatrix}\bar{f}'_\mathfrak{v}\\\bar{g}'_\mathfrak{v}\end{bmatrix}=A\begin{bmatrix}\bar{f}_\mathfrak{v}\\\bar{g}_\mathfrak{v}\end{bmatrix}, \quad A = \begin{bmatrix}-(I_n\otimes a) & -(I_n\otimes b)P_n\\(I_n\otimes b)P_n & (I_n\otimes a)\end{bmatrix},\label{eq:vectorizeddynamics}
\end{equation}
and boundary conditions
\begin{align}
0&=\bar{f}_\mathfrak{v}(0)-\bar{g}_\mathfrak{v}(0),\\
0&=\left(P_n-I_{n^2}\right)\bar{f}_\mathfrak{v}(\tau/2),\\
    0&=B_f \bar{f}_\mathfrak{v}(\tau/2) + B_g \bar{g}_\mathfrak{v}(\tau/2)+(\gamma\gamma^\mathsf{T})_\mathfrak{v},
\end{align}
where
\begin{equation}
    B_f = a \oplus a +\alpha\otimes\alpha+\beta\otimes\beta
\end{equation}
and
\begin{equation}
    B_g = \left(I_{n^2}+P_n\right)\left(b\otimes I_n+\beta\otimes\alpha\right).
\end{equation}

These equations have a unique solution in the form 
\begin{equation}
    \begin{bmatrix}\bar{f}_\mathfrak{v}(\vartheta)\\\bar{g}_\mathfrak{v}(\vartheta)\end{bmatrix} = e^{A\vartheta}\begin{bmatrix}\bar{f}_{\mathfrak{v},0}\\\bar{f}_{\mathfrak{v},0}\end{bmatrix},\label{eq:stat_sol_nonscalar}
\end{equation}
provided that the matrix
\begin{equation}
    \Psi = \begin{bmatrix}
        Q_n\begin{bmatrix}
            B_f & B_g \end{bmatrix} \\ 
        R_n\begin{bmatrix}
            P_n-I_{n^2} & 0_{n^2}
        \end{bmatrix} 
    \end{bmatrix}e^{A\tau/2}\begin{bmatrix}
            I_{n^2}\\I_{n^2}
    \end{bmatrix}\label{eq:sing_cond_mat}
\end{equation} 
is non-singular. Here $Q_n \in \mathbb{R}^{\frac{(n+1)n}{2} \times n^2}$ represent the projection of a vectorized symmetric matrix to its lower triangular elements, including the diagonal terms, while $R_n \in \mathbb{R}^{\frac{(n-1)n}{2} \times n^2}$ projects to the upper triangle elements, excluding the diagonal entries. The co-dimension-one surfaces implicitly defined by $\det\,(\Psi) = 0$ identify parameter combinations where uniqueness cannot be guaranteed, i.e., candidate second-moment stability boundaries.

If $\det\,(\Psi)\neq 0$, we may find the missing constants in \eqref{eq:stat_sol_nonscalar} by solving the linear algebraic equations
\begin{equation}
\Psi\bar{f}_{\mathfrak{v},0} =
    -\begin{bmatrix}
        Q_n(\gamma\gamma^\mathsf{T})_\mathfrak{v}\\
        R_n 0_{n,\mathfrak{v}}
    \end{bmatrix}.
\end{equation}
In this case, the vectorized stationary solution is given by
\begin{equation}
  \bar{\phi}_\mathfrak{v}(\vartheta) = \begin{bmatrix}
      I_{n^2} & 0_{n^2}
  \end{bmatrix}e^{A(\vartheta+\tau/2)}\begin{bmatrix}
            I_{n^2}\\I_{n^2}
    \end{bmatrix}\bar{f}_{\mathfrak{v},0}.
\end{equation}
In particular, we anticipate ascertaining whether the region of second-moment stability is found for $\det\,(\Psi)<0$ or $>0$ by investigating the sign of the first component of
\begin{equation}
  -\begin{bmatrix}
      I_{n^2} & 0_{n^2}
  \end{bmatrix}e^{A\tau/2}\begin{bmatrix}
            I_{n^2}\\I_{n^2}
    \end{bmatrix}\mathrm{adj}(\Psi)\begin{bmatrix}
        Q_n(\gamma\gamma^\mathsf{T})_\mathfrak{v}\\
        R_n 0_{n,\mathfrak{v}}\end{bmatrix}.\label{eq:signcondition}
\end{equation}
at a point where $\det\,(\Psi)=0$.

In the scalar case, $P_n=Q_n=1$ and $R_n=\emptyset$ imply that
\begin{equation}
    e^{A\vartheta}=\begin{pmatrix}\cosh\mu\vartheta-\frac{a}{\mu}\sinh\mu\vartheta & -\frac{b}{\mu}\sinh\mu\vartheta\\\frac{b}{\mu}\sinh\mu\vartheta & \cosh\mu\vartheta+\frac{a}{\mu}\sinh\mu\vartheta\end{pmatrix}
\end{equation}
and, consequently, that $\det\,(\Psi)$ equals $\chi$ given in Eq.~\eqref{eq:corr_chi}. Substitution in \eqref{eq:signcondition} yields $-\gamma^2\eta_{\tau/2}(a,-b)$, from which we again conclude that second-moment stability in the scalar case is associated with first-moment stability and the condition $\chi<0$.

\section{Numerical examples}\label{sec:numerical examples} 

As described in Sec.~\ref{sec:discretization} and consistent with the approach taken in \cite{torkamani2014numerical}, the stability of the second-moment dynamics may be assessed approximately by reference to the spectral abscissa $\rho(\mathbf{F})$ of the coefficient matrix $\mathbf{F}=\mathbf{A}\oplus\mathbf{A}+\mathbf{B}\otimes\mathbf{B}$ of the discretized, high-dimensional, surrogate covariance ODE in \eqref{eq:cov_disc}. Here, stability is guaranteed by $\rho(\mathbf{F})<0$. For a sufficiently large value of $M$, Figs.~\ref{fig:scalar_sol_ina}, \ref{fig:scalar_sol_ina_ex}, and \ref{fig:scalar_map_atau} validate that loss of stability is associated with a crossing of the zero-level set of $\chi$ in the increasing direction. In this section, we provide additional examples that validate this theory, including for non-scalar problems.

\subsection{Implementation}
Although stability assessment of the surrogate dynamics can produce reliable and accurate predictions, it should be noted that since $\mathbf{F}\in\mathbb{R}^{M^2n^2\times M^2n^2}$ and not necessarily sparse, the costs of evaluating its spectrum or even just its spectral abscissa scale rather poorly with both $M$ and $n$. Following \cite{torkamani2014numerical}, by exploiting the symmetry of $C_t(\theta,\vartheta)$, one may reduce consideration to a $M^2n^2(M^2n^2+1)/2$-dimensional subspace, but this does not help much for large values of $n$ or $M$, especially if the calculation has to be repeated for many different parameter combinations.

As an alternative, in light of the anticipated positive semi-definiteness of a discretization $\bar{\mathbf{C}}\in\mathbb{R}^{Mn\times Mn}$ of the stationary solution $\bar{C}$ for some nonzero $\gamma$, second-moment stability may be assessed by requiring that all eigenvalues of $\bar{\mathbf{C}}$ be greater than $-\varepsilon$ for some numerical threshold $0<\varepsilon\ll 1$. A candidate $\bar{\mathbf{C}}$ may be obtained using interpolation from the stationary solution to the surrogate correlation DDE \eqref{eq:corr_disc_eq} derived in \ref{app:num_corr}. Since the latter is obtained by solving a system of $Mn^2$ linear algebraic equations in $Mn^2$ unknowns, this stability assessment technique scales better by a factor of $M$ than direct spectral analysis of the coefficient matrix $\mathbf{F}$.

Given a parameter combination in the region of second-moment stability, loss of second-moment stability can be detected by a crossing through $0$ of $\det\,(\mathbf{F})$ or by a sign change of the largest eigenvalue of $\bar{\mathbf{C}}$. As discussed in \ref{app:num_corr}, it is computationally advantageous to monitor for a crossing through $0$ of $\det\,(\boldsymbol{\Xi})$, as both dimensions of the matrix $\boldsymbol{\Xi}$ are smaller than those of $\mathbf{F}$ by a factor of $M$. Indeed, if one is interested only in tracking approximate stability boundaries using a discretized problem formulation, we recommend using parameter continuation techniques applied to the \textit{zero problem}
\begin{equation}
    \begin{bmatrix}
        \boldsymbol{\Xi}\mathbf{v}\\\mathbf{v}^\mathsf{T}\mathbf{v}-1
    \end{bmatrix}=
        \mathbf{0}_{(Mn^2+1)\times 1}
\end{equation}
in terms of the unknown problem parameters and column matrix $\mathbf{v}\in\mathbb{R}^{Mn^2}$, assuming that an adequate initial solution guess has first been found~\cite{dankowicz2013recipes}.

\begin{figure*}[!h]
    \centering
    \includegraphics[width=0.4\linewidth]{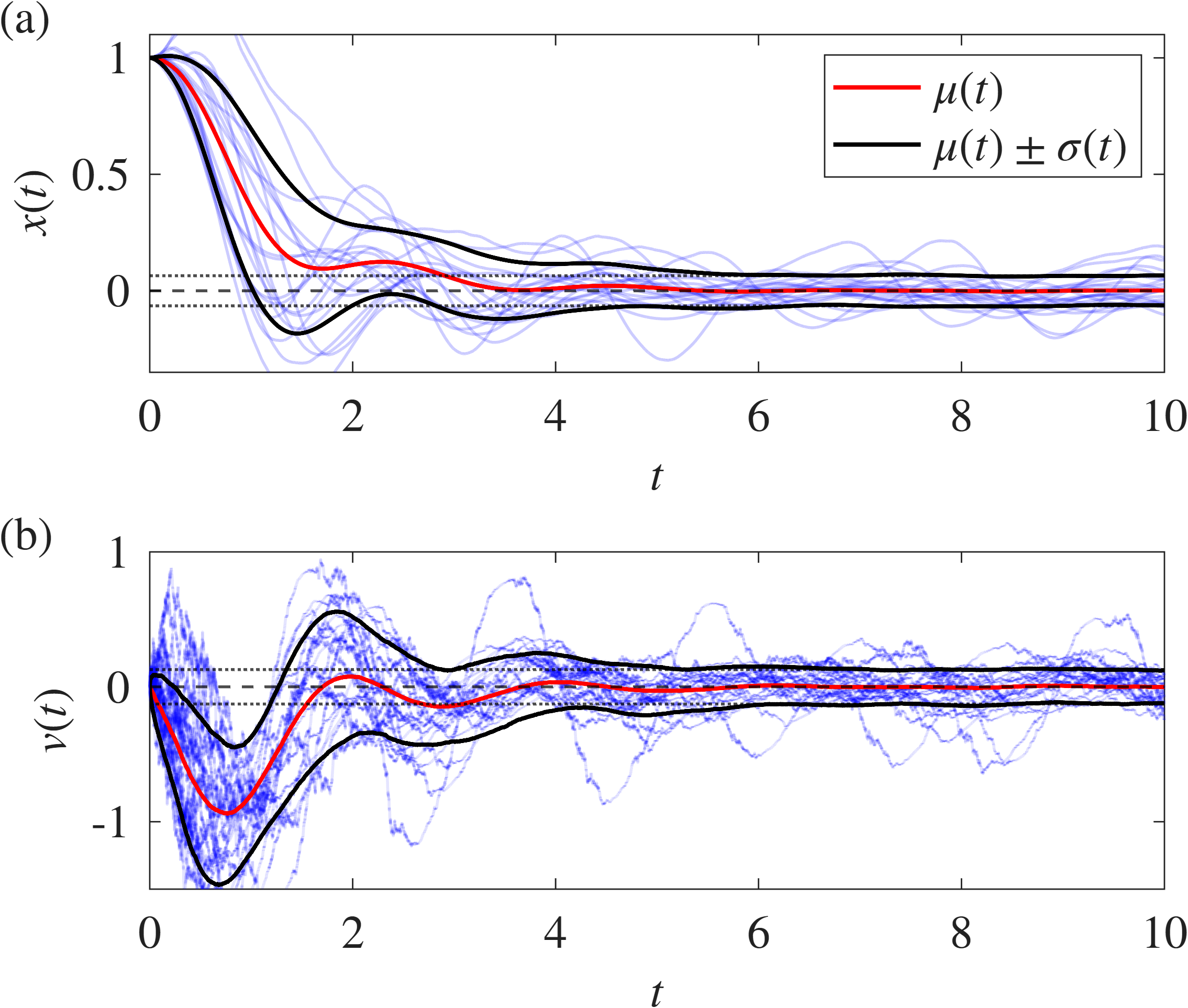}
    \includegraphics[width=0.4\linewidth]{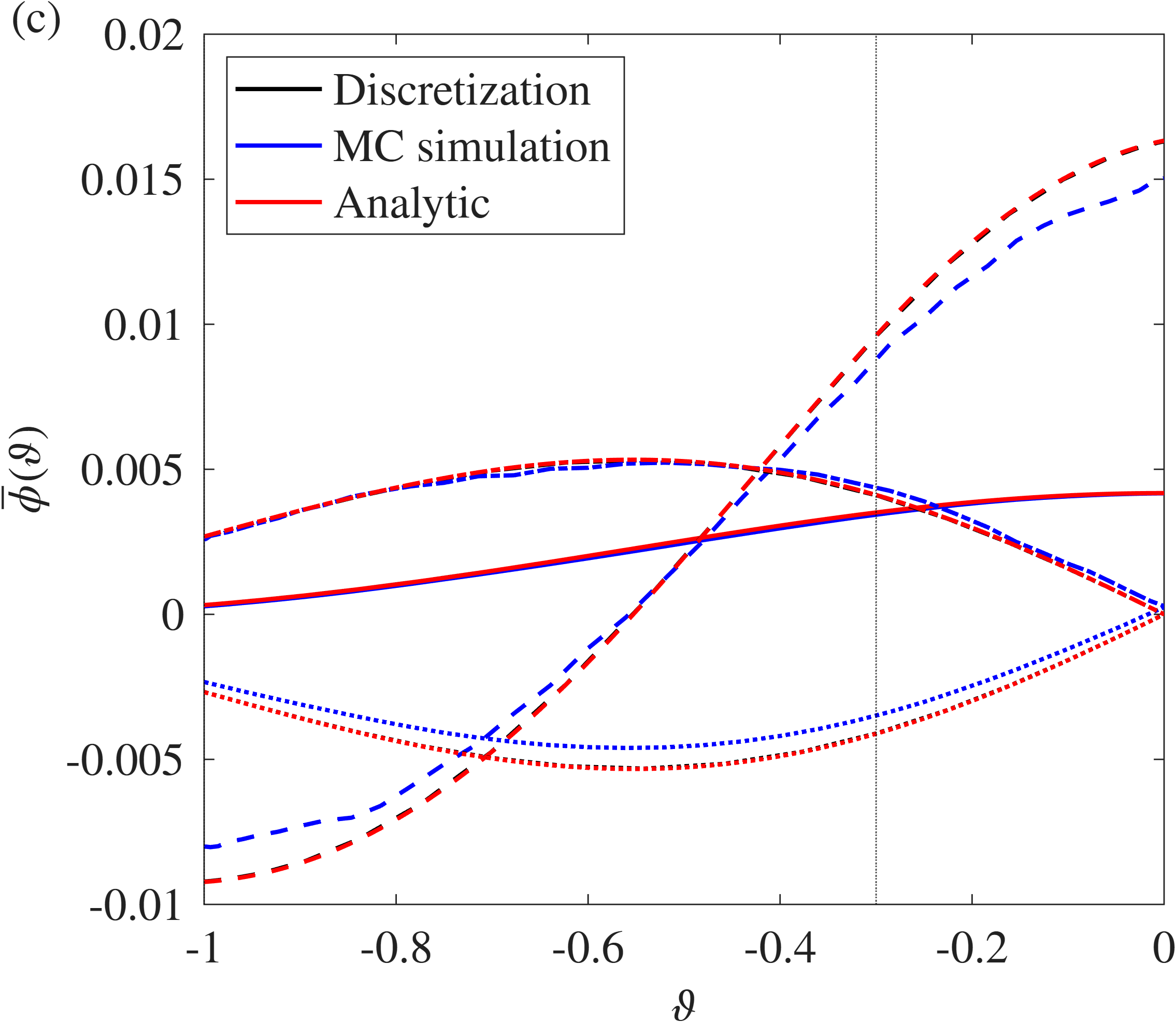}\\
    \includegraphics[width=0.4\linewidth]{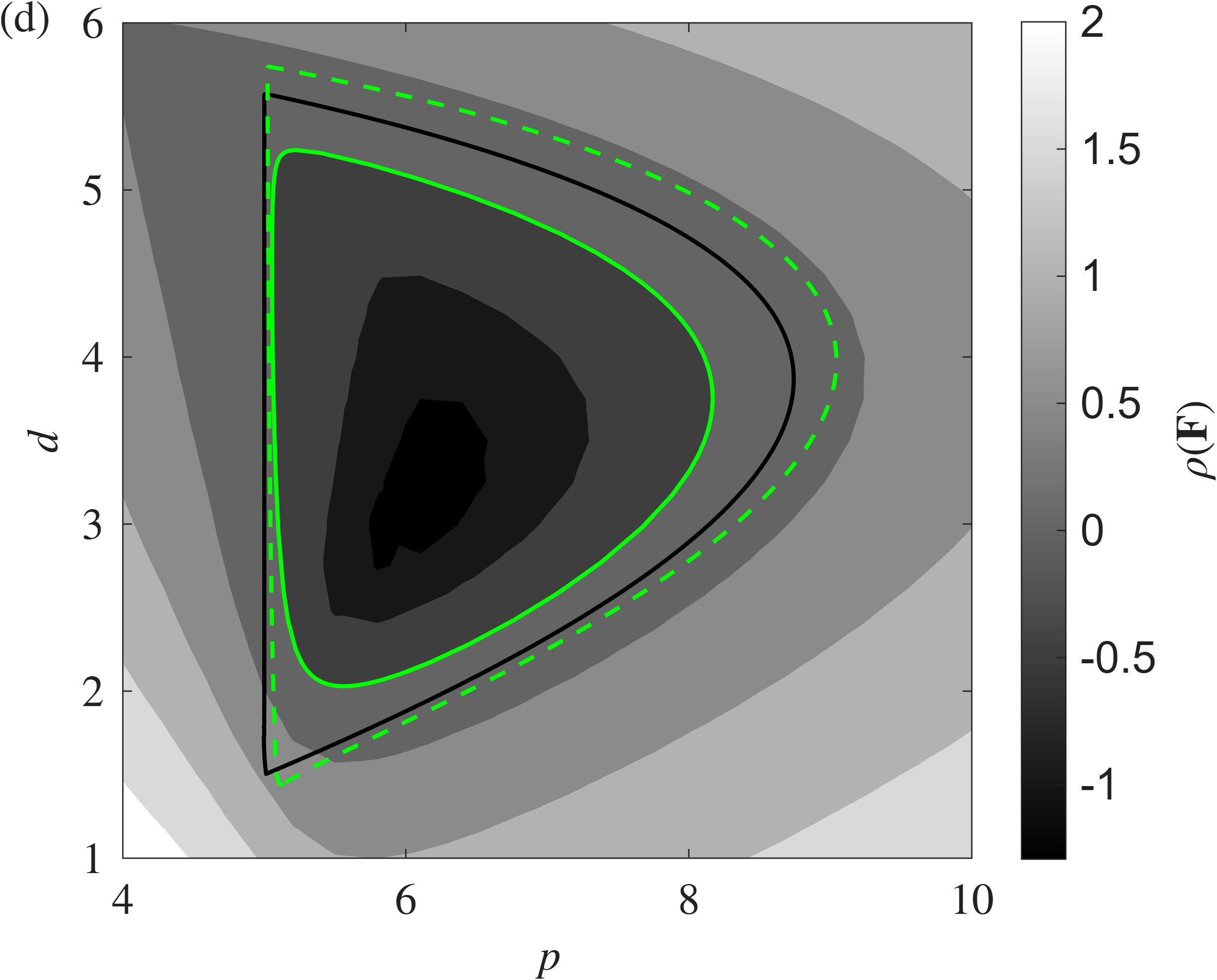}
    \includegraphics[width=0.4\linewidth]{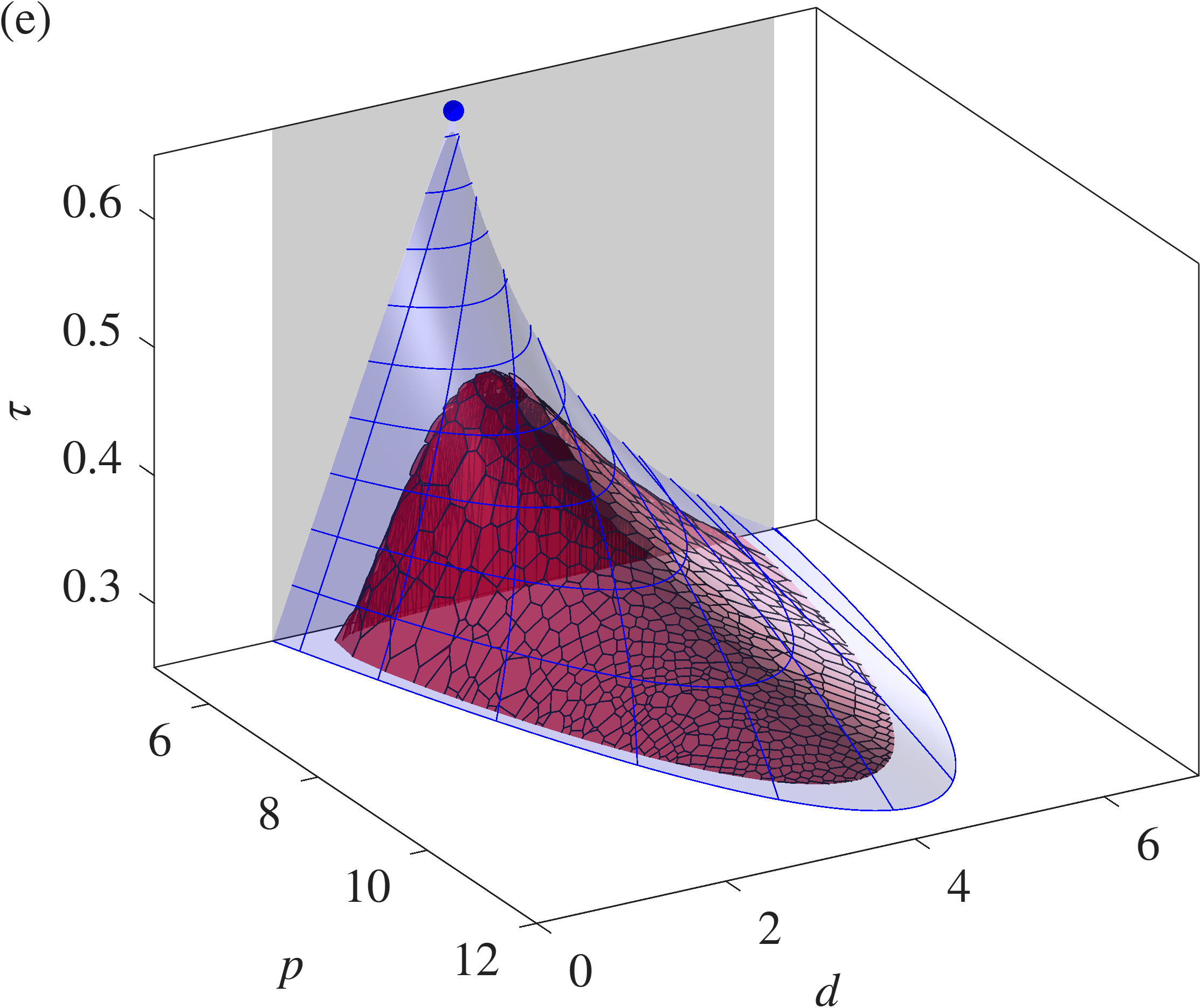}
    \caption{Numeric study of a noisy, single-degree-of-freedom model of an inverted pendulum subjected to delayed feedback control ($p=6.5$, $d=3.5$, $\sigma = 0.1$, $k=5$, $\tau = 0.3$). Panels (a)-(b) show a subset of simulated stochastic position and velocity trajectories (blue), time histories for the sample mean $\mu(t)$ (red), and deviations from the sample mean by one sample standard deviation (black) . In panel (c), a favorable comparison is made between predictions of the semi-analytic theory, problem discretization, and forward simulation, respectively, with solid lines marking $\mathbb{E}[x(T)x(T+\vartheta)]$, dashed lines marking $\mathbb{E}[v(T)v(T+\vartheta)]$, dash-dotted lines marking $\mathbb{E}[v(T)x(T+\vartheta)]$, and dotted lines marking $\mathbb{E}[x(T)v(T+\vartheta)]$. Panel (d) presents the spectral abscissa $\rho(\mathbf{F})$ under variations in $p$ and $d$, along with the first moment stability boundary (black line), the second-moment stability limits found via $\det\,(\Psi)=0$ (solid green line) and the condition associated with Eq.~\eqref{eq:Buckwar2} (dashed green line). Panel (e) illustrates the evolution of the stability bounds as $\tau$ is increased: blue surface: first-moment Hopf boundary; gray plane: first-moment saddle-node boundary; red surface: second-moment stability boundary obtained from multi-dimensional continuation using the \textsc{COCO} software package~\cite{dankowicz2020multidimensional}.}
    \label{fig:osc_stab}
\end{figure*}

In contrast to discretization-based techniques, the theoretical results in Sec.~\ref{sec:analytical results} necessitate the computation of only $2n^2$ eigenvalues of the matrix $A$ in Eq.~\eqref{eq:vectorizeddynamics} in order to compute the matrix exponential in the definition of $\Psi$. The particular form of $A$ implies that these come in either pairs that add to $0$ or quadruplets that add pairwise to either $0$ real part or $0$ imaginary part, or both. Alternatively, one may consider various forms of approximation of the matrix exponential $e^{tA}$ that eliminate the need to compute its eigenvalues and instead render an expression for $\Psi$ entirely in terms of the problem parameters. Most importantly, however, the matrix exponential is a smooth function of the matrix entries of $A$ and can therefore be integrated in a zero problem of the form
\begin{equation}
    \begin{bmatrix}
        \Psi v\\v^\mathsf{T}v-1
    \end{bmatrix}=
        0 \label{eq:bordered}
\end{equation}
in terms of an unknown column matrix $v\in\mathbb{R}^{n^2}$. Given an initial solution guess for $v$ and a parameter combination near the second-moment stability boundary, parameter continuation can then be applied to map out the stability boundary over a larger region of parameter space. A particularly compelling use of such techniques is the computation using the software package \textsc{COCO}~\cite{dankowicz2020multidimensional} of the two-dimensional second-moment stability boundary in Fig.~\ref{fig:osc_stab}.

\subsection{The scalar case}
\label{sec:numerics for scalar case}
\begin{figure}[t]
    \centering
    \includegraphics[width=0.9\linewidth]{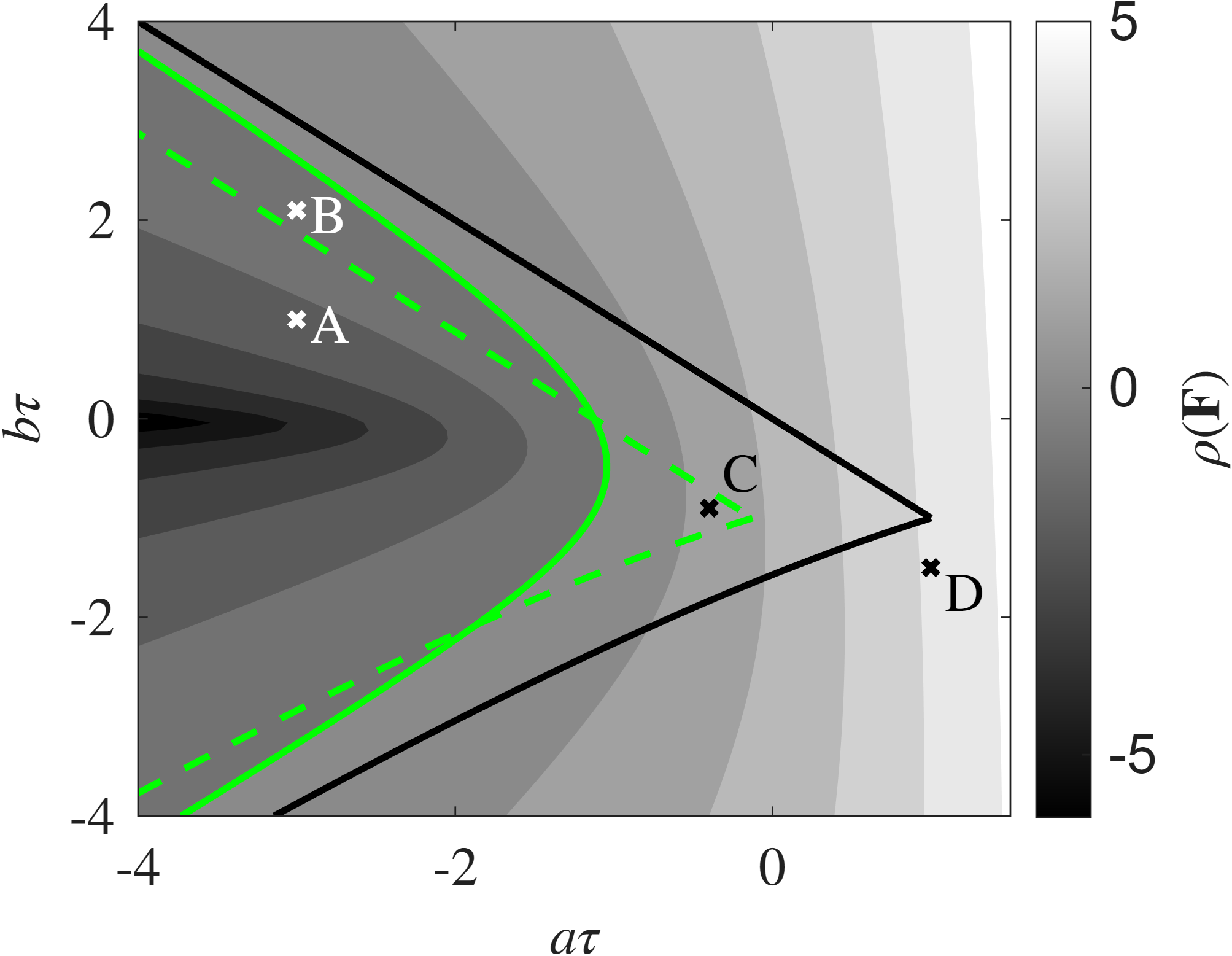}
    \caption{Candidate second-moment stability boundaries for the scalar problem with $\alpha=-1.5\sqrt{\tau}$, $\beta=0$ and the spectral abscissa $\rho(\mathbf{F})$ for $\mathbf{F}=\mathbf{A}\oplus\mathbf{A}+\mathbf{B}\otimes\mathbf{B}$ and $\mathbf{A}$ and $\mathbf{B}$ in \eqref{eq:discA} and \eqref{eq:discB}  under variations in $a$ and $b$  ($M=20$). Solid black line: first-moment stability boundary. Solid green line: limit from Eq.~\eqref{eq:corr_chi}. Dashed green line: limit from Eq.~\eqref{eq:Mackey}. The latter is seen to be neither necessary nor sufficient.}
    \label{fig:scalar_Stab_comp}
\end{figure}

We restrict attention again to the scalar case. Consider, for example, the parameter combination $\tau=1$, $\alpha=-1.5$, and $\beta=0$. Then, as shown in Fig.~\ref{fig:scalar_Stab_comp}, the curve $\chi=0$ bounds a region of second-moment stability that is a proper subset of the region of first-moment stability. Shown by a dashed curve is the prediction according to \cite{mackey1995solution}, such that all solutions of \eqref{eq:Mackey} have negative real part to the left of this curve. Clearly, the latter is neither a necessary nor sufficient condition for second-moment stability.

\begin{figure*}[th]
    \centering
    \includegraphics[width=0.9\linewidth]{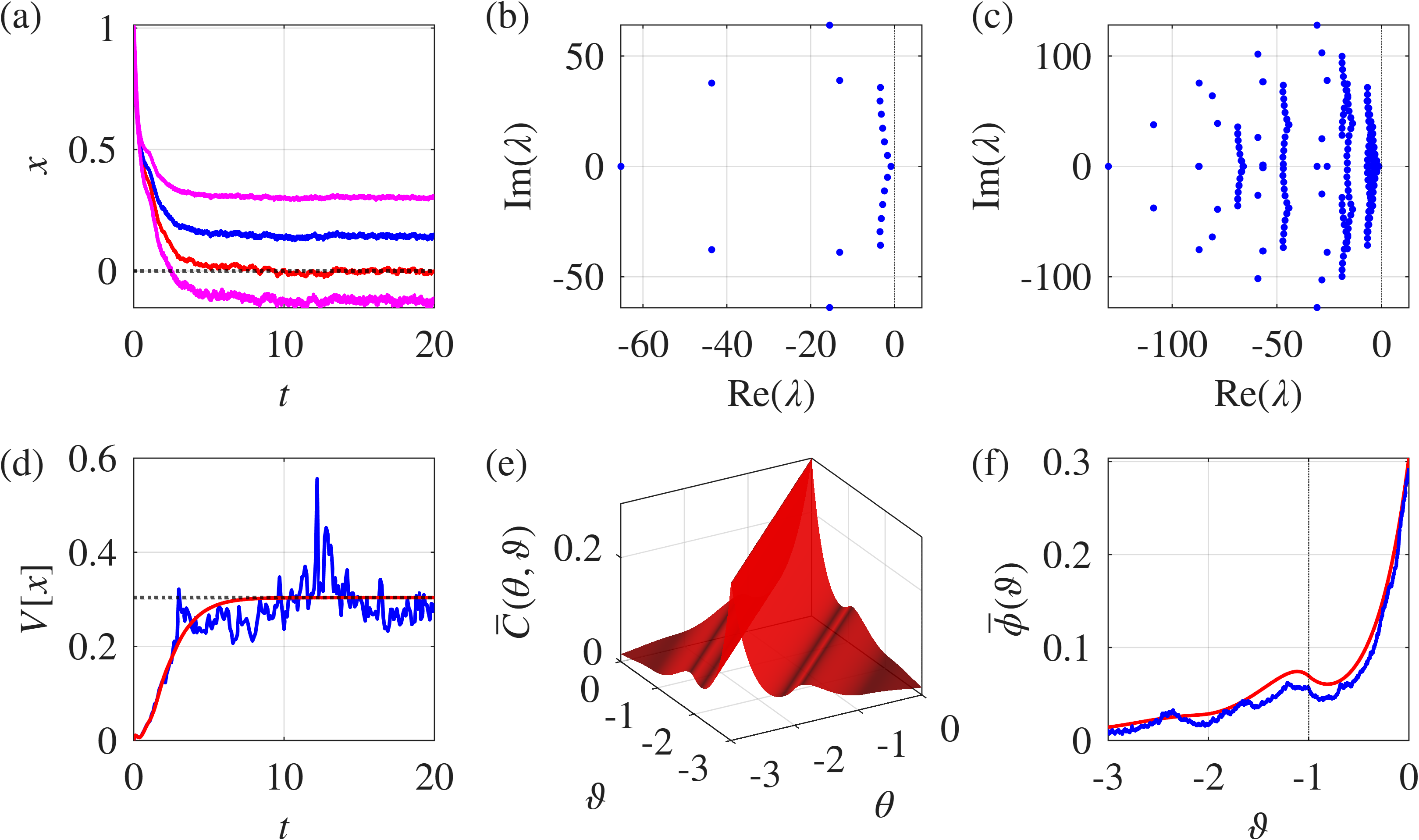}\\
    \caption{Numeric study for point A ($a=-3.0$, $b=1.0$) from Fig.~\ref{fig:scalar_Stab_comp}, for which Eq.~\eqref{eq:Ito_sdde} is predicted by the theory in this paper to be both first- and second-moment stable ($\rho(\mathbf{A})=-0.792$, $\rho(\mathbf{F})=-1.268$). (a) ensemble trajectory data from 5000 Euler-Maruyama simulations ($\Delta t = 0.001$), red line: mean, blue line: median, pink lines: 1st and 3rd quartile; (b) spectrum of $\mathbf{A}$ ($M=20$); (c) spectrum of $\mathbf{F}$ ($M=20$); (d) variance $C_t(0,0)-m_t^2(0)$, blue line: ensemble of stochastic trajectories, red line: direct integration of \eqref{eq:cov_disc}; (e) semi-analytic stationary correlation function $\bar{C}(\theta,\vartheta)$; (f) stationary correlation function $\bar{\phi}(\vartheta)$, blue line: ensemble of stochastic trajectories, red line: semi-analytic solution.
    \label{fig:scalar_stab_a}} 
\end{figure*}

\begin{figure*}[h]
    \centering
    \includegraphics[width=0.9\linewidth]{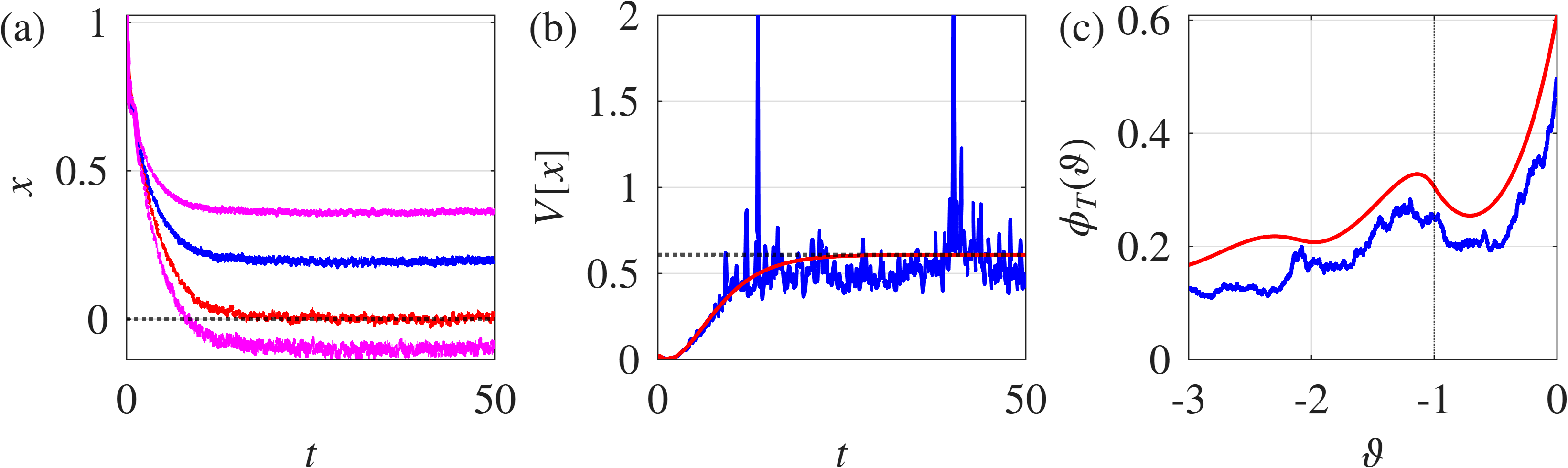}\\
    \caption{Numeric study of point B ($a=-3.0$, $b=2.1$) from Fig.\ref{fig:scalar_Stab_comp}, for which Eq.~\eqref{eq:Ito_sdde} is predicted by the theory in this paper to be both first- and second-moment stable ($\rho(\mathbf{A})=-0.264$, $\rho(\mathbf{F})=-0.303$). For legends, see the description of panels (a), (d), and (f) in the caption of Fig.\ref{fig:scalar_stab_a}.}
    \label{fig:scalar_stab_b} 
\end{figure*}

\begin{figure*}[h]
    \centering
    \includegraphics[width=0.6\linewidth]{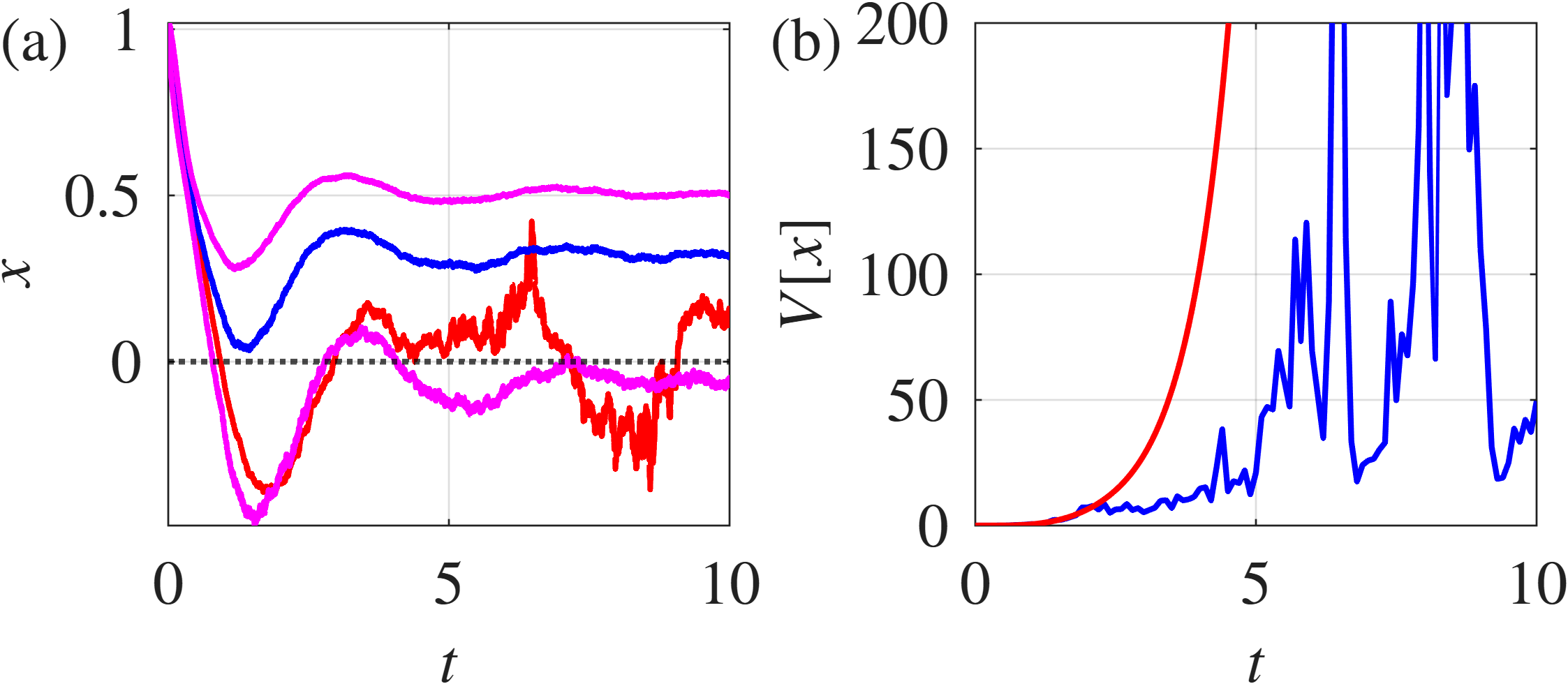}\\
    \caption{Numeric study of point C ($a=-0.4$, $b=-0.9$) from Fig.\ref{fig:scalar_Stab_comp}, for which Eq.~\eqref{eq:Ito_sdde} is predicted by the theory in this paper to be first- but not second-moment stable ($\rho(\mathbf{A})=-0.511$, $\rho(\mathbf{F})=1.288$). For legends, see the description of panels (a) and (d) in the caption of Fig.\ref{fig:scalar_stab_a}.}
    \label{fig:scalar_stab_c} 
\end{figure*}

\begin{figure*}[h]
    \centering
    \includegraphics[width=0.6\linewidth]{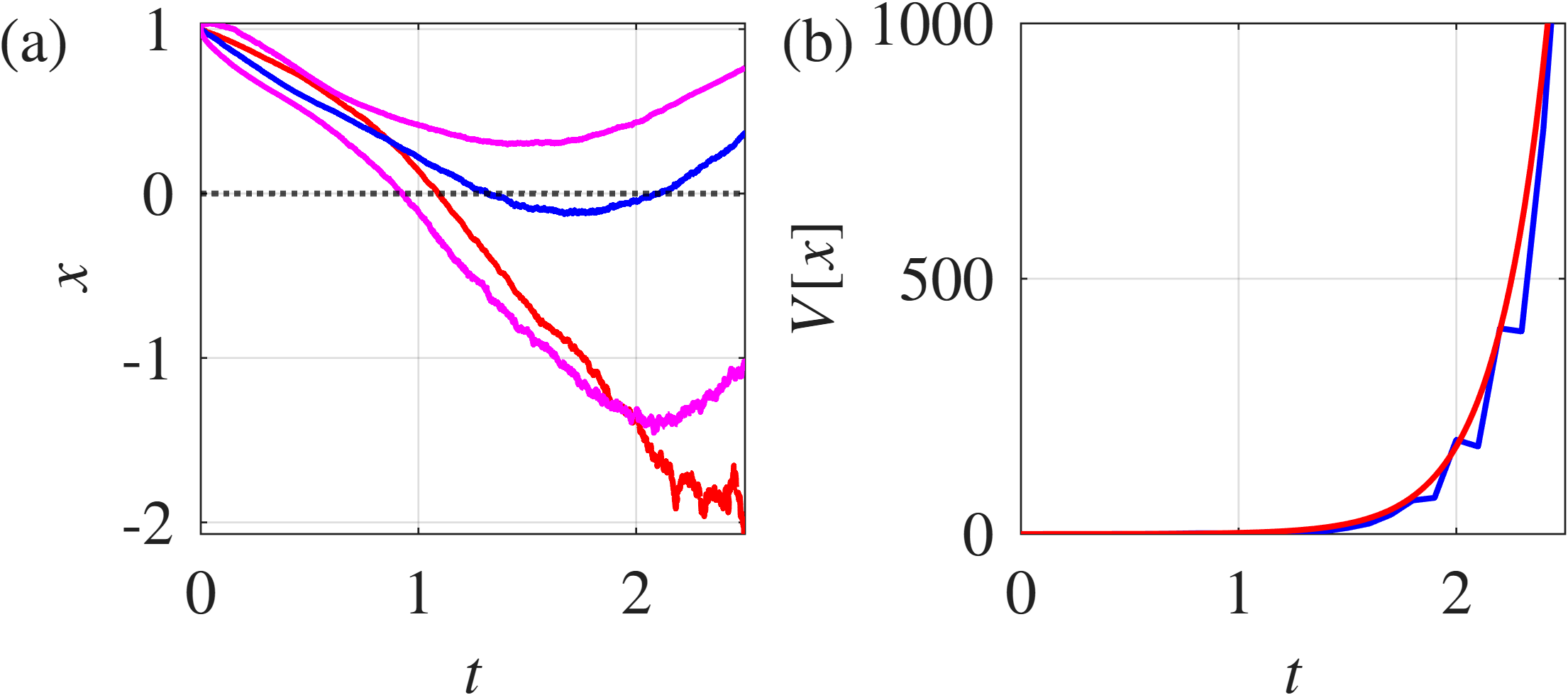}\\
    \caption{Numeric study of point D ($a=1.0$, $b=-1.5$) from Fig.\ref{fig:scalar_Stab_comp}, for which Eq.~\eqref{eq:Ito_sdde} is predicted by the theory in this paper to neither first- nor second-moment stable ($\rho(\mathbf{A})=0.273$, $\rho(\mathbf{F})=4.149$). For legends, see the description of panels (a) and (d) in the caption of Fig.\ref{fig:scalar_stab_a}.}
    \label{fig:scalar_stab_d} 
\end{figure*}

For the case that $\gamma=1$, Figs.~\ref{fig:scalar_stab_a} to \ref{fig:scalar_stab_d} validate these predictions and also provide more visual reference to the theoretical results of previous sections.

For point A in Fig.~\ref{fig:scalar_Stab_comp}, for which $a=-3$ and $b=1$ and both first-and second-moment stability is predicted by both the theory in this paper and that of \cite{mackey1995solution}, the panels of Fig.~\ref{fig:scalar_stab_a} show representations of (a) ensemble mean and quartile time signals for 5000 trajectories of Eq.~\eqref{eq:Ito_sdde} integrated using an Euler-Maruyama time discretization with a fixed stepsize of $\Delta t=0.001$ and a distribution of initial conditions with sample mean equal to $1$ and sample variance equal to $0$; (b)-(c) the spectra of $\mathbf{A}$ and $\mathbf{F}$ for the corresponding discretization with $M=20$, (d) the time evolution of the variance $C_t(0,0)-m_t^2(0)$ from integration of a discretized version of the covariance boundary-value problem starting with $C_0(\theta,\vartheta)= m_0(\theta)m_0(\vartheta)$ compared against an estimate from Monte-Carlo simulations; (e) the predicted stationary solution $\bar{C}(\theta,\vartheta)$ on an extension of the fundamental domain $[-3\tau,0]\times[-3\tau,0]$; (f) the predicted stationary solution $\bar{\phi}$ from the closed-form expression on $[-\tau,0]$ and integration of \eqref{eq:corr_stat_dde} on $[-3\tau, \, -\tau]$. Clearly, close agreement is demonstrated between brute-force integration and analysis of the covariance/correlation boundary-value problems.

Figure \ref{fig:scalar_stab_b} repeats this exercise but for point B in Fig.~\ref{fig:scalar_Stab_comp}, for which $a=-3$ and $b=2.1$ and second-moment stability is predicted only by the theory in this paper. The results closely resemble those in Fig.~\ref{fig:scalar_stab_a}, as expected. In contrast, for point C in Fig.~\ref{fig:scalar_Stab_comp}, for which $a=-0.4$ and $b=-0.9$ and second-moment stability is predicted only by the theory of \cite{mackey1995solution}, forward integration in Fig.~\ref{fig:scalar_stab_c} shows instability. In this case, we omit the predicted forms of $\bar{\phi}$ and $\bar{C}$, since the latter is no longer a valid covariance kernel. Analogous behavior is observed in Fig.~\ref{fig:scalar_stab_d} for point D in Fig.~\ref{fig:scalar_Stab_comp}, for which $a=1$ and $b=-1.5$ and both first- and second-moment instability are predicted by the theory in this paper as well as that of \cite{mackey1995solution}.

\subsection{An inverted pendulum}
We turn next to an example with $n=2$,  taken directly from \cite{fodor2023collocation}, for which
\begin{align}
    &a = \begin{pmatrix}
        0 & 1\\ k & 0
    \end{pmatrix},\, b = \begin{pmatrix}
        0 & 0\\ -p & -d
    \end{pmatrix},\,\gamma = \begin{pmatrix}
        0 \\ \sigma
    \end{pmatrix},
\end{align}
$\beta = -\sigma b$, and $\alpha=0$. For positive values of $k$, Eq.~\eqref{eq:Ito_sdde} then models the dynamics of an inverted pendulum near its equilibrium point, subjected to delayed and uncertain PD control, as well as additive external white noise. From the theory in Sec.~\ref{sec:analytical results}, we recall that the boundary of first-moment stability coincides with a subset of the zero-level set of $\det\,(\Psi)$ assuming that $\sigma=0$ and that second-moment stability is achieved within this region inside a boundary defined by $\det\,(\Psi)=0$ for positive $\sigma$. Panel (d) in Fig.~\ref{fig:osc_stab} shows the corresponding curves in the case that $k=5$, $\tau=0.3$, and $\sigma=0.1$, in close agreement with predictions in \cite{fodor2023collocation} obtained through a collocation-based discretization approach. If $\tau$ is allowed to increase, both regions of first- and second-moment stability shrink successively until that of second-moment stability disappears around $\tau=0.435$ and that of first-moment stability exactly at $\tau=\sqrt{2/k}=0.6325$, as shown in panel (e) of Fig.~\ref{fig:osc_stab}.

For a parameter combination securely within the predicted region of second-moment stability, panel (c) in Fig.~\ref{fig:osc_stab} compares the four components of the predicted stationary solution $\bar{\phi}$ against estimates obtained using Monte Carlo simulations of the original SDDE. The transient simulation results in panels (a) and (b) for the corresponding noisy position and velocity signals, along with their numerically estimated and analytically predicted ensemble statistical properties further serve to validate the theory.

This example also affords an opportunity to comment on the candidate boundaries implied by roots of either expression \eqref{eq:Buckwar} or \eqref{eq:Buckwar2} crossing the imaginary axis. Allowing for an arbitrary nonzero matrix $\alpha$, it is straightforward to show that the determinant in \eqref{eq:Buckwar} equals the determinant in \eqref{eq:Buckwar2} multiplied by the factor
\begin{equation}
    \frac{e^{\lambda\tau}\left(2\lambda-\det\,(\alpha)\right)}{d(1-\alpha_{11}\sigma)+e^{\lambda\tau}\left(2\lambda-\det\,(\alpha)\right)+p\alpha_{12}\sigma},
\end{equation}
which vanishes only when $\lambda=\det\,(\alpha)/2$. It follows that condition \eqref{eq:Buckwar} implies the existence of a root with positive real part for all parameter combinations if $\det\,(\alpha)>0$ and, if $\det\,(\alpha)<0$, roots with positive real part if and only if such roots exist for the determinant in \eqref{eq:Buckwar2}. For the parameter combination used in Fig.~\ref{fig:osc_stab} and $\alpha=0$, the boundary implied by \eqref{eq:Buckwar2} is shown as a dashed curve in panel (d). Clearly, this has no bearing on neither first- nor second-moment stability of the stochastic dynamics.

\subsection{Machine tool vibrations}
The final example in this paper, also with $n=2$, is a model of a regenerative turning process \cite{fodor2024efficient} given by
\begin{equation}
    m\ddot{x}(t) + c \dot{x}(t) + k x(t) = wK_c(h_0 - x(t) + x(t-T)), \label{eq:turning_model}
\end{equation}
where the cutting force coefficient $K_c=\bar{K}_{\mathrm{c}} + \sigma_{K_{\mathrm{c}}}\Gamma(t)$ has a stochastic component. Here, $m$, $c$, and $k$ are the modal mass, damping, and stiffness, respectively, measured at the tool tip, and the cutting force is linear in the chip width $w$ and the active chip thickness $h_0 - x(t) + x(t-T)$, defined in terms of the feed per revolution $h_0$ and the turning period $T = 2\pi/\omega$, where $\omega$ denotes the spindle speed. After non-dimensionalization, \eqref{eq:turning_model} reduces to \eqref{eq:Ito_sdde} with
\begin{equation}
    \begin{split}
    &a = \begin{pmatrix}
        0 & 1 \\  -1-\tilde{w} & -2\zeta
    \end{pmatrix}, \quad 
    b = \begin{pmatrix}
        0 & 0 \\ \tilde{w} & 0
    \end{pmatrix}, \\
    \alpha & = \begin{pmatrix}
        0 & 0 \\ -\sigma \tilde{w} & 0
    \end{pmatrix}, \quad 
    \beta = \begin{pmatrix}
        0 & 0 \\ \sigma \tilde{w} & 0
    \end{pmatrix}, \quad
    \gamma = \begin{pmatrix}
        0 \\ \sigma \tilde{w}
    \end{pmatrix}.
    \end{split}\label{eq:turningpars}
\end{equation}
Here, $\tilde{w} = \bar{K}_{\mathrm{c}}w/k$ is the equivalent chip width, and the modal damping ratio $\zeta = c/2m\omega_{\mathrm{n}}$, dimensionless time delay $\tau = \omega_{\mathrm{n}}T$, and rescaled noise coefficient $\sigma  = \sigma_{K_{\mathrm{c}}}/\bar{K}_{\mathrm{c}}\sqrt{\omega_n}$ are expressed in terms of the natural frequency of the machine tool assembly $\omega_{\mathrm{n}} = \sqrt{k/m}$. The non-dimensionalized spindle speed is given by $\tilde{\omega} = \omega/\omega_{\mathrm{n}}$.

When operating a turning machine, for example one whose dynamics is governed by \eqref{eq:turning_model}, machinists typically only have control over two parameters, namely the spindle speed $\omega$ and the chip width or axial feed $w$. A well-informed choice for these parameters is crucial for ensuring stable and reliable material removal without unwanted, harmful chatter vibrations \cite{munoa2016chatter}. Spindle speed vs chip width stability charts are regularly used by the manufacturing industry for this purpose. 

\begin{figure*}[ht]
    \centering
    \includegraphics[width=0.65\linewidth]{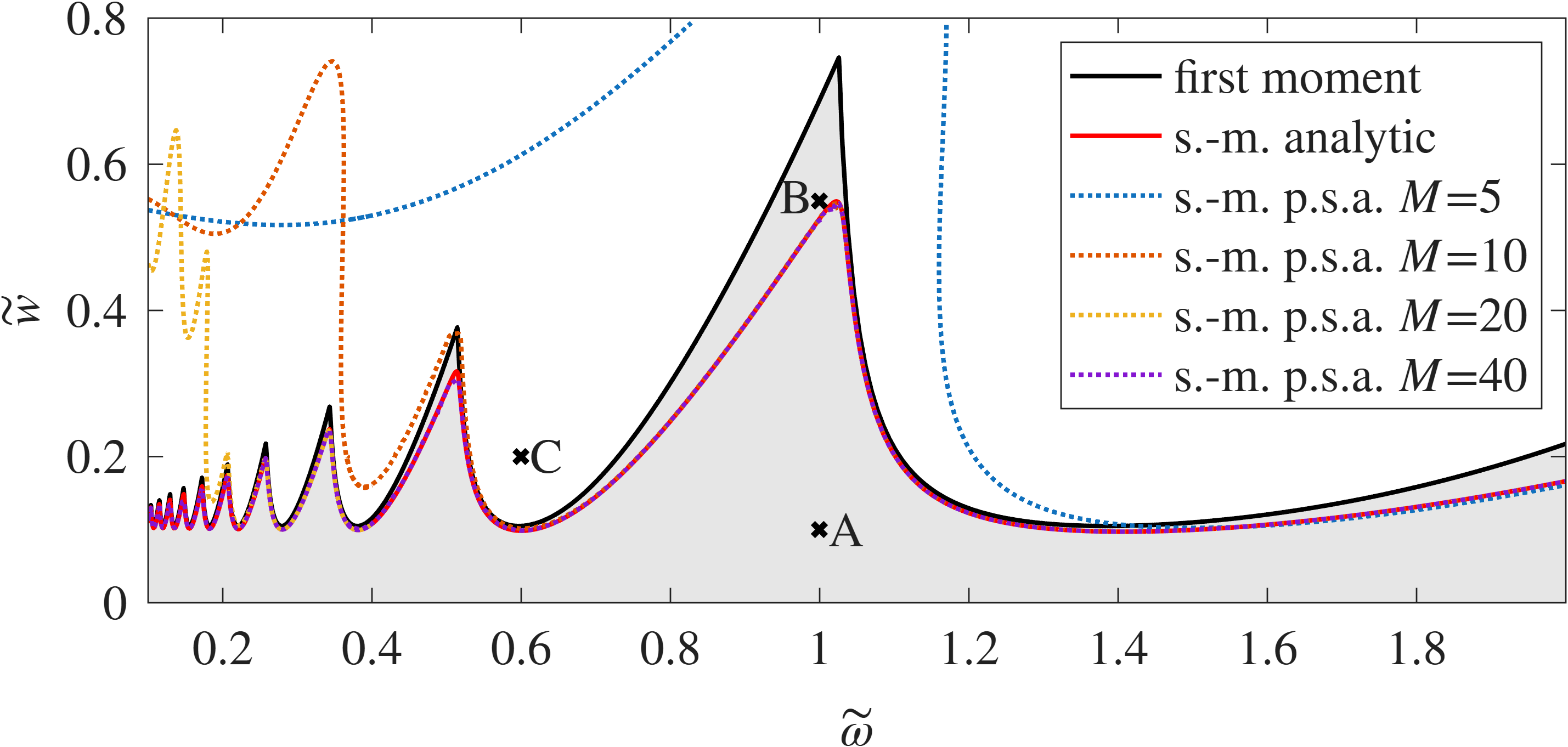}\\
    \vspace{0.01\textheight}
    \includegraphics[width=0.3\linewidth]{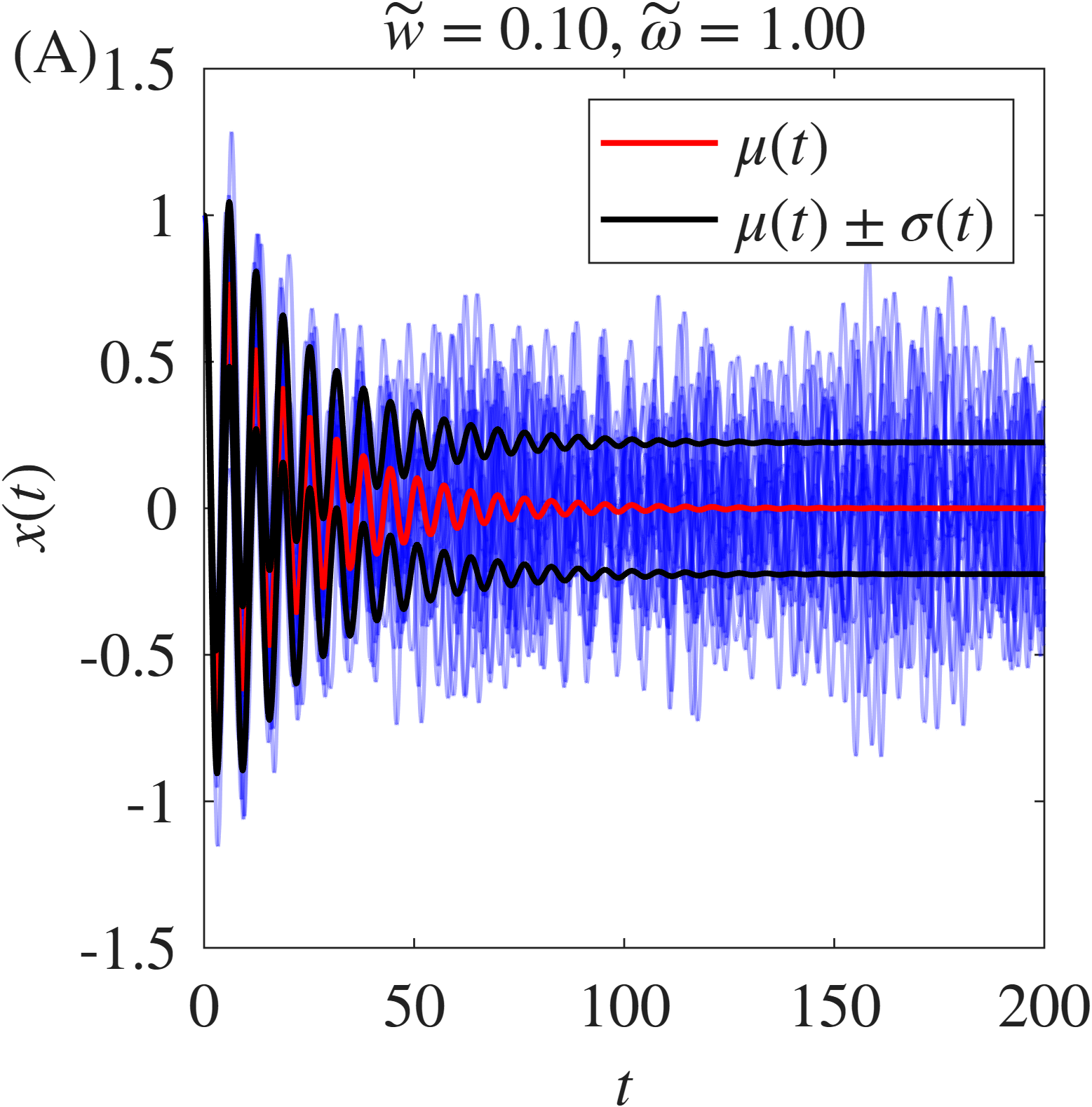}
    \includegraphics[width=0.3\linewidth]{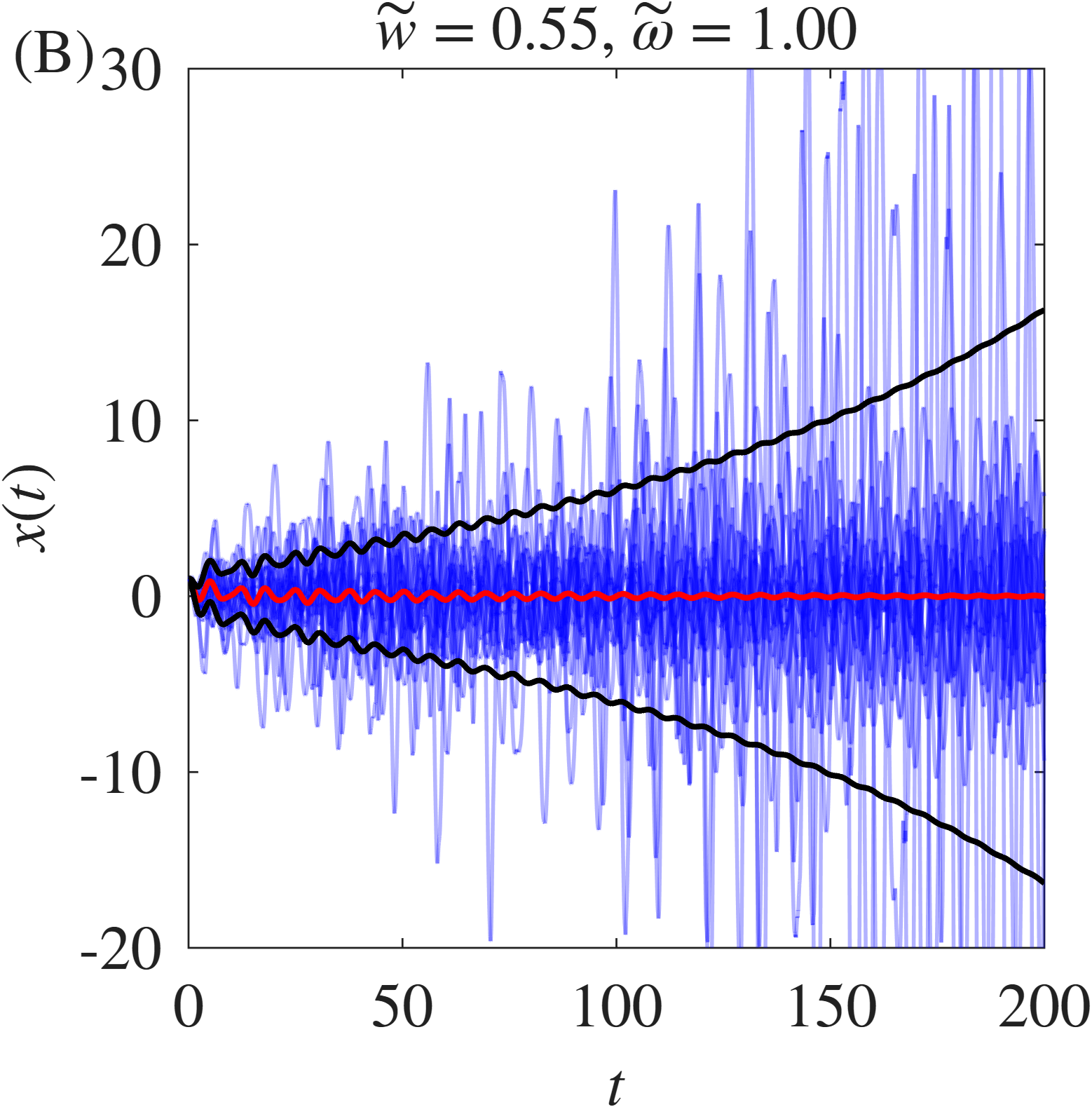}
    \includegraphics[width=0.3\linewidth]{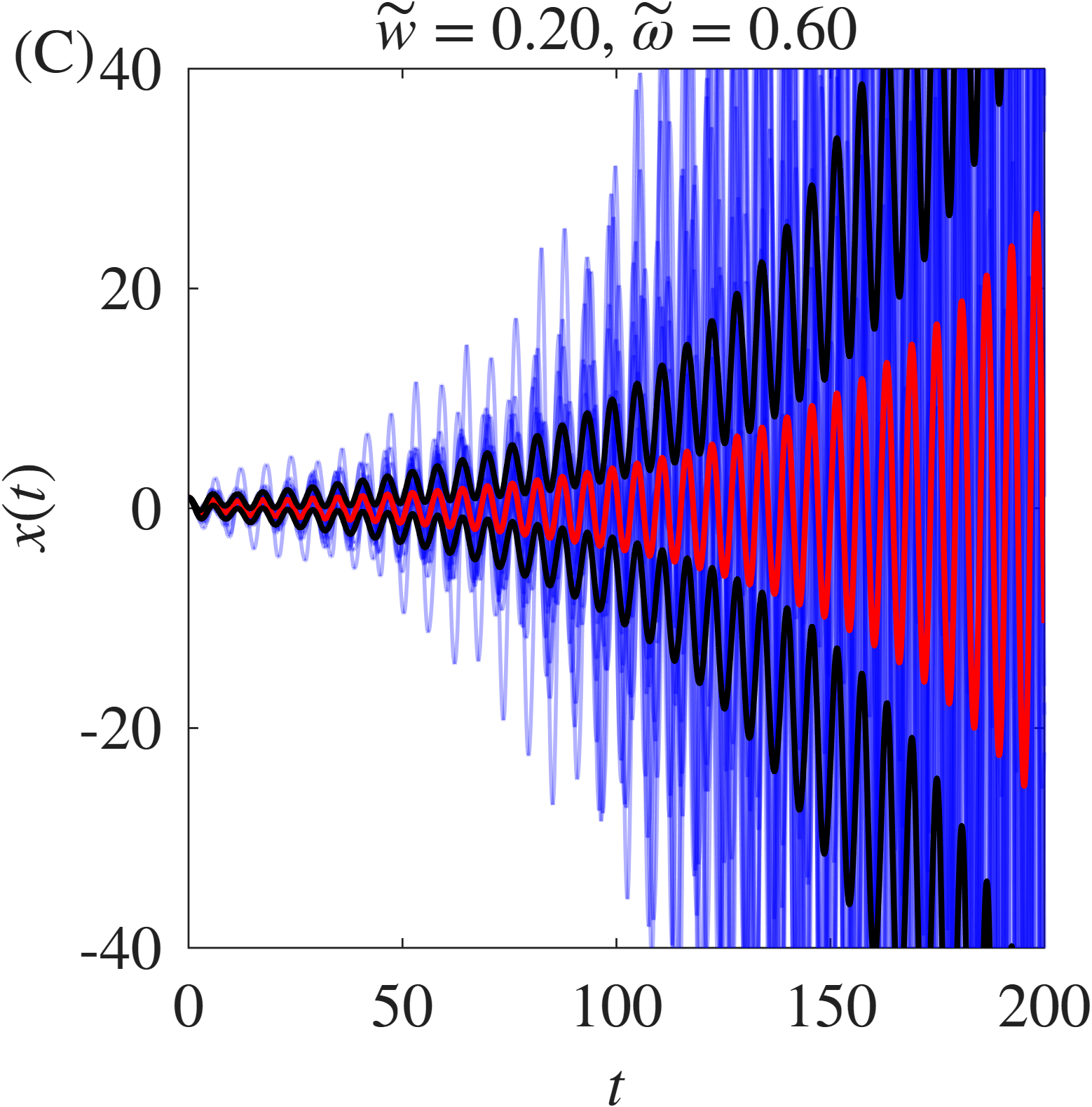}
    \caption{First and second-moment stability map and the results of stochastic integration for a single-degree-of-freedom mathematical model \eqref{eq:turning_model}-\eqref{eq:turningpars} of a turning process considering an uncertain regenerative cutting force ($\zeta=0.05$, $\sigma=1$). The lobe diagram shows a comparison of analytic and pseudo-spectral-approximation-based numeric stability predictions with different discretization resolutions $M$. Panels (A), (B), and (C) show a subset of simulated stochastic trajectories (blue) and time histories for the sample mean $\mu(t)$ (red) and deviations from the sample mean by one sample standard deviation (black) for the three labeled points in the lobe diagram.}
    \label{fig:turning_stab}
\end{figure*}

One such stability chart, a so-called lobe diagram, augmented with second-moment stability bounds is presented in Fig.~\ref{fig:turning_stab} for $\zeta=0.05$ and $\sigma=1$. Here the solid red curve was found by applying numerical continuation to the bordered matrix problem \eqref{eq:bordered}. For comparison, we similarly track zero crossings of the spectral abscissa of $\mathbf{F}$ for $M=[5,\,10,\,20,\,40]$. This shows that as $\tilde{\omega}$ gets smaller and the time delay grows larger, the convergence and accuracy of this numerical method degrades significantly. The stability bound found with $M=40$ is in good agreement with the solid red line, which suggests that the discretization-based stability predictions converge on the analytic limit derived in this paper as $M\rightarrow\infty$. For this example, the main advantage of the method proposed in this paper is its computational efficiency. On a personal computer, the generation of the second-moment stability maps using a pseudo-spectral discretization approach took 6 seconds with $M=5$, two minutes with $M=10$, 30 minutes with $M=20$, and roughly 22 hours with $M=40$. In comparison, tracking the $\det\,(\Psi)=0$ curve was accomplished in 10 seconds. The second-moment stability maps presented here are analogous with the numerical results presented in \cite{fodor2024efficient}, further substantiating the accuracy and reliability of the proposed methodology.

To help contextualize the implications of the stability chart in Fig.~\ref{fig:turning_stab} and to emphasize the importance of accounting for second-moment stability, the three lower panels show the result of Monte-Carlo simulations for the parameter pairs marked with capital letters in the lobe diagram. Panel (A) shows a desirable operation, where both the mean and variance settle to a stationary solution. In contrast, panel (B) displays an in-between scenario, which is first but not second-moment stable, while panel (C) illustrates an unstable process, where even the underlying deterministic DDE is unstable. The case of point B clearly demonstrates the importance of the insights gained from analyzing the effects of noise on the turning model \eqref{eq:turning_model}. If we were to neglect the effects of noise in the system, we would have deemed this parameter combination stable and thus desirable. But, as is clear from the individual trajectories, such an operation would lead to unacceptable levels of tool-tip oscillations.

\section{Conclusions}\label{sec:conclusions}

In this paper, we have shown that second-moment  stability boundaries for an arbitrary system of linear, stochastic delay-differential equations with a single constant delay, single noise source, and constant coefficients may be characterized without resorting to bounding estimates, Lyapunov theory, or various numerical approximations. A key enabling contribution is the discovery of a semi-analytic expression for the stationary solution to a delay-PDE boundary-value problem governing a two-variable correlation function. Using this expression, we have derived algebraic conditions that define second-moment stability boundaries and shown how these conditions can be solved numerically at significantly lower computational and memory costs, and without the scaling challenges of discretization-based approaches. Indeed, validation on several benchmark problems was able to reproduce approximate second-moment stability maps already available in the literature with ease and accuracy.

In contrast to ``discretization-first'' approaches to the analysis of (S)DDEs, in which time histories are discretized \textit{a priori} and high-dimensional systems of (S)ODEs approximate the system dynamics, our approach has been motivated by a preference to defer discretization until it can no longer be avoided. We have been rewarded with results that entirely avoid discretization and that we apply to the infinite-dimensional dynamics of the original SDDE. It is appropriate, however, to stress that several results in this paper are motivated, but not rigorously justified, by observations for finite-dimensional discretizations. Were such observations not to generalize to the original SDDEs, this would put into question the applicability of discretization-based approaches in the first place. We have great confidence that convergence results for DDEs \cite{breda2014stability} also carry over to SDDEs and have proceeded accordingly.

Although the novel analytical results presented in this paper were derived for a somewhat restricted class of SDDEs, we expect substantial parts of the developed theory to generalize beyond the form in Eq.~\eqref{eq:Ito_sdde}. For example, the formulas found for the stationary correlation functions and second-moment stability boundaries should carry over with minimal modification to stochastic systems with multiple, mutually independent sources of white noise. Apart from a few extra terms arising in the boundary conditions of the advection-PDEs, the analysis is unchanged, and thus the same steps may be taken to characterize the second-moment dynamics. Similarly, it should be straightforward to account for colored noise by including the linear or linearized form of the corresponding filtering equations in the system matrices.

Generalization to SDDEs with multiple or distributed delays, however, appears far less straightforward. Here, the correlation boundary-value problem takes on a considerably more complicated form with no clear analytical path forward. At present, we can offer no alternative to discretization, for example of the corresponding advection-type covariance boundary-value problem with suitably modified boundary conditions in terms of a system of surrogate covariance ODEs, albeit at significantly greater computational cost than afforded us by the semi-analytic treatment in this paper. We leave to future work the question of whether the characterization of stability boundaries in terms of a vanishing determinant also applies in this case, although this seems likely. Inspired by observations in this paper regarding the relative computational costs of analyzing the matrices $\mathbf{F}$ and $\mathbf{\Xi}$, we suspect that savings can also be achieved in the case of multiple or distributed delays by focusing on the corresponding surrogate correlation DDEs.

Finally, it is worth commenting on the possibility of relaxing the restriction to constant coefficient matrices, e.g., to periodic functions of time. While neither the covariance nor the correlation boundary-value problems change as a result, we again see no clear analytical path forward other than through discretization. Unfortunately, the scaling challenges that we observed already in the study of stationary solutions may put such an analysis for high-dimensional systems out of reach of all but specialized high-performance computing resources. Further work should be devoted to this issue.

\bibliographystyle{unsrt} 
\bibliography{cas-refs}

@article{mackey1995solution,
  title={Solution moment stability in stochastic differential delay equations},
  author={Mackey, Michael C and Nechaeva, Irina G},
  journal={Physical Review E},
  volume={52},
  number={4},
  pages={3366},
  year={1995},
  publisher={APS}
}

@article{torkamani2014numerical,
  title={Numerical stability analysis of linear stochastic delay differential equations using Chebyshev spectral continuous time approximation},
  author={Torkamani, Shahab and Samiei, Ehsan and Bobrenkov, Oleg and Butcher, Eric A},
  journal={International Journal of Dynamics and Control},
  volume={2},
  number={2},
  pages={210--220},
  year={2014},
  publisher={Springer}
}

@book{trefethen2000spectral,
  title={Spectral methods in MATLAB},
  author={Trefethen, Lloyd N},
  year={2000},
  publisher={SIAM}
}

@book{baxendale2007stochastic,
  title={Stochastic differential equations: theory and applications},
  author={Baxendale, Peter H and Lototsky, Sergey V},
  volume={2},
  year={2007},
  publisher={World Scientific}
}

@article{sykora2019stochastic,
  title={Stochastic semi-discretization for linear stochastic delay differential equations},
  author={Sykora, Henrik T and Bachrathy, Daniel and Stepan, Gabor},
  journal={International Journal for Numerical Methods in Engineering},
  volume={119},
  number={9},
  pages={879--898},
  year={2019},
  publisher={Wiley Online Library}
}

@article{fodor2023collocation,
  title={Collocation method for stochastic delay differential equations},
  author={Fodor, Gerg{\H{o}} and Sykora, Henrik T and Bachrathy, D{\'a}niel},
  journal={Probabilistic Engineering Mechanics},
  volume={74},
  pages={103515},
  year={2023},
  publisher={Elsevier}
}

@book{insperger2011semi,
  title={Semi-discretization for time-delay systems: stability and engineering applications},
  author={Insperger, Tam{\'a}s and St{\'e}p{\'a}n, G{\'a}bor},
  volume={178},
  year={2011},
  publisher={Springer Science \& Business Media}
}

@misc{orosz2010traffic,
  title={Traffic jams: dynamics and control},
  author={Orosz, G{\'a}bor and Wilson, R Eddie and St{\'e}p{\'a}n, G{\'a}bor},
  journal={Philosophical Transactions of the Royal Society A: Mathematical, Physical and Engineering Sciences},
  volume={368},
  number={1928},
  pages={4455--4479},
  year={2010},
  publisher={The Royal Society Publishing}
}

@article{munoa2016chatter,
  title={Chatter suppression techniques in metal cutting},
  author={Munoa, Jokin and Beudaert, Xavier and Dombovari, Z and Altintas, Yusuf and Budak, Erhan and Brecher, Christian and Stepan, Gabor},
  journal={CIRP annals},
  volume={65},
  number={2},
  pages={785--808},
  year={2016},
  publisher={Elsevier}
}

@article{glass2021nonlinear,
  title={Nonlinear delay differential equations and their application to modeling biological network motifs},
  author={Glass, David S and Jin, Xiaofan and Riedel-Kruse, Ingmar H},
  journal={Nature communications},
  volume={12},
  number={1},
  pages={1788},
  year={2021},
  publisher={Nature Publishing Group UK London}
}

@article{michiels2010control,
  title={Control design for time-delay systems based on quasi-direct pole placement},
  author={Michiels, Wim and Vyhl{\'\i}dal, Tom{\'a}{\v{s}} and Z{\'\i}tek, Pavel},
  journal={Journal of Process Control},
  volume={20},
  number={3},
  pages={337--343},
  year={2010},
  publisher={Elsevier}
}

@article{fodor2020stochastic,
  title={Stochastic modeling of the cutting force in turning processes},
  author={Fodor, Gerg{\H{o}} and Sykora, Henrik T and Bachrathy, Daniel},
  journal={The International Journal of Advanced Manufacturing Technology},
  volume={111},
  number={1},
  pages={213--226},
  year={2020},
  publisher={Springer}
}

@article{sykora2020moment,
  title={On the moment dynamics of stochastically delayed linear control systems},
  author={Sykora, Henrik T and Sadeghpour, Mehdi and Ge, Jin I and Bachrathy, D{\'a}niel and Orosz, G{\'a}bor},
  journal={International Journal of Robust and Nonlinear Control},
  volume={30},
  number={18},
  pages={8074--8097},
  year={2020},
  publisher={Wiley Online Library}
}

@article{babasola2023stochastic,
  title={Stochastic delay differential equations: a comprehensive approach for understanding biosystems with application to disease modelling},
  author={Babasola, Oluwatosin and Omondi, Evans Otieno and Oshinubi, Kayode and Imbusi, Nancy Matendechere},
  journal={AppliedMath},
  volume={3},
  number={4},
  pages={702--721},
  year={2023},
  publisher={MDPI}
}

@article{yang2022stochastic,
  title={Stochastic analysis of a time-delayed viscoelastic energy harvester subjected to narrow-band noise},
  author={Yang, Yong-Ge and He, Li-Li and Zeng, Yuan-Hui and Sun, Ya-Hui and Xu, Wei},
  journal={International Journal of Non-Linear Mechanics},
  volume={147},
  pages={104230},
  year={2022},
  publisher={Elsevier}
}

@article{frank2001stationary,
  title={Stationary solutions of linear stochastic delay differential equations: Applications to biological systems},
  author={Frank, TD and Beek, PJ},
  journal={Physical Review E},
  volume={64},
  number={2},
  pages={021917},
  year={2001},
  publisher={APS}
}

@article{rene2017mean,
  title={Mean, covariance, and effective dimension of stochastic distributed delay dynamics},
  author={Ren{\'e}, Alexandre and Longtin, Andr{\'e}},
  journal={Chaos: An Interdisciplinary Journal of Nonlinear Science},
  volume={27},
  number={11},
  year={2017},
  publisher={AIP Publishing}
}

@article{ohira2000delayed,
  title={Delayed stochastic systems},
  author={Ohira, Toru and Yamane, Toshiyuki},
  journal={Physical Review E},
  volume={61},
  number={2},
  pages={1247},
  year={2000},
  publisher={APS}
}

@article{kuske2010competition,
  title={Competition of noise sources in systems with delay: the role of multiple time scales},
  author={Kuske, R},
  journal={Journal of Vibration and Control},
  volume={16},
  number={7-8},
  pages={983--1003},
  year={2010},
  publisher={SAGE Publications Sage UK: London, England}
}

@article{wang2014moment,
  title={Moment boundedness of linear stochastic delay differential equations with distributed delay},
  author={Wang, Zhen and Li, Xiong and Lei, Jinzhi},
  journal={Stochastic Processes and their Applications},
  volume={124},
  number={1},
  pages={586--612},
  year={2014},
  publisher={Elsevier}
}

@article{rackauckas2017differentialequations,
  title={Differentialequations. jl--a performant and feature-rich ecosystem for solving differential equations in julia},
  author={Rackauckas, Christopher and Nie, Qing},
  year={2017}
}

@article{higham2001algorithmic,
  title={An algorithmic introduction to numerical simulation of stochastic differential equations},
  author={Higham, Desmond J},
  journal={SIAM review},
  volume={43},
  number={3},
  pages={525--546},
  year={2001},
  publisher={SIAM}
}

@incollection{oksendal2003stochastic,
  title={Stochastic differential equations},
  author={{\O}ksendal, Bernt},
  booktitle={Stochastic differential equations: an introduction with applications},
  pages={38--50},
  year={2003},
  publisher={Springer}
}

@book{sun2006stochastic,
  title={Stochastic dynamics and control},
  author={Sun, Jian-Qiao},
  volume={4},
  year={2006},
  publisher={Elsevier}
}

@article{shaikhet2021one,
  title={About one method of stability investigation for nonlinear stochastic delay differential equations},
  author={Shaikhet, Leonid},
  journal={International Journal of Robust and Nonlinear Control},
  volume={31},
  number={8},
  pages={2946--2959},
  year={2021},
  publisher={Wiley Online Library}
}

@article{tuncc2019asymptotic,
  title={On the asymptotic stability of solutions of stochastic differential delay equations of second order},
  author={Tun{\c{c}}, Osman and Tun{\c{c}}, Cemil},
  journal={Journal of Taibah University for Science},
  volume={13},
  number={1},
  pages={875--882},
  year={2019},
  publisher={Taylor \& Francis}
}

@article{mao1992robustness,
  title={Robustness of stability of nonlinear systems with stochastic delay perturbations},
  author={Mao, Xuerong},
  journal={Systems \& control letters},
  volume={19},
  number={5},
  pages={391--400},
  year={1992},
  publisher={Elsevier}
}

@article{samiei2013lyapunov,
  title={On Lyapunov stability of scalar stochastic time-delayed systems},
  author={Samiei, Ehsan and Torkamani, Shahab and Butcher, Eric A},
  journal={International Journal of Dynamics and Control},
  volume={1},
  number={1},
  pages={64--80},
  year={2013},
  publisher={Springer}
}

@article{novivcenko2012phase,
  title={Phase reduction of weakly perturbed limit cycle oscillations in time-delay systems},
  author={Novi{\v{c}}enko, V and Pyragas, K},
  journal={Physica D: Nonlinear Phenomena},
  volume={241},
  number={12},
  pages={1090--1098},
  year={2012},
  publisher={Elsevier}
}

@book{dankowicz2013recipes,
  title={Recipes for continuation},
  author={Dankowicz, Harry and Schilder, Frank},
  year={2013},
  publisher={SIAM}
}

@article{hayes1950roots,
  title={Roots of the transcendental equation associated with a certain difference-differential equation},
  author={Hayes, N.D.},
  journal={Journal of the London Mathematical Society},
  volume={1},
  number={3},
  pages={226--232},
  year={1950},
  publisher={Wiley Online Library}
}

@book{mao2007stochastic,
  title={Stochastic differential equations and applications},
  author={Mao, Xuerong},
  year={2007},
  publisher={Elsevier}
}

@book{breda2014stability,
  title={Stability of linear delay differential equations: A numerical approach with MATLAB},
  author={Breda, Dimitri and Maset, Stefano and Vermiglio, Rossana},
  year={2014},
  publisher={Springer}
}

@article{fodor2024efficient,
  title={Efficient approximation of stochastic turning process based on power spectral density},
  author={Fodor, Gerg{\H{o}} and Bachrathy, D{\'a}niel},
  journal={The International Journal of Advanced Manufacturing Technology},
  volume={133},
  number={11},
  pages={5673--5681},
  year={2024},
  publisher={Springer}
}

@book{hale2013introduction,
  title={Introduction to functional differential equations},
  author={Hale, Jack K and Lunel, Sjoerd M Verduyn},
  volume={99},
  year={2013},
  publisher={Springer Science \& Business Media}
}

@article{higham2000mean,
  title={Mean-square and asymptotic stability of the stochastic theta method},
  author={Higham, Desmond J},
  journal={SIAM journal on numerical analysis},
  volume={38},
  number={3},
  pages={753--769},
  year={2000},
  publisher={SIAM}
}

@article{nandanoori2018mean,
  title={Mean square stability analysis of stochastic continuous-time linear networked systems},
  author={Nandanoori, Sai Pushpak and Diwadkar, Amit and Vaidya, Umesh},
  journal={IEEE Transactions on Automatic Control},
  volume={63},
  number={12},
  pages={4323--4330},
  year={2018},
  publisher={IEEE}
}

@book{williams2006gaussian,
  title={Gaussian processes for machine learning},
  author={Williams, Christopher KI and Rasmussen, Carl Edward},
  volume={2},
  number={3},
  year={2006},
  publisher={MIT press Cambridge, MA}
}

@book{milton2021mathematics,
  title={Mathematics as a laboratory tool},
  author={Milton, John and Ohira, Toru},
  year={2021},
  publisher={Springer}
}

@incollection{doering2018modeling,
  title={Modeling complex systems: stochastic processes, stochastic differential equations, and Fokker-Planck equations},
  author={Doering, Charles R},
  booktitle={1990 Lectures in Complex Systems},
  pages={3--52},
  year={2018},
  publisher={CRC Press}
}

@book{ghanem2003stochastic,
  title={Stochastic finite elements: a spectral approach},
  author={Ghanem, Roger G and Spanos, Pol D},
  year={2003},
  publisher={Courier Corporation}
}

@article{buckwar2013note,
  title={A NOTE ON THE ANALYSIS OF ASYMPTOTIC MEAN-SQUARE STABILITY PROPERTIES FOR SYSTEMS OF LINEAR STOCHASTIC DELAY DIFFERENTIAL EQUATIONS.},
  author={Buckwar, Evelyn and Notarangelo, Girolama},
  journal={Discrete \& Continuous Dynamical Systems-Series B},
  volume={18},
  number={6},
  year={2013}
}

@article{erban2016cucker,
  title={A Cucker--Smale model with noise and delay},
  author={Erban, Radek and Haskovec, Jan and Sun, Yongzheng},
  journal={SIAM Journal on Applied Mathematics},
  volume={76},
  number={4},
  pages={1535--1557},
  year={2016},
  publisher={SIAM}
}

@article{havskovec2022asymptotic,
  title={Asymptotic and exponential decay in mean square for delay geometric Brownian motion},
  author={Ha{\v{s}}kovec, Jan},
  journal={Applications of Mathematics},
  volume={67},
  number={4},
  pages={471--483},
  year={2022},
  publisher={Springer}
}

@article{appleby2009geometric,
  title={Geometric Brownian motion with delay: Mean square characterisation},
  author={Appleby, John and Mao, Xuerong and Riedle, Markus},
  journal={Proceedings of the American Mathematical Society},
  volume={137},
  number={1},
  pages={339--348},
  year={2009}
}

@article{Schneider65,
    title={Positive operators and an inertia theorem},
    author={Schneider, Hans},
    journl={Numerische Mathematik},
    volume={7},
    pages={11–17.},
    year={1965}
}

@article{dankowicz2020multidimensional,
  title={Multidimensional manifold continuation for adaptive boundary-value problems},
  author={Dankowicz, Harry and Wang, Yuqing and Schilder, Frank and Henderson, Michael E},
  journal={Journal of Computational and Nonlinear Dynamics},
  volume={15},
  number={5},
  pages={051002},
  year={2020},
  publisher={American Society of Mechanical Engineers}
}

\appendix

\section{Discretization of the correlation boundary-value problem}\label{app:num_corr}

We considered in Sec.~\ref{sec:discretization} the discretization of the covariance boundary-value problem in terms of the coefficient matrix $\mathbf{F}$, as this yielded insights regarding second-moment stability that we extrapolated to the infinite-dimensional case. In this appendix, we derive the discretization of the correlation boundary-value problem, partly for completeness and partly as it affords us a computationally more efficient condition for approximate second-moment stability boundaries than the vanishing of $\det\,(\mathbf{F})$.

To discretize the correlation boundary-value problem in Eqs.~\eqref{eq:corr_pde} and \eqref{eq:corr_pde_bc}, consider, again, the mesh $-\tau=s_1<\cdots<s_M=0$ from Sec.~\ref{sec:discretization}. Define $\boldsymbol{\varphi}(t)$ as the vectorization of the block matrix
\begin{equation}
    \boldsymbol{\phi}_t=\begin{bmatrix}
        \phi_t(s_1) &\cdots & \phi_t(s_M)
    \end{bmatrix}.
\end{equation}
We use Lagrange interpolation to model the values of $\phi_{t+s_i}(-\tau-s_i)$ as these are not assumed to be contained in $\boldsymbol{\phi}_{t+s_i}$ other than for $i=1$ and $i=M$. Specifically, let $\boldsymbol{\ell}(\sigma)$ denote a row vector of Lagrange interpolation coefficients corresponding to the query point $\sigma$ on the given mesh. Then, the matrix
\begin{equation}
    P_n \left(\boldsymbol{\ell}(-\tau-s_i) \otimes I_{n^2}\right)\boldsymbol{\varphi}(t+s_i).
\end{equation}
represent the vectorization of
$\phi^\mathsf{T}_{t+s_i}(-\tau-s_i)$. For $i=1$ and $i=M$, this reduces exactly to the vectorizations of $\phi_{t-\tau}^\mathsf{T}(0) = \phi_{t-\tau}(0)$ and $\phi_t^\mathsf{T}(-\tau)$, respectively. It follows that the vectorization of the block matrix
\begin{align}
        &\begin{bmatrix}
        b\phi_{t+s_1}^\mathsf{T}(-\tau-s_1) &\cdots & b\phi_{t+s_{M}}^\mathsf{T}(-\tau-s_{M})
    \end{bmatrix}\nonumber\\&\qquad+\begin{bmatrix}
        \mathbf{0}_{1\times n(M-1)} & \phi_t(-\tau)b^\mathsf{T}
    \end{bmatrix}
\end{align}
is given by
\begin{equation}
        \begin{bmatrix}
        (I_n\otimes b)P_n \left(\boldsymbol{\ell}(-\tau-s_1) \otimes I_{n^2}\right)\boldsymbol{\varphi}(t+s_1)\\\vdots\\(I_n\otimes b)P_n \left(\boldsymbol{\ell}(-\tau-s_{M-1}) \otimes I_{n^2}\right)\boldsymbol{\varphi}(t+s_{M-1})\\\left(\mathbf{e}_1^\mathsf{T}\otimes\left((I_n\otimes b)P_n +(b\otimes I_n)\right)\right)\boldsymbol{\varphi}(t)
    \end{bmatrix}.\label{eq:gammas}
\end{equation}

Similarly, the vectorizations of the block matrices
\begin{align}
    &\begin{bmatrix}
        a\phi_t(s_1) & \cdots & a\phi_t(s_M)
        \end{bmatrix}+\begin{bmatrix}\mathbf{0}_{1\times n(M-1)} & \phi_t(0)a^\mathsf{T}
    \end{bmatrix},
\end{align}
\begin{equation}
    \begin{bmatrix}
        \mathbf{0}_{1\times n(M-1)} & \alpha \phi_t(-\tau)\beta^\mathsf{T}+\beta \phi_t^\mathsf{T}(-\tau)\alpha^\mathsf{T}
    \end{bmatrix},
\end{equation}
and
\begin{equation}
    \begin{bmatrix}
        \mathbf{0}_{n\times n(M-1)} & \alpha \phi_t(0)\alpha^\mathsf{T}+\beta \phi_{t-\tau}(0)\beta^\mathsf{T}
    \end{bmatrix}
\end{equation}
are given by
\begin{equation}
    \left(\left(\boldsymbol{I}_{Mn}\otimes a\right)+\left(\mathbf{e}_M\mathbf{e}_M^\mathsf{T}\right)\otimes(a\otimes I_n)\right)\boldsymbol{\varphi}(t),\label{eq:phicoeff1}
\end{equation}
\begin{align}
    \left(\left(\mathbf{e}_M\mathbf{e}_1^\mathsf{T}\right)\otimes\left((\beta\otimes\alpha)+(\alpha\otimes\beta)P_n\right)\right)\boldsymbol{\varphi}(t),\label{eq:phicoeff2}
\end{align}
and
\begin{align}
    &\left(\left(\mathbf{e}_M\mathbf{e}_M^\mathsf{T}\right)\otimes(\alpha\otimes\alpha)\right)\boldsymbol{\varphi}(t)\nonumber\\
    &\qquad+\left(\left(\mathbf{e}_M\mathbf{e}_M^\mathsf{T}\right)\otimes (\beta\otimes\beta) \right)\boldsymbol{\varphi}(t+s_1),\label{eq:gammaonemore}
\end{align}
respectively.

We refer to the equation
\begin{equation}
    \dot{\boldsymbol{\varphi}}(t)=\boldsymbol{\Lambda}\boldsymbol{\varphi}(t)+\sum_{i=1}^{M-1}\boldsymbol{\Gamma}_i\boldsymbol{\varphi}(t+s_i)+\boldsymbol{\gamma}(t)\label{eq:corr_disc_eq}
\end{equation}
for suitably constructed matrices $\boldsymbol{\Lambda}$, $\boldsymbol{\Gamma}_i$, $i=1,\ldots, M-1$, and $\boldsymbol{\gamma}(t)$ as the \textit{surrogate correlation DDE} associated with the correlation boundary-value problem. For example, $\boldsymbol{\gamma}(t)$ is the vectorization of
\begin{equation}
\begin{split}
    \mathbf{e}_{M}^\mathsf{T}\otimes\left(\boldsymbol{\eta}^\mathsf{T}\mathbf{\mathbf{m}}(t)\gamma^\mathsf{T} + \gamma \mathbf{m}^\mathsf{T}(t)\boldsymbol{\eta} + \gamma\gamma^\mathsf{T}\right)
    \end{split},\label{eq:cov_disc_const}
\end{equation}
where $\mathbf{m}(t)$, $\boldsymbol{\eta}$, and $\mathbf{e}_i$ were all defined in Sec.~\ref{sec:discretization}. Moreover, it follows from \eqref{eq:gammas} and \eqref{eq:gammaonemore} that
\begin{align}
    \boldsymbol{\Gamma}_1&= \mathbf{e}_1\otimes\left((I_n\otimes b)P_n \left(\boldsymbol{\ell}(-\tau-s_1) \otimes I_{n^2}\right)\right)\nonumber\\&\qquad+\left(\mathbf{e}_M\mathbf{e}_M^\mathsf{T}\right)\otimes (\beta\otimes\beta)
\end{align}
while, for $i\ne 1$,
\begin{equation}
    \boldsymbol{\Gamma}_i = \mathbf{e}_i\otimes\left((I_n\otimes b)P_n \left(\boldsymbol{\ell}(-\tau-s_i) \otimes I_{n^2}\right)\right).
\end{equation} 
Finally, from \eqref{eq:gammas}, \eqref{eq:phicoeff1}-\eqref{eq:gammaonemore}, we conclude that
\begin{align}
    \boldsymbol{\Lambda}&=\begin{bmatrix}\hat{\mathbf{D}}_M \\ \mathbf{0}_{1 \times M} \end{bmatrix} \otimes I_{n^2} + \boldsymbol{I}_{Mn}\otimes a+\left(\mathbf{e}_M\mathbf{e}_M^\mathsf{T}\right)\otimes(a\otimes I_n)\nonumber\\
    &\qquad+\left(\mathbf{e}_M\mathbf{e}_1^\mathsf{T}\right)\otimes\left((I_n\otimes b)P_n +(b\otimes I_n)\right)\nonumber\\
    &\qquad+\left(\mathbf{e}_M\mathbf{e}_1^\mathsf{T}\right)\otimes\left((\beta\otimes\alpha)+(\alpha\otimes\beta)P_n\right)\nonumber\\
    &\qquad+\left(\mathbf{e}_M\mathbf{e}_M^\mathsf{T}\right)\otimes(\alpha\otimes\alpha),
\end{align}
where $\mathbf{D}_M$ is the approximate differentiation matrix from Sec.~\ref{sec:discretization} and the accent $\hat{\;}$ again indicates omission of the last row of the matrix.

It is evident from the surrogate correlation DDE that a unique stationary solution is obtained when $\mathbf{m}\equiv 0$ if and only if $\det\,(\boldsymbol{\Xi})\ne 0$, where $\boldsymbol{\Xi}=\boldsymbol{\Lambda}+\sum_{i=1}^{M-1}\boldsymbol{\Gamma}_i$. Since $\mathbf{F}\in\mathbb{R}^{M^2n^2\times M^2n^2}$ and $\boldsymbol{\Xi}\in\mathbb{R}^{Mn^2\times Mn^2}$, it is computationally advantageous (by an order of magnitude in the discretization parameter $M$) to replace the condition $\det\,(\mathbf{F})=0$ with the condition $\det\,(\boldsymbol{\Xi})=0$ if one wishes to rely on discretization for mapping out second-moment stability boundaries.

Given $\mathbf{m}(t)$ and an initial history function for $t\in[-\tau,0]$ that is consistent with the definition of $\phi_t(\vartheta)$ in terms of $\mathbb{E}[x_t(0)x_t(\vartheta)]$, the surrogate correlation DDE \eqref{eq:corr_disc_eq} may be integrated forward in time using a standard DDE solver. For $\mathbf{m}\equiv 0$, the stability of the stationary solution may be investigated either by i) requiring that the roots of a corresponding high-dimensional, matrix-valued characteristic quasipolynomial all have negative real parts, or by ii) first discretizing the advection form of \eqref{eq:corr_disc_eq} and then requiring that all the eigenvalues of the coefficient matrix of the resultant high-dimensional system of ODEs have negative real part. In both cases, the initial advantage of having reduced the system dimension by a factor of $M$ is lost and we arrive at the same level of complexity as when analyzing the spectrum of $\mathbf{F}$.
\end{document}